\documentclass[11pt]{article}
\usepackage[square,sort,comma,numbers]{natbib}
\usepackage{amsthm,amssymb,amsmath}
\usepackage[pdftex]{thumbpdf}       
\usepackage{float}
\usepackage{graphicx}
\usepackage{verbatim}
\usepackage[usenames,dvipsnames]{color}
\usepackage{dsfont}
\usepackage{fullpage}
\usepackage{ifpdf}
\ifpdf
\usepackage[pdftex]{hyperref} 
\usepackage{authblk} 
\usepackage[utf8]{inputenc}
\usepackage[T1]{fontenc}
\usepackage{blindtext}
\usepackage{booktabs}
\usepackage{tabularx}
\usepackage[dvipsnames]{xcolor}

\usepackage{natbib}
 \bibpunct[, ]{(}{)}{,}{a}{}{,}%

\hypersetup{%
pdftitle={Matching policies for two-sided matching},
pdfauthor={A. Aveklouris},
pdfkeywords={matching},
bookmarks=true,%
bookmarksnumbered=true,
pdfstartview={FitH},
pdfpagelayout={OneColumn},
colorlinks=true,
citecolor=NavyBlue,
linkcolor=RedOrange,
urlcolor=BrickRed,
bookmarksopen=true,%
bookmarksopenlevel=1,
unicode=true,%
breaklinks=true,%
}%
\else
  \newcommand\phantomsection\relax
  \newcommand{\url}[1]{#1}
  \newcommand{\href}[2]{#2}
  \usepackage{nohyperref}
\fi

\usepackage{url}
\mathchardef\ordinarycolon\mathcode`\: \mathcode`\:=\string"8000 \begingroup
\catcode`\:=\active \gdef:{\mathrel{\mathop\ordinarycolon}} \endgroup
 \theoremstyle{plain}              

\usepackage{caption}
\usepackage{subcaption}
\usepackage{algorithm,algpseudocode}
\usepackage[hang,flushmargin]{footmisc}

\usepackage{tikz}
\usetikzlibrary{calc}

\newtheorem{proposition}{Proposition}
\newtheorem{lemma}{Lemma}
\newtheorem{theorem}{Theorem}
\newtheorem{corollary}{Corollary}
\newtheorem{remark}{Remark}
\newtheorem{definition}{Definition}

\newtheorem{claim}{Claim}

\newtheorem{examp}{Example}

\newcommand{\bbJ}{{\mathbb J}}
\newcommand{\bbK}{{\mathbb K}}
\newcommand{\calM}{{\mathcal M}}

\newcommand{\pfa}[1]{{\noindent\emph{\textbf{#1.\ }}}}
\newcommand{\eof}{\hfill\qedsymbol}

\newcommand{\p}{{v}}
\newcommand{\ad}{{A}^D}
\newcommand{\as}{{A}^{S}}
\newcommand{\asn}{{A}^{S,n}}
\newcommand{\qd}{{Q}^D}
\newcommand{\qdn}{{Q}^{D,n}}

\newcommand{\bbS}{\mathbb S}

\newcommand{\amy}[1]{{\color{black}#1}}
\newcommand{\levi}[1]{{\color{black}#1}}

\newcommand{\E}[1]{\mathbb{E}\left[#1\right]}  
\newcommand{\Prob}[1]{\mathbb{P}\left(#1\right)} 
\newcommand{\ind}[1]{1{\left\{#1\right\}}}   
\newcommand{\bld}{\boldsymbol}
\newcommand{\R}{\mathbb{R}}
\newcommand{\Rp}{\R_+}
\newcommand{\mc}{\mathcal}
\newcommand{\ovl}{\overline}
\newcommand{\lran}[1]{\left\langle#1\right \rangle}

\newcommand{\calP}{\mathcal{P}}
\newcommand{\calQ}{\mathcal{Q}}
\newcommand{\mstar}{m^\star}
\graphicspath{{Figures/}}
\newcommand{\arcset}{\mathcal{E}}
\newcommand{\demset}{\mathcal{D}}
\newcommand{\suppset}{\mathcal{S}}

\begin{document}
\title{Matching Impatient and Heterogeneous Demand~and~Supply}
\author{Angelos\ Aveklouris}
\author{Levi\ DeValve}
\author{Maximiliano Stock}
\author{Amy\ R. Ward}
\affil{The University of Chicago Booth School of Business}
\date{}
\maketitle

 \begin{abstract}
Service platforms must determine rules for matching heterogeneous demand (customers) and supply (workers) that arrive randomly over time and may be lost if forced to wait too long for a match.  Our objective is to maximize the cumulative value of matches, minus costs incurred when demand and supply wait.  We develop a fluid model, that approximates the evolution of the stochastic model, and captures explicitly the nonlinear dependence between the amount of demand and supply waiting and the distribution of their patience times, \amy{also known as reneging or abandonment times in the literature}.  The fluid model invariant states approximate the steady-state mean queue-lengths in the stochastic system, and, therefore, can be used to develop an optimization problem whose optimal solution provides matching rates between demand and supply types that are asymptotically optimal (on fluid scale, as demand and supply rates grow large).  We propose a discrete review matching policy that asymptotically achieves the optimal matching rates.  We further show that when the aforementioned matching optimization problem has an optimal extreme point solution, which occurs when the patience time distributions have increasing hazard rate functions, a state-independent priority policy, that ranks the edges on the bipartite graph connecting demand and supply, is asymptotically optimal. A key insight from this analysis is that the ranking critically depends on the patience time distributions, and may be different for different distributions even if they have the same mean, demonstrating that models assuming, e.g., exponential patience times for tractability, may lack robustness.  Finally, we observe that when holding costs are zero, a discrete review policy, that does not require knowledge of inter-arrival and patience time distributions, is asymptotically optimal.

\end{abstract}

\textbf{Keywords:} matching; two-sided platforms; bipartite graph; discrete review policy; high-volume setting; fluid model; \amy{impatience}; reneging; \amy{abandonment} \\

\section{Introduction}\label{sec:Introduction}

Service platforms exist to match demand (customers) and supply (workers); see, for example, \cite{Hu2019} for a broader perspective on the intermediary role of platforms in the sharing economy.  The challenge in operating these platforms is that demand and supply often come from heterogeneous customers and workers that arrive randomly over time, and do not like to wait.  Then, there is a trade-off between making valuable matches quickly, and waiting in order to enable better matches.  Our goal is to develop matching policies that optimize a long-run objective function consisting of the total value of matches made minus holding costs incurred when demand and supply wait.

The dislike of waiting can lead to demand and supply abandoning before being matched, e.g., as occurs in a ride hailing application \amy{(\cite{wang2019demand})}, in which customers must be matched to \amy{workers} (drivers).  Moreover, the willingness of \amy{agents\footnote{\levi{We henceforth use the term \emph{agent} when collectively referring to customers and workers together.}}} to wait can change over time, with some \amy{agents} becoming more patient (for example, as a result of having already paid a waiting cost), and some becoming less patient (for example, as a result of becoming more and more frustrated while waiting).  In a transplant application setting, where donors must be matched with needy patients, the survival time of a patient depends on the history, and not only on the current state of health (\cite{bertsimas2019development}). \amy{Also}, in a call center setting, empirical evidence does not usually support using an exponential distribution to fit customer patience times (\cite{mandelbaum2013data, brown2005statistical}).
As a result, using an exponential distribution to model \amy{agent} patience times may be too restrictive.
Still, the exponential distribution is more analytically tractable, and, if optimal or near-optimal matching policies are robust to distributional characteristics, then the exponential distribution is a good modeling choice.

\begin{figure}[t]
\centering
\begin{subfigure}[t]{.3\textwidth}
  \centering
  \includegraphics[width=\linewidth]{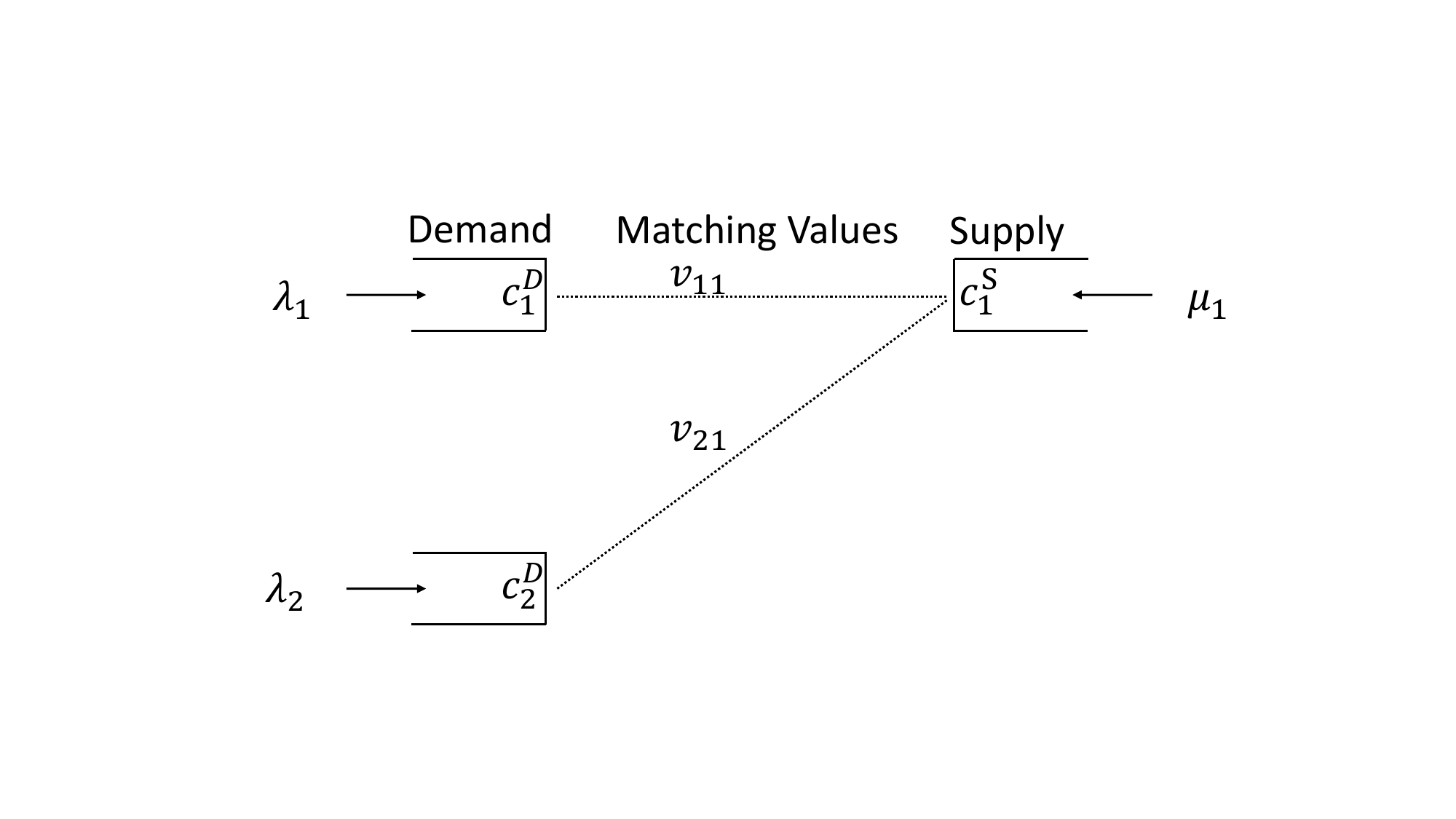}
  \caption{The network.}
  \label{figure:badmodel}
\end{subfigure}
\begin{subfigure}[t]{.3\textwidth}
  \centering
  \includegraphics[width=.8\linewidth]{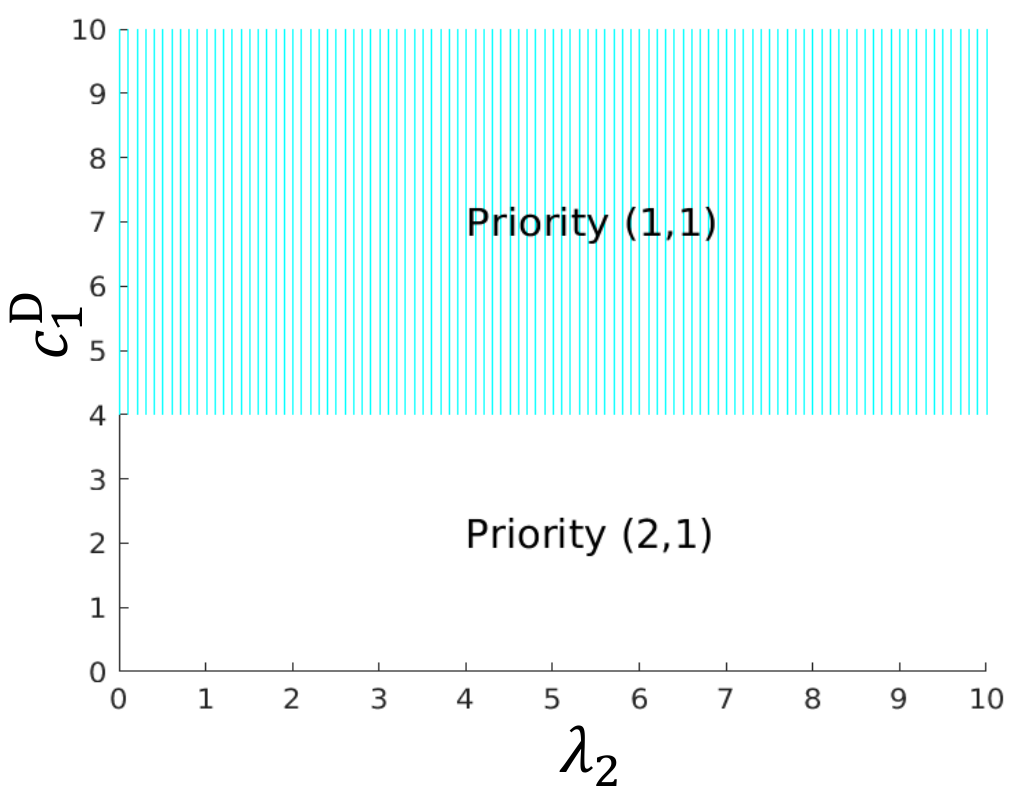}
  \caption{Exponential patience times.}
  \label{figure:badexample}
\end{subfigure}
\begin{subfigure}[t]{.3\textwidth}
  \centering
  \includegraphics[width=.8\linewidth]{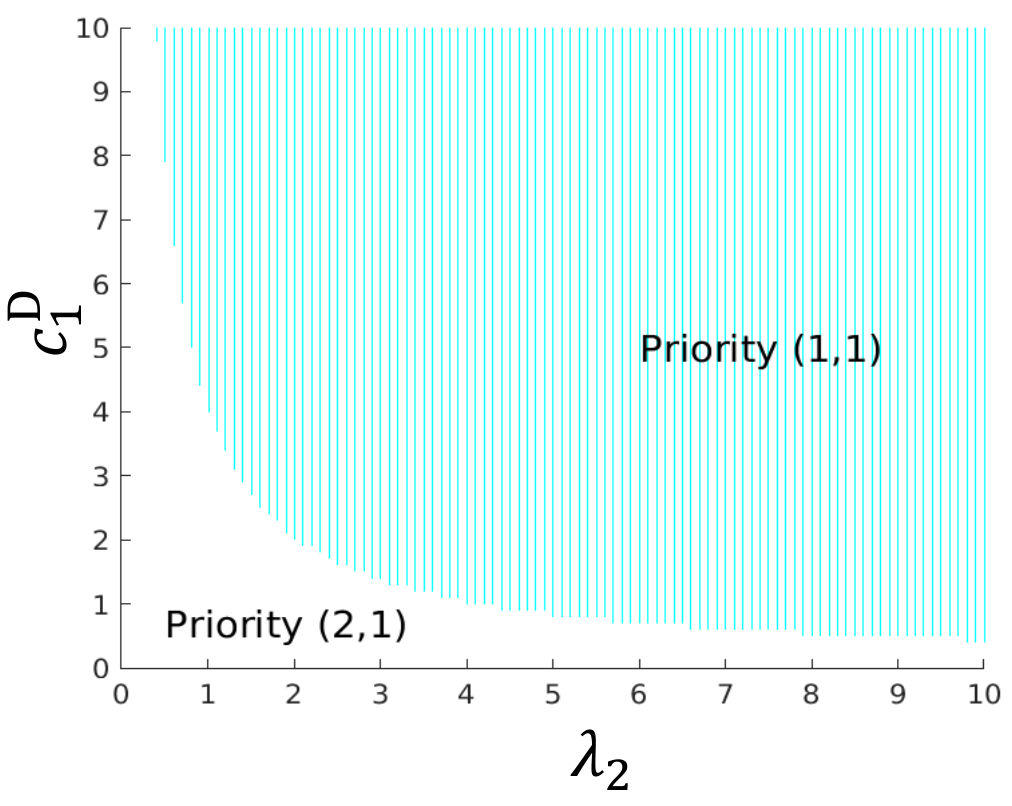}
  \caption {Uniform patience times.}
  \label{figure:badexample1}
\end{subfigure}
\caption{\amy{The priority matches given by the MP solution} for a network
with one supply node and two demand nodes, having matching values  $\p_{11}=\p_{21}=1$ (i.e., matching a supply unit with a demand unit from either node gives a value of one), identical patience time distributions for all types having mean equal to one, holding costs per job per unit time $c_2^D=4$, $c_1^S=1$, and Poisson arrivals with rates per unit time $\lambda^D_1=\lambda^S_1=1$.}
\label{fig:badexample}
\end{figure}

\amy{The issue when testing whether or not the aforementioned robustness holds is that solving for an optimal matching policy appears intractable for non-exponential patience time distributions.  This is because without the memoryless property, in order for the system to be Markovian, the state space must track the amount of time each \levi{agent} has waited, resulting in a very complex state space.  To overcome this, we study a deterministic fluid model, that serves as an approximation to the stochastic model when arrival rates are large.  The invariant states of the fluid model (that is, the fixed points of the fluid model equations) approximate the steady-state mean amount of demand and supply waiting to be matched, and depend nonlinearly on the patience time distributions.  We use the aforementioned invariant states to define an optimization problem, termed the \emph{matching problem} (MP),
whose objective aims to maximize the value of matches made minus holding costs.  Then we study the effect of the patience time distribution on the MP solution.}

Figure~\ref{fig:badexample} shows that the choice of distribution used to model \levi{agent} patience times can change \levi{how} the \amy{ MP solution} \levi{prioritizes matches}.   For the simple matching network in
Figure~\ref{figure:badmodel}, the \amy{MP solution} \levi{either prioritizes making matches on edge $(1,1)$ before using any leftover supply}
for making matches on the edge (2,1), or vice versa.  Figures~\ref{figure:badexample} and \ref{figure:badexample1} show that the same parameters lead to different priority orderings for the exponential and uniform patience time distributions.  As a result, we are motivated to develop methodology that allows us to analyze as general a class of patience time distributions as possible; in particular, the assumptions we require \amy{to justify the MP are that the patience time distributions have a density, a finite mean, and a strictly increasing cdf\footnote{\amy{The proof that the solution to the MP motivates an asymptotically optimal matching policy in a high-volume regime (specifically, Theorems~\ref{proposition:DRconvergence2} and~\ref{Thm:aopriority} in Section~\ref{sec:Performance Analysis}) requires the additional assumption that the hazard rates associated with the patience time distribution functions are bounded. However, our simulation results in Section~\ref{sec:SimJK} suggest that our proposed priority ordering policies perform well even when hazard rates are unbounded (as is true for the uniform distribution). See also Remark~\ref{rem:bounded-hazard}.}}.}

The matching model we consider generalizes that shown in Figure~\ref{fig:badexample} to allow for an arbitrary number of demand and supply types, each with a possibly different patience
time distribution and arrival rate.  The matching values between demand and supply types are the edge values on the underlying bipartite graph, which could be positive or zero.  Since different edges have different values, some matches are preferable to others.  Since different demand and supply types have different holding costs and patience time distributions, waiting is more costly for some types than others. Our objective captures the trade-off between these opposing priorities by incorporating both matching values and holding costs.



\subsection{Contributions of this paper}
We contribute a number of theoretical results and practical insights for matching on service platforms. First, the MP provides an asymptotic upper bound on the objective function value as arrival rates become large (Theorem~\ref{prop:upper}).  We then develop a discrete review matching policy whose matching rates mimic an MP optimal solution (Theorem~\ref{proposition:DRconvergence2}), and  is asymptotically optimal as arrival rates become large. As a consequence, because an MP optimal solution depends on the patience \levi{time} distributions, so does the proposed matching policy. While this policy is applicable for general patience time distributions, our analysis provides further insight for the important special case when there exists an optimal extreme point solution to the MP (this is ensured when the patience time hazard rates are increasing, i.e., supply and demand become more impatient the longer they wait), which we discuss next.

When the patience time hazard rate functions are increasing, the MP becomes a convex maximization problem, and so has an optimal extreme point solution (Theorem~\ref{lem:concave}).  In this case, we propose a discrete review policy that is based on a static ranking of the edges, and which makes matches by prioritizing the use of demand and supply according to that ranking. We show that this policy also achieves the matching rates given by the MP (Theorem~\ref{Thm:aopriority}) and is asymptotically optimal as arrival rates become large. The static ranking depends on the patience time distribution, and can be different for different distributions with the same mean, as observed earlier in Figure~\ref{fig:badexample}.
\amy{As a result, the service platform must obtain distributional information to determine edge priorities -- knowing only information about the mean is not enough.}


In the case that holding costs are zero, we observe that an MP optimal solution does not depend on the patience time distribution.  Then, a discrete review policy based on a linear problem (LP), that does not need to know information about the arrival rates or the patience time distribution, is asymptotically optimal as arrival rates become large (Theorem~\ref{theorem:blindOptimal}).

The MP is based on the invariant states of a fluid model (Definition~\ref{Definition:IM}), as mentioned earlier.
\amy{Our asymptotic optimality results mentioned in the previous paragraphs require showing that }
 the fluid model starting from any initial state converges to an invariant state as time becomes large (Theorem~\ref{Th:steadySMB}).
Our convergence proof is non-trivial when viewed in light of the fact that such a result for the fluid model approximating the single class many-server queue with \amy{impatience}  was only shown recently in \cite{atar2021large} (despite being worked on for many years), because dealing with the measure tracking the age-in-service is complicated. In our analysis, we observe that the aforementioned difficulty is \amy{not present} in matching systems,
\amy{because when service is instantaneous the need for an age-in-service measure is eliminated (an observation also leveraged in \cite{DCN2018}).  }

The remainder of the paper is organized as follows. We end Section~\ref{sec:Introduction} by reviewing related literature. In Section~\ref{sec:Model description}, we provide a detailed model description. In Section~\ref{sec:DRMUP}, we propose our matching policies. A fluid model is presented in Section~\ref{sec:fluid model}. In particular, we present asymptotic approximations for the stationary mean queue-lengths, and we show the impact of the patience time distributions. In Section~\ref{sec:Performance Analysis}, we introduce our high-volume setting and we study the asymptotic behavior of our proposed matching policies. \levi{In Section~\ref{sec:SimJK} we demonstrate numerically that our policies perform well in both high volume and non-asymptotic simulations.} In Section~\ref{sec:PSB}, we study a platform without holding costs and we propose a matching policy that does not require the knowledge of the system parameters. Concluding remarks are \amy{made} in Section~\ref{sec:conclusion}. All the proofs are given in \amy{the Appendix.}

\subsection{Literature review}
We focus on on-demand service platforms that aim to facilitate matching. We refer the reader to \cite{BenjaafarHu2020}, \cite{ChenEtAl2020}, and \cite{Hu2020}  for  excellent higher-level perspectives on how such platforms fit into the sharing economy and to~\cite{Hu2019} for a survey of recent sharing-economy research in operations management.  Three important research questions for such platforms identified in \cite{BenjaafarHu2020} are how to price services, pay workers, and match requests.  These decisions are ideally made jointly; however, because the joint problem is difficult, the questions are often attacked separately.  In this paper, we focus on the matching question.

Our basic matching model is a bipartite graph with demand on one side and supply on the other side.  There is a long history of studying two-sided matching problems described by bipartite graphs, beginning with the stable matching problem introduced in the groundbreaking work of~\cite{GaleShapley1962} and continuing to this day; see~\cite{AS2013} and \cite{RS1990} for later surveys and~\cite{ABKS2020} for more recent work.  In the aforementioned literature, much attention is paid to eliciting agent preferences because the outcomes of the matching decisions, made at one prearranged point in time, can be life-changing events (for example, the matches between medical schools and potential residents).  In contrast, many platform matching applications are {\em not} life-changing events,  so there is less need to focus on eliciting agent preference.  Moreover,  supply and demand often arrive randomly and continuously over time,  must be matched dynamically over time, as the arrivals occur (as opposed to at one prearranged point in time), \amy{and will leave the system without being matched (renege) if made to wait too long}.

\amy{
Recent work has  used dynamic two-sided matching models with reneging in the context of organ allocation (\cite{Zenios1999, BoxEtAl2011, DCN2018, KL2020}), trading systems known as crossing networks  (\cite{AfecheEtAl2014}), online dating and labor markets (\cite{AJK2020,KS2020}),  ridesharing (\cite{BanerjeeKanoriaQian2020, OzkanWard2020, ozkan2020joint}), and quantum switches (\cite{ZubEtAl2023}).  Simultaneously and motivated by the aforementioned applications, there have been many works that begin with a more loosely motivated modeling abstraction, such as \cite{buke2017fluid, liu2019}, \cite{blanchet2020asymptotically, CastroEtAl2020, JonkEtAl2023, MM2023, KG2023}.
The aforementioned works focus on a range of issues; however, except for \cite{BoxEtAl2011, DCN2018, KL2020} all assume that agents' waiting times are either deterministic or exponentially distributed. Allowing for more general willingness-to-wait distributions is important because in practice the time an agent is willing to wait to be matched may depend on how long that agent has already waited.  The difficulty is that tracking how the system evolves over time requires tracking the remaining time each agent present in the system will continue to wait to be matched, resulting in an infinite-dimensional state space.  In contrast to our work, \cite{BoxEtAl2011,KL2020} has one demand and one supply type and \cite{DCN2018} considers an index policy class.
}

Taking an approach similar in spirit to how stochastic processing networks are controlled in the queueing literature developed in \cite{Har:95b, harrison2006correction,maglaras1999dynamic, maglaras2000discrete}, we use a high-volume asymptotic regime to prove that our matching policies are asymptotically optimal.
In particular, we begin by solving a static matching problem and we focus on fluid-scale asymptotic optimality results. Adopting an approach reminiscent of the one \cite{GW14} use in a dynamic multipartite matching model that assumes agents will wait forever to be matched, we use a discrete review policy to balance the trade-off between having agents wait long enough to build up matching flexibility but not so long that many will leave.


We consider an objective function that involves the queue-lengths and is an infinite horizon objective. As a result, we need to consider the steady-state performance of the system. This is the focus of the infinite bipartite matching queueing models considered for example  in
\cite{adan2012exact, AW2014, ABMW2018, ACG2019, BGM2013, CKW2009, FK2018, diamant2019double}.

\section{Model description}
\label{sec:Model description}

\amy{We use the following notation conventions.  All vectors and matrices are denoted by bold letters.
Further, $\R$ is the set of real numbers, $\R_+$ is the set of non-negative real numbers, $\mathbb{N}$ is the set of strictly positive integers, and $\mathbb{Z}_+=\mathbb{N}\cup \{0\}$.} \\

\begin{figure}[htb]
\centering
  {\includegraphics[width=0.7\textwidth]{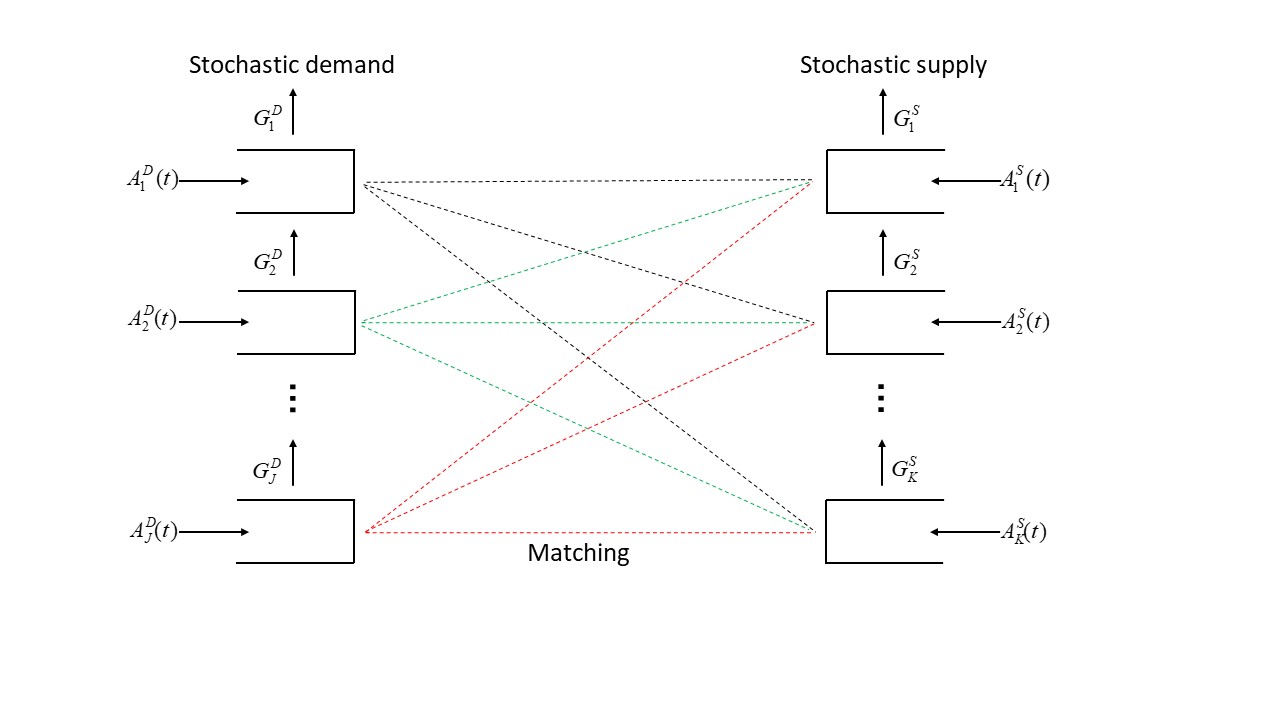}}
\caption{The two-sided matching model with  general patience time distributions.}
\label{figure:Matching}
\end{figure}

\textbf{Primitive inputs:}
The set of demand nodes
$\bbJ :=\{ 1,\ldots,J\}$ (representing customer types)
and  supply nodes (representing worker types)
$\bbK := \{ 1,\ldots, K\}$ form a bipartite graph, as shown in Figure~\ref{figure:Matching}.
The set $\arcset \subseteq \bbJ \times \bbK$ denotes the set of compatible matches between demand and supply nodes, i.e., demand type $j\in \bbJ$ can be matched with supply type $k \in \bbK$ if and only if $(j,k) \in \arcset$. Further, for demand type $j\in \bbJ$ let $\suppset_j = \{k|(j,k) \in \arcset\}$ denote the set of supply types compatible with $j$, and similarly for supply type $k\in \bbK$ let $\demset_k = \{j|(j,k) \in \arcset\}$ denote the set of demand types compatible with $k$.
The value of matching demand type $j\in \bbJ$ and supply type $k \in \bbK$ is $\p_{jk}\ge 0$.
The holding costs $c_j^D\ge 0$ and $c_k^S\ge 0$ are incurred for each unit of time demand type $j\in \bbJ$ or supply type $k \in \bbK$ waits.
\levi{Agents} arrive according to renewal processes, denoted by $\ad_j(\cdot)$,  $j \in \bbJ $  and $\as_k(\cdot)$, $k \in \bbK$, having respective rates $\lambda_j^D> 0$ and $\lambda_k^S> 0$, \amy{and inter-arrival distributions with finite 5th moment (an assumption convenient for our asymptotic analysis, but not strictly required)}.  Upon arrival, each
type $j\in \bbJ$ customer and type $k\in \bbK$ worker independently samples their patience time from respective distributions $G_j^D(\cdot)$ and $G_k^S(\cdot)$ having support on $[0,H^D_j)$ and $[0,H^S_k)$ for some $H^D_j \in [0, \infty] $ and $H^S_k \in  [0, \infty]$, which represents the maximum amount of time they are willing to wait to be matched. \levi{Agents} that wait longer than their patience time leave the system without being matched.
The patience times are absolutely continuous random variables with density functions $g_j^D(\cdot)$ and $g_k^S(\cdot)$ \amy{that are mutually independent of each other and of the arrival processes $\bld{A}^D$ and $\bld{A}^S$.}
We let $h_j^D(x) = g_j^D(x) / (1-G_j^D(x))$ for $x \in [0,H_j^D)$ and $h_k^S(x) = g_k^S(x) / (1-G_k^S(x))$ for $x \in [0,H_k^S)$ be the associated hazard functions.

\amy{
Throughout the paper, we assume that the patience time distributions satisfy the following conditions:
\begin{enumerate}
\item[(1)] $\int_0^\infty (1-G_j^D(x))dx=1/\theta_j^D\in(0,\infty)$, $j\in\bbJ$, and $\int_0^\infty (1-G_k^S(x))dx=1/\theta_k^S\in(0,\infty)$, $k\in\bbK$;
\item[(2)] $G_j^D$, $j \in \bbJ$, and $G_k^S$, $k \in \bbK$ are strictly increasing.
\end{enumerate}
\noindent The condition (1) above ensures that the patience time distributions have finite positive mean and (2) ensures that the inverse functions
$\left( G_j^D \right)^{-1}: [0,1) \rightarrow [0,H_j^D)$ and $\left( G_k^S \right)^{-1}: [0,1) \rightarrow [0,H_k^S)$
are well-defined for each $j \in \bbJ$ and $k \in \bbK$.
The excess life distributions of the patience time distributions can be defined under (1) above and are as follows:
\begin{equation*}
G^D_{e,j}(x)=\int_{0}^{x} \theta_j^D (1-G^D_j(u))du, \mbox{ for }
 j\in \bbJ \mbox{ and } x\in \R_+
\end{equation*}
and
\begin{equation*}
G^S_{e,k}(x)=\int_{0}^{x} \theta_k^S (1-G^S_k(u))du, \mbox{ for }
 k\in \bbK \mbox{ and } x\in \R_+.
\end{equation*}

The patience time distributions are also known as reneging distributions in the literature.  Then, the cumulative number of \levi{agents} that leave without being matched is more succinctly termed the cumulative reneging.  Throughout this paper, we fluctuate between the terms patience and reneging.
}
\\


\noindent \textbf{Objective and admissible policies:}  The objective is to maximize the long-run average value of matches made, minus holding costs incurred for waiting.  A matching \amy{process}  is a $|\arcset|$-dimensional stochastic process $\bld{M}(\cdot) := \{ M_{jk}(\cdot),  (j,k) \in \arcset \}$ that  specifies the cumulative number of matches made  between types $j\in \bbJ$ and $k\in \bbK$ in the time interval $[0,t]$, denoted by $M_{jk}(t)$, for each $t \geq 0$, under the assumption that matches are made first-come-first-served (FCFS) within each type.  Then, if $\qd_j(t)$ represents the number of type $j\in \bbJ$ customers waiting and $Q^S_k(t)$ represents the number of type $k \in \bbK$ workers waiting at time $t \geq 0$, we want to maximize
\begin{equation} \label{eq:Profit}
 V_{\bld{M}}:= \liminf_{t\rightarrow \infty}
 \frac{1}{t}V_{\bld{M}}(t),
\end{equation}
where $$V_{\bld{M}}(t):=\sum_{(j,k) \in \arcset}
 \p_{jk} M_{jk}(t)-
 \sum_{j\in \bbJ} \int_{0}^{t}c_j^D Q^D_{j}(s)ds -
 \sum_{k\in \bbK} \int_{0}^{t} c_k^S Q^S_{k}(s)ds.$$

The objective (\ref{eq:Profit}) captures the trade-off between making matches quickly and waiting for better matches.    The difficulty when deciding whether or not to incur holding costs for a short period of time in order to enable a potential future higher value match is that the matching policy does not know how much longer each \levi{agent} currently waiting will remain waiting without being matched. \\

\begin{figure}[htb]
\centering
  {\includegraphics[width=0.8\textwidth]{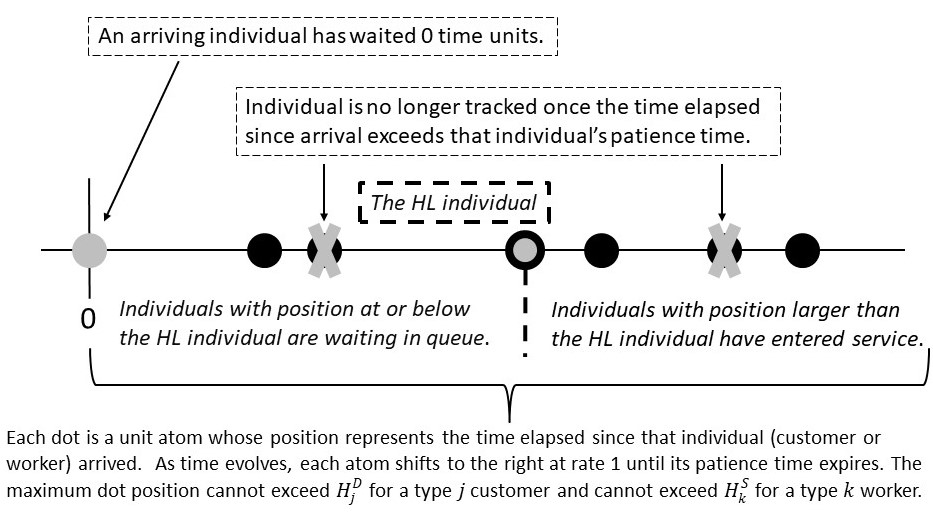}}
\caption{\amy{The measure-valued process $\eta_j$ (for given customer type) or $\eta_k$ (for given worker type).}}
\label{figure:eta}
\end{figure}

\amy{
\noindent \textbf{System state:}
The system state in our model is complex, because the state must track how long each \levi{agent} in queue has been waiting, as well as the time passed since the last arrival of each \levi{agent} type, in order to be Markovian.  In particular, tracking how long each \levi{agent} in the queue has been waiting requires the use of a measure.  For $H \in \R_+ \cup \{\infty\}$, let $\mc{M}[0,H)$ be the set of finite nonnegative Borel measures on $[0,H)$, endowed with the topology of weak convergence.  A state at time $t \in \R_+$ in our model is described as follows:  for each $j \in \bbJ$ and $k \in \bbK$,
\begin{itemize}
\item $\alpha_j^D(t) \in \R_+$ and $\alpha_k^S(t) \in \R_+$ are the times elapsed since the last type $j$ customer arrived and last type $k$ worker arrived;
\item $Q_j^D(t) \in \mathbb{Z}_+$ and $Q_k^S(t) \in \mathbb{Z}_+$ are the number of type $j$ customers and type $k$ workers waiting in queue;
\item $\eta_j^D(t) \in  \mc{M}[0,H_j^D)$ and $\eta_k^S(t) \in  \mc{M}[0,H_k^S)$ store the amount of time that has passed between each type $j$ customer's arrival time, and each type $k$ worker's arrival time, up until that customer's or worker's sampled patience time.
\end{itemize}

The measures $\eta_j^D$ and $\eta_k^S$  track the evolution of unit atoms over time, where each atom is associated with a particular type $j$ customer's or type $k$ worker's time-since-arrival, as shown in Figure~\ref{figure:eta}.  All \levi{agents} tracked are ``potentially'' waiting in queue, because their wait is less than their sampled patience time.  The term ``potential'' refers to the fact that such \levi{agents} may or may not have been matched.  The FCFS matching assumption implies that all potential customers (workers) that have waited longer than the customer (worker) at the head-of-the-line for their type have already been matched, and all potential customers (workers) that have waited less than that customer are in queue.  As such, at any time $t \in \R_+$, for each $j \in \bbJ$ and $k \in \bbK$,
\begin{equation} \label{eq:state-space-restrictions}
    Q_j^D(t) \leq \int_0^{H_j^D} \eta_j^D(t) dt \mbox{ and }  Q_k^S(t) \leq \int_0^{H_k^S} \eta_k^S(t) dt.
\end{equation}
 The measures $\eta_j^D(t)$ and $\eta_k^S(t)$ are independent of the matching process policy, which facilitates the analysis.
}

\noindent
\textbf{Queue evolution equations:}
 For a given matching \amy{process}, $\bld{M}(\cdot)$, the demand and supply queue-lengths at time $t\ge 0$ are given by
\begin{equation} \label{eq:posQ}
  \qd_j(t) :=   \qd_j(0) +  \ad_j(t) - R_j^D(t)-\sum_{k \in \suppset_j} M_{jk}(t)\ge 0
\end{equation}
and
 \begin{equation} \label{eq:posI}
 Q^S_k(t) :=  Q^S_k(0) + \as_k(t) - R_k^S(t)- \sum_{j\in \demset_k}  M_{jk}(t)\ge 0,
 \end{equation}
for all $j\in \bbJ$ and $k\in \bbK$, where $R_j^D(t)$ and $R_k^S(t)$ denote the cumulative number of type $j\in \bbJ$ customers and type  $k\in \bbK$ workers that left the system without being matched in $[0,t]$.   The equations \eqref{eq:posQ} and \eqref{eq:posI} are balance equations that say that the number of type $j \in \bbJ$ customers and type $k \in \bbK$ workers waiting to be matched at time $t\ge 0$ equals those present initially plus the cumulative arrivals minus the cumulative reneging minus the cumulative matches.

We refer to $\bld{R}^D(\cdot)$ and $\bld{R}^S(\cdot)$ as the reneging processes, because in the queueing literature leaving the system before being served \amy{(matched in our context)} is commonly known as reneging.  \amy{We do not explicitly construct these processes, nor do we specify the system dynamics, because that involves substantial mathematical overhead that is not relevant to our goal of maximizing (\ref{eq:Profit}).  Instead we refer the reader to the companion technical paper (\cite{aveklouris2021MM}) for that.  To advance our goal of maximizing (\ref{eq:Profit}), we first accept that finding an exact solution does not appear possible, and second focus on studying a deterministic matching problem (see \eqref{eq:max_A_w_v3G} below)  whose objective function approximates (\ref{eq:Profit}).} \\

\textbf{Admissible matching policy:}
\amy{
A matching policy induces a matching process $\bld M$, and that may restrict the set of achievable system states.  For example, if the matching policy prioritizes matches between type $j \in \bbJ$ customers and type $k \in \bbK$ workers, then the matching policy will disallow states in which both type $j$ customers and type $k$ workers are present, meaning that $Q_j^D(t) Q_k^S(t) = 0$ for all $t \geq 0$.

A matching policy is admissible if the induced matching process $\bld M$ satisfies the following natural assumptions.
 First, no partial matches can be made, and, if a match occurs, then it cannot be taken back.  That is, a matching process is integer-valued and non-decreasing.  Next, matches are recorded exactly when made (not before), meaning a matching process is right-continuous with left limits.
Third, \levi{agents} must be present in the system to be matched, which is equivalent to having non-negative \levi{agent} queue-lengths, as required in (\ref{eq:posQ}) and (\ref{eq:posI}).   Fourth, a matching process is non-anticipating, meaning the process cannot know the exact times of future arrivals.
We call matching processes that satisfy all the aforementioned properties \textbf{\it admissible}, and the maximization in (\ref{eq:Profit}) is over the class of admissible matching processes\footnote{For a precise mathematical specification of a matching policy and the admissible policy class, see Definitions 1 and 2 in (\cite{aveklouris2021MM}).  A more rigorous mathematical specification is not possible in this paper because we do not fully specify the system dynamics.}.
}

\section{Proposed matching policies}
\label{sec:DRMUP}
In this section, we develop matching policies with the goal of \amy{asymptotically} maximizing our objective \eqref{eq:Profit}. To do so, we first define an optimization problem that determines the optimal matching rates (Section~\ref{sec:SMP}) and we refer to it as the \emph{matching problem} (MP). Then, we propose a discrete review matching policy that asymptotically achieves these rates (Section~\ref{sec:matching rate}), and we call this policy the \emph{matching-rate-based policy}. When the existence of an optimal extreme point solution \amy{to the MP} is ensured, we are able to propose a discrete review \emph{priority-ordering policy}  that prioritize the edges and also asymptotically achieves the optimal matching rates (Section~\ref{sec:priotiry}).

A discrete review policy decides on matches at review time points and does nothing at all other times  \citep[see, for example,][]{GW14}. Longer review periods allow more flexibility in making matches but risk losing impatient \levi{agents}.  Shorter review periods prevent customer/worker loss \amy{and \levi{reduce} holding cost} but may not have sufficient numbers of customers/workers of each type to ensure the most valuable matches can be made.

Let $t>0$. We let $l, 2l, 3l,\ldots$ be the discrete review time points, where $l \in (0,t)$ is the review period length.
\amy{Since a discrete review policy does not make matches between review time points, at time $il$ for $i \in \{1,\ldots, \lfloor t/l \rfloor \}$,
from \eqref{eq:posQ},  the number of type $j$ customers  available to be matched
is}
\begin{equation} \label{eq:queueDRdemand}
Q_j^D(il-)=\qd_j((i-1)l)+ \ad_j(il) - \ad_j((i-1)l)- R^D_j(il) + R^D_j((i-1)l),
\end{equation}
\amy{for each $j\in \bbJ$,  and, from \eqref{eq:posI},  the number of  type $k$ workers available to be matched is}
 \begin{equation} \label{eq:queueDRsupply}
 Q^S_k(il-)=Q^S_k((i-1)l) + \as_k( i l )- \as_k((i-1)l)- R^S_k(il) + R^S_k((i-1)l),
 \end{equation}
\amy{for each $k\in \bbK$.}

\subsection{A matching problem}\label{sec:SMP}
Suppose we ignore the discrete and stochastic nature of \levi{agent} arrivals, and assume that these flow at their long-run average rates.
Then, an upper bound on the long-run average matching value in (\ref{eq:Profit}) follows by solving an optimization problem to find the optimal instantaneous matching values, which we refer to as the matching problem (MP).  For $\bld m = (m_{jk}: (j,k) \in \arcset)$, that denotes the instantaneous matching rate, the MP can be defined through functions $q_j^{D,\star}( \bld m)$, $j\in \bbJ$ and $q_k^{S,\star}( \bld m)$, $k\in \bbK$ (defined in detail below), that approximate the steady-state mean queue-lengths, as follows:
\begin{align} \label{eq:max_A_w_v3G}
\begin{aligned}
\max \ &
\sum_{(j,k) \in \arcset}
 \p_{jk} m_{jk}-
 \sum_{j \in \bbJ} c_{j}^D q_j^{D,\star}( \bld m)-
\sum_{k \in \bbK} c_{k}^S q_k^{S,\star} (\bld m)  \\
\text{s.t.} \ & \sum_{j\in \demset_k} m_{jk} \le \lambda^S_{k}, \ k \in \bbK,\\
&\sum_{ k \in \suppset_j} m_{jk} \le \lambda_j^D, \ j \in \bbJ, \\
& m_{jk} \ge 0,  \ (j,k) \in \arcset.
\end{aligned}
\end{align}
The objective function in (\ref{eq:max_A_w_v3G}) is the instantaneous version of the long-run average matching value in (\ref{eq:Profit}) and the constraints in (\ref{eq:max_A_w_v3G}) prevent us from matching more demand or supply than is available.

An optimal solution to the MP (\ref{eq:max_A_w_v3G}), $\bld m^\star = (m_{jk}^\star: (j,k) \in \arcset)$, can be interpreted as the optimal instantaneous rate of matches between demand $j \in \bbJ$ and supply $k \in \bbK$.  An admissible matching policy  should maximize the long-run average matching value in (\ref{eq:Profit}) if its associated matching rates equal $\bld m^\star$.

Unfortunately, having closed form expressions for the exact steady-state mean queue-lengths appears intractable.  Fortunately, when demand and supply arrival rates become large, we can approximate the steady-state mean queue-lengths by developing a fluid model to approximate the evolution of the stochastic model defined in Section~\ref{sec:Model description}, and finding its invariant states, which we do in Section~\ref{sec:fluid model}
(see \eqref{eq:InvDeQ} and \eqref{eq:InvSuI} in Proposition~\ref{prop:Chofinvariant}, duplicated below for convenience), in order to develop the approximating functions
\begin{equation*}
q_j^{D,\star}(\bld m) =
\begin{cases}
	\frac{\lambda_j^D}{\theta^D_j}, & \mbox{ if }
\sum_{k\in \suppset_j} m_{jk}=0, \\
\frac{\lambda_j^D}{\theta^D_j}	G^D_{e,j} \Big(
	(G^D_j)^{-1}
	\Big( 1-\frac{ \sum_{k\in \suppset_j} m_{jk}}{\lambda_j^D}\Big) \Big),
	& \mbox{ if }  \sum_{k\in \suppset_j} m_{jk} \in (0,\lambda_j^D],
\end{cases}
\end{equation*}
\begin{equation*}
q_k^{S,\star} (\bld m)=
\begin{cases}
	\frac{\lambda_k^S}{\theta^S_k}, & \mbox{ if }
 \sum_{j\in \demset_k} m_{jk}=0, \\
\frac{\lambda_k^S}{\theta^S_k}	G^S_{e,k} \left(
	(G^S_k)^{-1}
	\Big(1-\frac{ \sum_{j\in \demset_k} m_{jk}}{\lambda_k^S}\right)\Big),
	& \mbox{ if }  \sum_{j\in \demset_k} m_{jk} \in (0,\lambda_k^S],
\end{cases}
\end{equation*}
\amy{for each $j \in \bbJ$ and $k \in \bbK$, recalling that  $G^D_{e,j}$ and  $G^S_{e,k}$ are the excess life distributions of the patience times.}

\begin{figure}[H]
\centering
\begin{subfigure}[t]{0.45\textwidth}
  \centering
  \includegraphics[width=\linewidth, trim = {4cm 3cm 4cm 0}, clip]{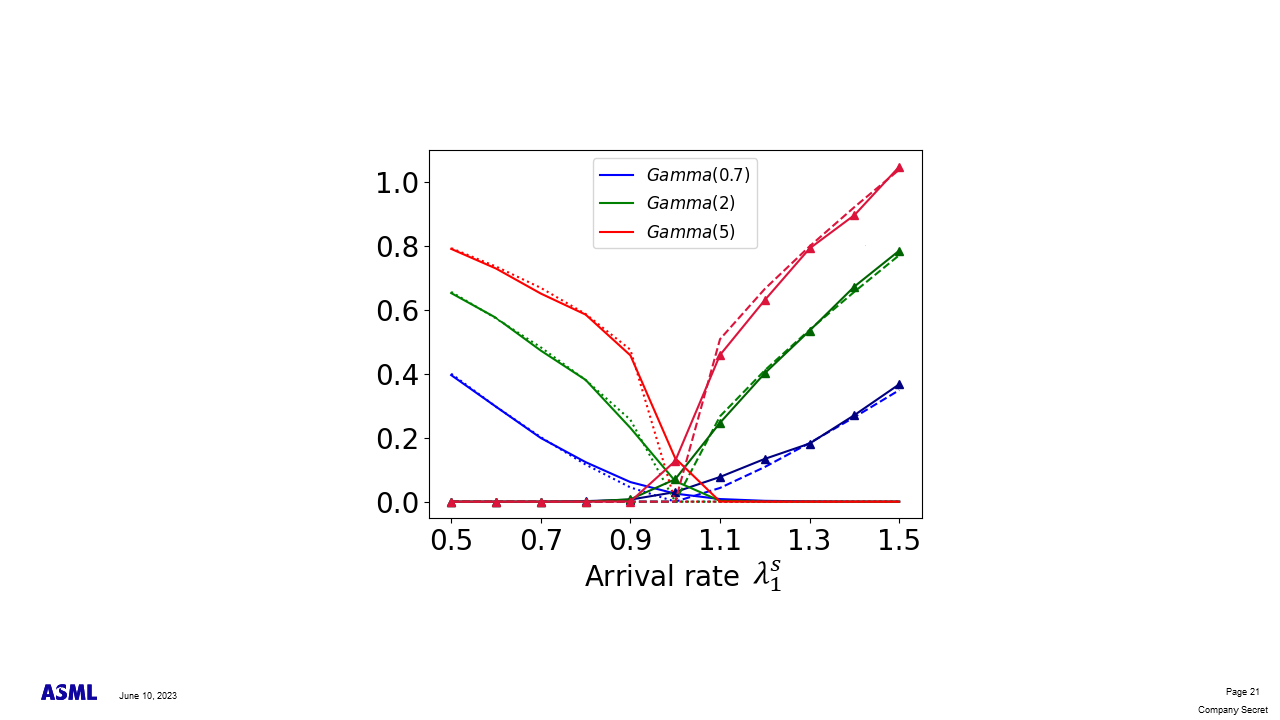}
  \caption{Gamma distributed patience times with unit mean and variance $1/x$, $x=0.7,2,5$. \amy{Average queue-lengths found via simulation.}}
  \label{figure:GammaQueues}
\end{subfigure}\hfill%
\begin{subfigure}[t]{.45\textwidth}
  \centering
  \includegraphics[width=\linewidth, trim = {4cm 3cm 4cm 0}, clip]{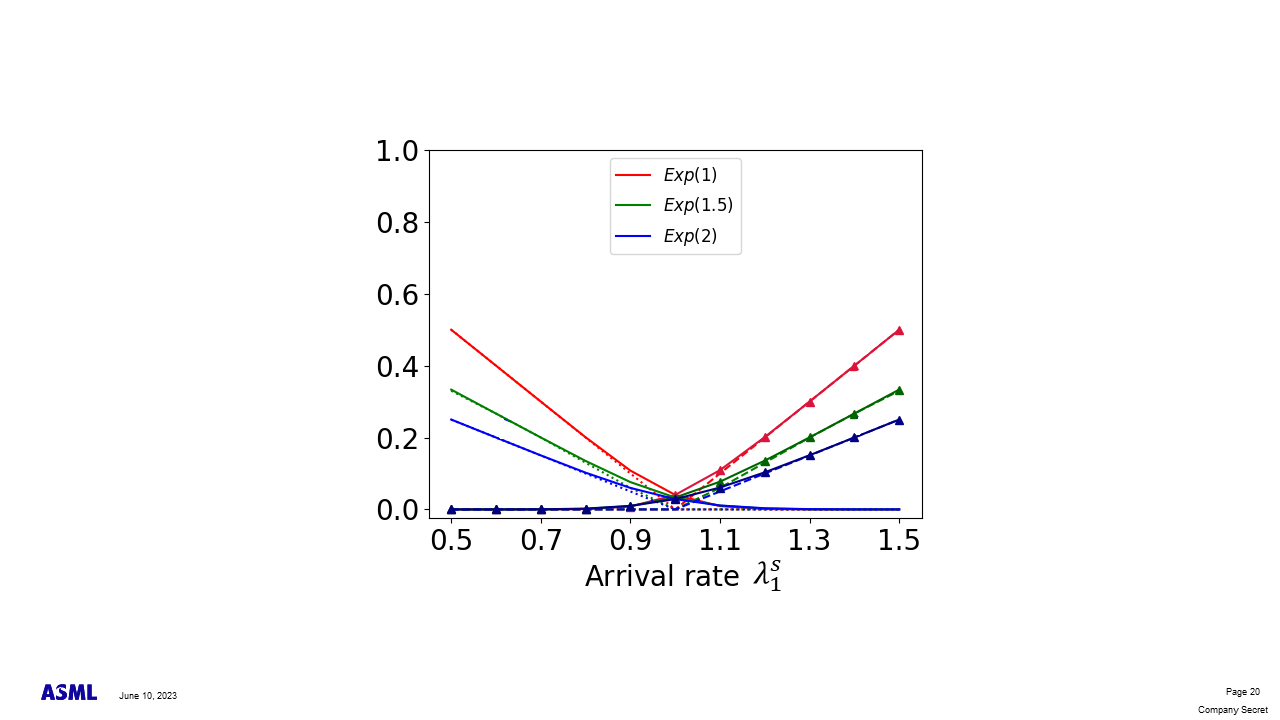}
  \caption{Exponential distributed patience times with rate $\theta = 1,1.5,2$.  \amy{Average queue-lengths found via exact calculation  (see Appendix~~\ref{sec:markovChain}).} }
  \label{figure:ExpQueues}
\end{subfigure}
\caption{The solid lines represent the average queue-lengths \amy{in the discrete-event stochastic model} and the dotted/triangle marker lines \amy{ plot the functions $q_1^{D,\star}(\mbox{min}\{\lambda_1^D,\lambda_1^S \})$ and $q_1^{S,\star}(\mbox{min}\{\lambda_1^D,\lambda_1^S\})$ for a network with one demand node and one supply node.} The parameters used are $\lambda_1^D=1$ per unit time, $t=100$, and the patience time distributions for both nodes is the same.  }
\label{figure:invariant}
\end{figure}

An optimal solution to the MP (\ref{eq:max_A_w_v3G}) provides a good approximation of the optimal instantaneous matching values if $q_j^{D,^\star}(\bld m)$ and $q_k^{S,\star}(\bld m)$  provide good approximations to the steady-state mean queue-lengths.  Figure~\ref{figure:invariant} provides supporting evidence for this claim in a network with one demand and one supply node that makes every possible match \amy{(in which case $m$ is one-dimensional and equal to $\mbox{min}\{\lambda^D_1, \lambda^S_1 \}$)}.

In general the MP is non-convex because the approximations of the queue-lengths depend on the patience time distributions through a nonlinear relationship. The only time the relationship is linear is when the patience times follow an exponential distribution, namely
\begin{equation*}
q_j^{D,\star}(\bld{m})=
\frac{1}{\theta_j^D}
\left( \lambda_j^D-  \sum_{k\in \suppset_j} m_{jk} \right)
\mbox{ and }
q_k^{S,\star}(\bld{m})=
\frac{1}{\theta_k^S}
\left( \lambda_k^S-
\sum_{j\in \demset_k} m_{jk} \right).
\end{equation*}
The last formulas have an intuitive explanation: The expressions in the parentheses represent the expected demand (supply) that is not matched, and so to find the expected demand (supply) waiting to be matched, one needs to multiply by the mean patience time.
The relationship is quadratic when the patience times are uniformly distributed in $[0,\frac{2}{\theta_j^D}]$ and  $[0,\frac{2}{\theta_k^S}]$, and is
\begin{equation*}
q_j^{D,\star}(\bld{m})=
\frac{\lambda_j^D}{\theta_j^D}
\left( 1- \left ( \frac{\sum_{k\in \suppset_j} m_{jk}}{\lambda_j^D}\right)^2 \right)
\mbox{ and }
q_k^{S,\star}(\bld{m})=
\frac{\lambda_k^S}{\theta_k^S}
\left( 1-
\left ( \frac{\sum_{j\in \demset_k} m_{jk}}{\lambda_k^S}\right)^2 \right).
\end{equation*}
\noindent
In some cases, the relationship cannot be written in a closed form, for instance when the patience times follow a gamma distribution.

Despite \eqref{eq:max_A_w_v3G} being non-convex in general, we are able to characterize the structure of the objective function of \eqref{eq:max_A_w_v3G} for a rather wide class of patience time distributions.
\begin{theorem}\label{lem:concave}
If the hazard rate functions of the patience time distributions are (strictly) increasing (decreasing), then
$q_j^{D,\star}(\bld{m})$ and
$q_k^{S,\star}(\bld{m})$
are (strictly) concave (convex) functions of the vector of matching rates
$\bld{m}$, for all $j\in \bbJ$ and $k\in \bbK$.
\end{theorem}
\begin{corollary}\label{col:convexity}
If the hazard rate functions of the patience time distributions are increasing (decreasing), then the objective function of \eqref{eq:max_A_w_v3G} is a convex (concave) function.
\end{corollary}

Having introduced the MP and studied its properties, we are now ready to introduce our proposed matching policies.

\subsection{A matching-rate-based policy}\label{sec:matching rate}
Here, we define a policy that can mimic the optimal matching rates given by \eqref{eq:max_A_w_v3G} as closely as possible, subject to the available stochastic demand and supply. Let $\bld{m}$ be a feasible point of \eqref{eq:max_A_w_v3G}. Then, to mimic the matching rates defined by $\bld{m}$, let the amount of type $j \in \bbJ$ demand we match with type $k \in \bbK$ supply, for $(j,k) \in \arcset$, in review period $i \in \{1,\ldots, \lfloor t/l \rfloor \}$ be
\begin{equation} \label{eq:DR-matching-per-period2}
\mathcal{M}_{ijk}^r := \left \lfloor m_{jk}
\min\left( l, \frac{\qd_j(il-)}{\lambda_j^D},
\frac{Q^S_k(il-)}{\lambda^S_{k}} \right)\right \rfloor\ge 0,
\end{equation}
which implies the cumulative number of matches made for  $(j,k) \in \arcset$, and
$t\geq 0$ is
\begin{equation} \label{eq:DRcumulative}
M_{jk}^r(t) :=
 \sum_{i=1}
 ^{ \lfloor t/l \rfloor } \mathcal{M}_{ijk}^r.
\end{equation}
The matching-rate-based policy is defined for each feasible point of \eqref{eq:max_A_w_v3G} and for all patience time distributions, and tries to achieve the feasible matching rates given by $\bld m$. The targeted rate during a review period between demand $j\in \bbJ$ and supply $k\in \bbK$ is $m_{jk} l$. Further, the policy takes into account the fraction of available demand and supply matched, \amy{which ensures the queue-lengths in (\ref{eq:posQ}) and (\ref{eq:posI}) satisfy the non-negativity condition.  }
\begin{proposition}\label{proposition:Madmissible2}
The proposed matching policy given by \eqref{eq:DR-matching-per-period2} and \eqref{eq:DRcumulative}
\amy{ satisfies $Q_j^D(t) \geq 0$ and $Q_k^S(t) \geq 0$ for each $j \in \bbJ$ and $k \in \bbK$ and for all $t \in \R_+$,}
 for any feasible point $\bld{m} $ to \eqref{eq:max_A_w_v3G}.
\end{proposition}
 \noindent \amy{As a result of Proposition~\ref{proposition:Madmissible2}, the matching-rate-based policy satisfies the conditions given at the end of Section~\ref{sec:Model description} for the associated matching process in (\ref{eq:DRcumulative}) to be admissible.  In particular, this is immediate to see that the matching process is integer-valued, non-decreasing, right-continuous with left limits, and does not use information about future arrivals to make matches (is non-anticipating). }

The matching-rate-based policy tries to mimic the matching rates at each edge and  discrete review point by proportionally splitting the available demand (supply) at a node between the available supply (demand) at other nodes in the same proportions as the feasible matching rates in $\bld{m}$. \amy{Then, by using $\bld{m}^\star$ to define the amount matched in (\ref{eq:DR-matching-per-period2}), we expect to achieve very close to the optimal matching rates (and we show that this is the case asymptotically, in Theorem~\ref{proposition:DRconvergence2} in Section~\ref{sec:Performance Analysis}).  However, the policy}  does not give any insight about the order in which edges should be prioritized for matching.
\levi{Prioritizing edges can be useful for platforms wishing to implement simple rules for matching, i.e., each demand (supply) type  has a fixed ranking of the supply (demand) types which it checks for availability to match at each discrete review point. Such a static priority ordering is appealing because the ordering is state-independent,
and easy to describe to a practitioner. Thus, we would like to understand when the matching-rate-based policy reduces to a static priority ordering of the edges. In the next section, we show how to interpret an optimal extreme point solution as a static priority ordering of the edges.}

\subsection{A priority-ordering policy}\label{sec:priotiry}
	In this section, we propose a \emph{priority-ordering} matching policy, which defines an order of edges over which matches are made in a greedy fashion at each discrete review point using available demand and  supply. We will show that, with a properly chosen priority-ordering, such a policy is able to imitate the matching rates determined  by an optimal extreme point solution to the MP when there exists one. \amy{From Corollary~\ref{col:convexity} to Theorem~\ref{lem:concave}}, the existence of an optimal extreme point solution is guaranteed when the patience time distributions have increasing hazard rates. In words, if supply and demand become more impatient the longer they wait in queue, then a simple priority ordering is enough to balance the trade-off between matching values and holding costs.

	
In the rest of this section, we assume the existence of an optimal extreme solution to the MP \eqref{eq:max_A_w_v3G} which we denote by $\bld\mstar$.
	
	Before constructing the priority-ordering, we need to establish a few properties of an optimal extreme point solution $\bld\mstar$. It will be helpful to define the following: we say a node $j \in \bbJ$ ($k\in \bbK$) is \emph{slack} if
	$\sum_{k \in \suppset_j}   m_{jk}^\star < \lambda_j^D$ ($\sum_{j\in \demset_k}   m_{jk}^\star < \lambda_k^S$), and \emph{tight} if $\sum_{k \in \suppset_j}  m_{jk}^\star = \lambda_j^D$ ($\sum_{ j \in \demset_k}  m_{jk}^\star = \lambda_k^S$); let $E_0 = \{(j,k)  \in \arcset: \mstar_{jk}>0\}$ denote the set of edges with positive matching rates, and let $(\{\bbJ,\bbK\},E_0)$ denote the \emph{induced graph} of $\bld\mstar$, i.e., the bipartite graph between node sets $\bbJ$ and $\bbK$ with edges corresponding to positive matching rates associated with $\bld\mstar$. Moreover, a tree is a connected acyclic undirected graph and
a forest is  an acyclic undirected graph; see \cite{gross2018}.
Then, we have the following properties of the induced graph of  an optimal extreme point solution.

\begin{lemma}\label{lem:properties}
Let $\bld\mstar$ be an optimal extreme point solution to \eqref{eq:max_A_w_v3G}. Then, the induced graph of $\bld\mstar$ has the following properties:
\begin{enumerate}
	\item It has no cycles, i.e., it is a forest, or collection of trees;
	\item Each distinct tree has at most one node with a slack constrain.
\end{enumerate}
\end{lemma}

With Lemma \ref{lem:properties} in hand we are now ready to describe the intuition for the priority-ordering we construct. Lemma \ref{lem:properties} guarantees that the set of edges, $E_0$, on which $\bld\mstar$ makes matches can be partitioned into disjoint trees, each of which have all nodes but at most one with a tight constraint (i.e., all their capacity is used up by $\bld\mstar$). Intuitively, this allows us to construct the priority sets by iteratively following a path traversing the tight nodes of each tree, removing an edge at each iteration, and updating the remaining capacities by deducting the match that $\bld\mstar$ made on that edge. The edges removed in the $h{\text{th}}$ iteration become the $h{\text{th}}$ highest priority in our ordering, and the process results in a smaller tree after each iteration so works through all edges in finite time. This process continues until it reaches the final priority set that contains all the edges with $\mstar_{jk}=0$. Before we move to the general algorithm that constructs the priority sets, we present an example.

\begin{figure}[H]
	\centering
	{\includegraphics[width=0.7\textwidth, trim={0 5cm 0 2cm},clip]{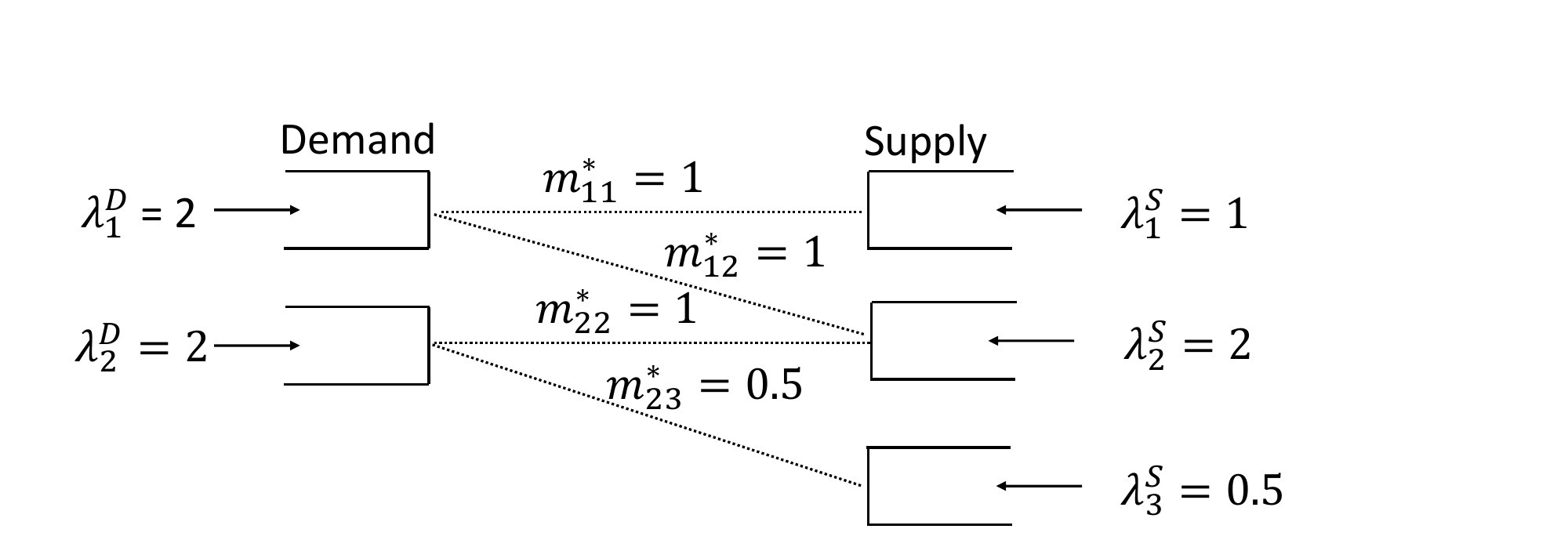}}
	\caption{A matching model with two demand and three supply nodes.}
	\label{figure:Exmpriority}
\end{figure}

\begin{examp}\label{examp:Exmpriority}
	Consider a system with two demand nodes and three supply nodes, with arrival rates as depicted in Figure~\ref{figure:Exmpriority}. The patience time distributions, values, and holding costs are such that an optimal extreme point solution to the MP for this system sets $\mstar_{11} = \mstar_{12} = \mstar_{22} = 1$ and $\mstar_{23} = 0.5$, as depicted in the figure, with all other edges having $\mstar_{jk} = 0$ \levi{(e.g., this is the case if all matching values are 1 and all holding costs are 0)}. The edges with positive matching rate, i.e., $E_0$, are drawn with dotted lines in the figure, note that these edges form a tree. For this optimal extreme point solution, we construct
	the following priority sets:
	$\calP_0(\bld{m}^\star)=\{(1,1),(2,3)\}$, $\calP_1(\bld{m}^\star)=\{(1,2)\}$, $\calP_2(\bld{m}^\star)=\{(2,2)\}$, and $\calP_3(\bld{m}^\star)=\{(1,3),(2,1)\}$. Observe that in the first iteration, the algorithm selects two edges with tight leaf nodes from the tree $($i.e., $(1,1)$ and $(2,3))$, since they are disjoint, while in the second iteration the algorithm picks only one of two available edges with tight leaf nodes
$($i.e., $(1,2)$ and $(2,2))$, since they share the end point $k = 2$, leaving the second edge until the next iteration, \levi{which is the last iteration processing edges with positive matching value. Note that if edge $(1,2)$ (resp. $(2,2)$) were chosen in the first set $\calP_0(\bld{m}^\star)$, a priority policy would make too many matches along this edge and risk leaving agents at supply node 1 (resp. 3) unmatched.
	Then the final step of the algorithm puts all edges with zero matching value in the set $\calP_3(\bld{m}^\star)$.}
\end{examp}

We formally present the construction of the priority sets in Algorithm~\ref{alg:priority_sets}, and here we provide a more detailed explanation of the intuition of the algorithm. The order we traverse the trees in $E_0$ is determined by the following observation. In each tree, at least one of the tight nodes must be a ``leaf'' of the tree, i.e., it is connected to only one edge in $E_0$ (or, equivalently only one element of $\bld\mstar$ involving this node is positive). For an example, assume that $j$ is such a leaf node in $E_0$ with a tight constraint, i.e., $\sum_k \mstar_{jk} = \lambda_j^D$, and $\mstar_{jk'}>0$ for some $k'$, while $\mstar_{jk}=0$ for all other $k \ne k'$. Then, together this implies that  $\mstar_{jk'}= \lambda_j^D$, i.e., there is only one edge, $(j,k')$, that uses up all of node $j$'s capacity. Our algorithm gives this edge highest priority, since greedily matching all available supply and demand along this edge will match $\min(\lambda_j^D,\lambda^S_{k'}) = \lambda_j^D$ (since we must have had $\lambda_j^D\le \lambda^S_{k'}$ for the initial value of $\mstar_{jk'}$ to be feasible), and this replicates the value of $\mstar_{jk'}= \lambda_j^D$. Then, removing this edge and reducing the capacities by $\mstar_{jk'}$ results in a smaller tree with the same properties as $E_0$, so we can repeat the process until we work through the whole tree. \levi{We note that \cite{kerimov2021} establishes existence of a similar priority ordering (for a model with no reneging), without explicitly constructing the priority sets as we do in Algorithm~\ref{alg:priority_sets}.} \levi{Further, we make the} following adjustment to allow for a potentially fewer number of priority sets in the final construction: more than one leaf node with a tight constraint can be put in the same priority set in a given iteration as long as their associated edges are disjoint (i.e., they do not share any endpoints).

Let $H+1$ denote the final value of the counter $h$ in Algorithm~\ref{alg:priority_sets}, so that $\calP_H(\bld{m}^\star)$ denotes the last priority set with edges $(j,k)$ such that $\mstar_{jk}>0$, and $\calP_{H+1}(\bld{m}^\star)$ denotes the final priority set with all edges $\mstar_{jk}=0$.

\begin{algorithm}[t]
	\caption{Priority Set Construction}\label{alg:priority_sets}
	\begin{algorithmic}[1]
		\State Input: supply, demand rates $\bld{\lambda}^D\in \R_+^{J}$, $\bld{\lambda}^S\in \R_+^{K}$, and $\bld{\mstar}$ an optimal extreme point solution to \eqref{eq:max_A_w_v3G}
		\State Initialize priority set counter $h=0$, edge set
$E_0 = \{(j,k) \in \arcset: \mstar_{jk}>0\}$, capacities $\bld{d} = \bld{\lambda}^D$, $\bld{s} = \bld{\lambda}^S$, and $E=E_0$
		\While{$|E|>0$}\label{alg:priority_sets_out_whl}
			\State Initialize priority set: $\calP_h(\bld{m}^\star) = \emptyset$
			\State Initialize consideration set (copy $E$): $C = E$
			\While{$|C|>0$}\label{alg:priority_sets_in_whl}
				\State Choose any $(j,k) \in C$\label{alg:priority_sets_choose_jk}
				\If{$\mstar_{jk}= d_j$ or $\mstar_{jk}= s_k$}\label{alg:priority_sets_if_tight}
					\State Add to priority set: $\calP_h(\bld{m}^\star) \leftarrow \calP_h(\bld{m}^\star) \cup \{(j,k)\}$\label{alg:priority_sets_add_Ph}
					\State Update capacities: $s_k \leftarrow s_k - \mstar_{jk}$, $d_j \leftarrow d_j - \mstar_{jk}$ \label{alg:priority_sets_updt_cap}
					\State Remove $(j,k)$ and neighbors:
$C\leftarrow C\setminus \{(j',k')\in C:j'= j \text{or } k'= k\}$ \label{alg:priority_sets_remove_nbrs}
				\Else{ Remove $(j,k)$ from consideration: $C\leftarrow C\setminus (j,k)$}
				\EndIf
			\EndWhile
			\State Remove priority set from edge set: $E \leftarrow E\setminus \calP_h$\label{alg:priority_sets_remove_Ph}
			\State Increment counter: $h \leftarrow h+1$
		\EndWhile
		\State Last priority set: $\calP_h(\bld{m}^\star) = \{(j,k) \in \arcset:\mstar_{jk}=0\}$
	\end{algorithmic}
\end{algorithm}

Note that Algorithm~\ref{alg:priority_sets} returns non-overlapping priority sets and it terminates as the following result states.
\begin{lemma}\label{prop:Alg1}
	Algorithm \ref{alg:priority_sets} runs in finite time $O\left(|E_0|^2\right)$.
\end{lemma}

The sets $\mathcal{P}_h(\bld{m}^\star)$ form a partition of the edges, $\arcset$, of the bipartite graph.  Let $\mathcal{Q}_h(\bld{m}^\star) = \cup_{l=0}^{h}\mathcal{P}_l(\bld{m}^\star)$ denote the set of all priority edges chosen up to iteration $h$. Given this partition, define the following solution to \eqref{eq:max_A_w_v3G} recursively for $(j,k) \in \mathcal{P}_h(\bld{m}^\star)$ as
\begin{align*}
y^p_{jk}(\bld{m}^\star) =
\min\left( \lambda_j^D-\sum_{k':(j,k')\in \mathcal{Q}_{h-1}(\bld{m}^\star)}
y^p_{jk'}(\bld{m}^\star),\lambda_k^S-\sum_{j':(j',k)\in \mathcal{Q}_{h-1}(\bld{m}^\star)}
y^p_{j'k}(\bld{m}^\star)
\right),
\end{align*} where an empty sum is defined to be zero.

\begin{proposition}\label{prop:property extreme}
	For $\bld{m}^\star$ an optimal extreme point solution to \eqref{eq:max_A_w_v3G}, we have $y^p_{jk}(\bld{m}^\star) = m^\star_{jk}$ for all $j\in \bbJ$ and $k\in\bbK$.
\end{proposition}

We are ready now to define recursively the priority-ordering matching policy.
The amount of type $j\in \bbJ$ demand we match with type $k\in \bbK$ supply for $(j,k)\in\mathcal{P}_0$ at the end of review period $i \ge 1$ is
\begin{equation}\label{eq:priority0}
\begin{split}
\mathcal{M}_{ijk}^p :=
\min\left( \qd_j(il-),Q^S_k(il-)
 \right).
\end{split}
\end{equation}
For $(j,k)\in\mathcal{P}_h$, $h=1,\ldots,JK$ and $i \in \{1,\ldots, \lfloor t/l \rfloor \}$, define recursively
\begin{equation}\label{eq:priorityh}
\begin{split}
\mathcal{M}_{ijk}^p :=
\min\left( \qd_j(il-)-
\sum_{k':(j,k')\in \mathcal{Q}_{h-1}(\bld{m}^\star)}
\mathcal{M}_{ijk'}^p,Q^S_k(il-)
-\sum_{j':(j',k)\in \mathcal{Q}_{h-1}(\bld{m}^\star)}
\mathcal{M}_{ij'k}^p
 \right).
\end{split}
\end{equation}
The aggregate number of matches at time $t\geq 0$ is given by
\begin{equation}\label{eq:PR-matching-cumulative}
M_{jk}^p(t) :=
\sum_{i=1}^{ \lfloor t/l \rfloor } \mathcal{M}_{ijk}^p.
\end{equation}

\amy{This is straightforward to see that the priority-ordering matching policy satisfies the admissibility conditions given at the end of Section~\ref{sec:Model description}.}

The implementation of the priority-ordering policy does not require the knowledge of the optimal extreme point solution once the priority sets have been constructed. In other words, a computer program needs to know only the priority sets to be able to implement the aforementioned policy.

\begin{remark}[Connection to the $\frac{c \mu}{\theta}$ rule.]\label{rem:cm tule}
The priority ordering studied in this section can be connected to the well-known $\frac{c \mu}{\theta}$ rule studied for multiclass  queueing systems where $\mu$ denotes the service rate; \emph{\cite{smith1956various, atar2011asymptotic}}. To see this, observe that for a system with a single supply node and exponential patience time distributions, the objective function of the MP takes the form:
$\sum_j a_j m_{j1}$ where $a_j = \p_{j1}+\frac{c_j^D}{\theta_j^D}+\frac{c_1^S}{\theta_1^S}$. Hence, its optimal solution is given by a priority rule that orders the weights $a_j$. These weights do not depend on the arrival rates, and are similar to that of the $\frac{c \mu}{\theta}$ rule but are modified to include the matching values and the fact that there is not service rate in our model. However, in the general case, the priority ordering cannot be written as a simple priority rule because of the dimension of the model and the general patience time distributions.
\end{remark}

We have proposed matching policies that use the information of an optimal solution to the MP,
which involves an approximation of the demand and supply queue-lengths. The study of the asymptotic performance of the aforementioned policies requires \amy{showing} a connection between \eqref{eq:Profit}, which involves the stochastic queue-lengths, and the objective of the MP in \eqref{eq:max_A_w_v3G}.
\amy{We do this by developing and analyzing a fluid model that approximates the demand and supply queue-lengths in
 Section~\ref{sec:fluid model}.  We then use that analysis in Section~\ref{sec:Performance Analysis} to develop an asymptotic connection between \eqref{eq:Profit} and \eqref{eq:max_A_w_v3G}.}

\section{The fluid model and its convergence to invariant states}
\label{sec:fluid model}

\amy{Section~\ref{sec:fl1} presents a fluid model that be seen as an approximation to the discrete-event stochastic model developed in Section~\ref{sec:Model description}. The invariant states of the fluid model are fixed points of the fluid model equations, and their characterization is required to develop the proposed matching policies in Section~\ref{sec:DRMUP}.  Section~\ref{sec:fl2} characterizes the fluid model invariant states and shows that a fluid model solution converges to an invariant state as time becomes large.
The fluid model presented in Section~\ref{sec:fl1} is exactly the fluid model presented in Section 3.1 in~\cite{aveklouris2021MM}, except restricted to linear arrival functions and linear matching functions. }

\amy{In this section, we require the following notation.  Given a Polish space $\bbS$, we use the notation $\mc C(\bbS)$ (with no subscript) to denote the set of $\bbS$ valued functions with domain $\Rp$ that are continuous.   We endow  $\mc C(\bbS)$  with the usual Skorokhod $J_1$-topology~\citep{billingsley1999convergence}.   Recall that for $L\in [0, \infty]$,  $\mc{M}[0,L)$
denotes the set of finite, non-negative Borel measures on $[0, L)$
endowed with the topology of weak convergence, \amy{which is a Polish space}. Given a measure
$\nu \in \mc M[0,L)$, we write $\nu[0,x] := \int_0^x \nu(dy)$. Given a measure
$\nu \in \mc M[0,L)$
and a Borel measurable function $f: [0,L] \to \R$ that is integrable with respect to $\nu$ define
$\lran{ f,\nu }:=\int_{[0,L)}f(x)\nu(dx)$.  }

\subsection{The fluid model}\label{sec:fl1}

\amy{
The fluid model state space parallels the state space for the discrete-event stochastic model in Section~\ref{sec:Model description}, except that we do not need to track the times elapsed since the last arrival of each \levi{agent} type (because arrivals occur continuously in time).  Recall that $H_j^D$ and $H_k^S$ are the right edges of the support of the patience time distribution functions for any $j \in \bbJ$ and $k \in \bbK$.  The fluid model state at time $t \in \Rp$ is described as follows:  for each $j \in \bbJ$ and $k \in \bbK$,
\begin{itemize}
\item $\overline{Q}_j^D(t)$ and $\overline{Q}_k^S(t)$ are the fluid queue-lengths at time t;
\item $\overline{\eta}_j^D(t) \in \mc{M}[0,H_j^D)$ and $\overline{\eta}_k^S(t) \in \mc{M}[0,H_k^S)$  are measure-valued functions that provide upper bounds on the fluid queue-lengths at every time $t$; that is,
analogous to \eqref{eq:state-space-restrictions},
\begin{equation} \label{in:QBs}
\ovl Q^D_j(t)\leq \lran{ 1, \ovl{\eta}_j^D(t) }, \mbox{ for all } j \in \bbJ, \mbox{ and }
\ovl Q^S_k(t)\leq \lran{ 1, \ovl{\eta}_k^S(t) }, \mbox{ for all } k \in \bbK.
\end{equation}
\end{itemize}
In summary, a fluid model solution is a vector $(\ovl{\bld{Q}}^D,\ovl{\bld{Q}}^S,
\ovl{\bld{\eta}}^D, \ovl{\bld{\eta}}^S )
\in  \mathcal{C}(\overline{\mathbb{Y}}) $
 having components $\left( \ovl Q^D_j(t),   j \in \bbJ \right)$,  $\left( \overline{Q}_k^S(t) \right)$, $\left(\overline{\eta}_j^D(t), j \in \bbJ \right)$, and $\left(  \overline{\eta}_k^S(t), k \in \bbK \right)$
that take values in the subset $\overline{\mathbb{Y}}$ of
$$\R_+^{J}\times \R_+^{K} \times
\times_{j=1}^{J}\textbf{M}[0,H^D_j)
\times \times_{j=1}^{K}\textbf{M}[0,H^S_j),$$
for which  (\ref{in:QBs})  holds for each $t \in \R_+$.
We require that a fluid model solution $(\ovl{\bld{Q}}^D,\ovl{\bld{Q}}^S,
\ovl{\bld{\eta}}^D, \ovl{\bld{\eta}}^S )
\in  \mathcal{C}(\overline{\mathbb{Y}}) $ satisfies finiteness conditions  such that for all $t \geq 0$
\begin{equation}\label{in:abrates}
\int_0^{t}   \lran{ h^D_j, \ovl{\eta}^D_j(u) } du <\infty,\mbox{ for all } j \in \bbJ, \mbox{ and }
\int_0^{t}\lran{ h^S_k, \ovl{\eta}^S_k(u) } du<\infty, \mbox{ for all } k \in \bbK,
\end{equation}
and has initial potential queue measures with no atoms; i.e.,
\begin{equation} \label{eq:fluid-no-atom-initial-condition}
\lran{ 1_{\{x\}},\overline{\eta}_j^D(0) } = 0 \mbox{ for all } x \in [0,H_j^D), j \in \bbJ \mbox{ and } \lran{ 1_{\{x\}},\overline{\eta}_k^S(0) } = 0 \mbox{ for all } x \in [0,H_K^S), k \in \bbK.
\end{equation}

The fluid analogues of the cumulative reneging processes  are explained in words below, and
are, for $t \geq 0$ and each $j \in \bbJ$ and $k \in \bbK$,
\begin{equation}\label{eq:AlRD}
\ovl{R}_j^D(t)  =
\int_0^t \int_0^{H^D_j} h^D_j(x)
\ind{\ovl{\eta}_j^D(u)[0,x]< \ovl{Q}^D_j(u)}
\ovl{\eta}_j^D(u)(dx)du
\end{equation}
and
\begin{equation}\label{eq:AlRS}
\ovl{R}_k^S(t) =
\int_0^t \int_0^{H^S_k} h^S_k(x) \ind{\ovl{\eta}_k^S(u)[0,x]< \ovl{Q}^S_k(u)}
\ovl{\eta}_k^S(u)(dx)du.
\end{equation}
The inner integrals in \eqref{eq:AlRD} and \eqref{eq:AlRS}
represent the instantaneous reneging rate, which is determined by the hazard rate function and fluid age.  Then, integrating over the instantaneous reneging rate in $[0,t]$ gives the cumulative reneging up to time $t$.
When the patience times are exponentially distributed with rates $\theta^D_j$ and  $\theta^S_k$ so that the $h_j^D(x) = \theta^D_j$ and $h_k^S(x) = \theta_k^S$ for $x \in \Rp$, we have that the instantaneous rate at which fluid leaves without being matched is linear in the queue-length; that is,
\begin{equation*}
\overline{R}_j^D(t)=
\int_0^t \theta^D_j \overline{Q}^D_j(u)du\ \mbox{ and }\
\overline{R}_k^S(t)=
\int_0^t \theta^S_k \overline{Q}^S_k(u)du.
\end{equation*}
  The condition (\ref{eq:fluid-no-atom-initial-condition}) ensures that $\ovl{R}_j^D(t)$ and $\ovl{R}_k^S(t)$ in \eqref{eq:AlRD} and \eqref{eq:AlRS} are finite for all $t \in \Rp$.
}

\amy{
The measure-valued functions $\overline{\eta}_j^D \in \mc{M}[0,H_j^D)$ and $\overline{\eta}_k^S \in \mc{M}[0,H_k^S)$ must satisfy the equations given below, and interpreted in words after, for given input to the fluid model $\bld{\lambda}^D \in (0,\infty)^{J}$ and $\bld{\lambda}^S \in (0,\infty)^{K}$.  We interpret $\lambda_j^D$ as the type $j \in \bbJ$ fluid arrival rate, and $\lambda_k^S$ as the type $k \in \bbK$ fluid arrival rate, meaning that $\lambda_j^D t$  and $\lambda_k^S t$ are the cumulative amounts of type $j$ and $k$ fluid to arrive by time $t \in \Rp$.
For any bounded function
$f\in \mathcal{C}(\R_+)$  the following integral equations hold for each $j\in \bbJ$, $k \in \bbK$, and $t \geq 0$,
\begin{equation}\label{eq:FetaD}
\lran{ f,  \ovl{\eta}_j^D(t) } =
\int_0^{H^D_j} f(x+t)
\frac{1-G^D_j(x+t)}{1-G^D_j(x)} \ovl{\eta}^D_j(0)(dx) +
 \lambda_j^D \int_0^t f(t-u)(1-G^D_j(t-u))
du,
\end{equation}
and
\begin{equation}\label{eq:FetaS}
\lran{ f,  \ovl{\eta}_k^S(t) } =
\int_0^{H^S_k} f(x+t)
\frac{1-G^S_k(x+t)}{1-G^S_k(x)} \ovl{\eta}^S_k(0)(dx)
 +
\lambda_k^S \int_0^t f(t-u)(1-G^S_k(t-u))
du.
\end{equation}
The first term of the integral equations for measures
$\overline{\eta}_j^D(\cdot)$ and $\overline{\eta}_k^S(\cdot)$ in (\ref{eq:FetaD}) and (\ref{eq:FetaS})  tracks when the patience time of fluid demand/supply initially present in the system at time zero expires, and the second term has an analogous meaning for the newly arriving fluid.

The specification of a fluid model solution for  given  arrival rates  $\bld{\lambda}^D \in (0,\infty)^{J}$ and $\bld{\lambda}^S \in (0,\infty)^{K}$ and  an initial condition
requires the specification of the matching rates, which must lie in the set

$$ \mathbb{M}:=\left\{ \bld m \in \R_+^{J \times K }:
 \sum_{k\in \suppset_j} m_{jk} \leq \lambda_j^D,
 \sum_{j\in \demset_k} m_{jk} \leq \lambda_k^S, \mbox{ and } m_{jk} = 0 \ \forall (j,k) \notin \arcset \right\}.$$
\noindent For a given matching rate matrix $\bld m \in \mathbb{M}$, $m_{jk}$ is the rate at which type $j \in \bbJ$ customer fluid is matched with type $k \in \bbK$   worker fluid, so that $m_{jk} t$ is the cumulative amount of fluid of those types matched by time $t \in \Rp$.
 Then, the fluid queue-lengths evolve as follows: for all $j \in \bbJ$, $k \in \bbK$ and $t \geq 0$,
\begin{equation}\label{eq:QFD}
\ovl{Q}^D_j(t)=\ovl{Q}^D_j(0)+\lambda^D_j t-\ovl{R}_j^D(t)
-\sum_{k\in \suppset_j} m_{jk}t,
\end{equation}
and,
\begin{equation}\label{eq:QFS}
\ovl{Q}^S_k(t)=\ovl{Q}^S_k(0)+ \lambda^S_k t-\ovl{R}_k^S(t)
-\sum_{j\in \demset_k} m_{jk}t.
\end{equation}
The equations \eqref{eq:QFD} and \eqref{eq:QFS} are the fluid analogues of the queue-length evolution equations (\ref{eq:posQ}) and (\ref{eq:posI}) in the stochastic model.

\begin{definition} \label{def:fluid-model}
Let $\bld{\lambda}^D \in (0,\infty)^{J}$ and $\bld{\lambda}^S \in (0,\infty)^{K}$ be given.
A \textbf{fluid model solution} for
arrival rates  $\left( \bld{\lambda}^D, \bld{\lambda}^S \right)$
 is
$(\ovl{\bld{Q}}^D,\ovl{\bld{Q}}^S,
\ovl{\bld{\eta}}^D, \ovl{\bld{\eta}}^S )
\in \mc \amy{\mathcal{C}}(\overline{\mathbb{Y}})$
that satisfies conditions \eqref{in:abrates} and \eqref{eq:fluid-no-atom-initial-condition},
 the integral equations \eqref{eq:FetaD} and \eqref{eq:FetaS},
and is such that there exists a matching rate matrix $\bld m \in \mathbb{M}$ for which \eqref{eq:QFD} and \eqref{eq:QFS} hold,
 with  $\ovl{\bld{R}}^D$ and $\ovl{\bld{R}}^S$  given by \eqref{eq:AlRD} and \eqref{eq:AlRS}.
\end{definition}

\begin{remark}
A fluid model solution for given arrival rates  $\left( \bld{\lambda}^D, \bld{\lambda}^S \right)$ and given initial state is unique; see Theorem 1 in~\cite{aveklouris2021MM}.
\end{remark}
}

\subsection{Invariant analysis}\label{sec:fl2}

\amy{
The matching problem (\ref{eq:max_A_w_v3G}) that approximates the objective (\ref{eq:Profit}) to maximize the long-run average matching value in the discrete event stochastic system follows by identifying the invariant states of the fluid model. In particular, the approximations for the mean steady-state queue-lengths used in (\ref{eq:max_A_w_v3G}) are invariant states of the fluid model.

\begin{definition} \label{Definition:IM}
Let $\bld{\lambda}^D \in (0,\infty)^{J}$ and $\bld{\lambda}^S \in (0,\infty)^{K}$.
A tuple
$(\bld{q}^{D,\star}, \bld{q}^{S,\star}, \bld{\eta}^{D,\star}, \bld{\eta}^{S,\star}) \in\overline{\mathbb{Y}}$ is an \textbf{invariant state} for
$(\bld{\lambda}^D, \bld{\lambda}^S)$ if the constant function
$(\ovl{\bld{Q}}^D, \ovl{\bld{Q}}^S,
\ovl{\bld{\eta}}^{D}, \ovl{\bld{\eta}}^{S})$
given by
\begin{equation}\label{eq:constant}
(\ovl{\bld{Q}}^D(t), \ovl{\bld{Q}}^S(t),
\ovl{\bld{\eta}}^{D}(t), \ovl{\bld{\eta}}^{S}(t))=
(\bld{q}^{D,\star}, \bld{q}^{S,\star}, \bld{\eta}^{D,\star}, \bld{\eta}^{S,\star}), \mbox{ for all } t \geq 0
\end{equation}
 is a fluid model solution for  $\left( \bld{\lambda}^D, \bld{\lambda}^S \right)$. The \textbf{invariant manifold} for
 $\bld \lambda = (\bld{\lambda}^D, \bld{\lambda}^S)$ is the set of all
invariant states for $(\bld{\lambda}^D, \bld{\lambda}^S)$, which we denote by $\mathcal{I}_{\bld{\lambda}}$.
\end{definition}

Recall that
$G^D_{e,j}(x)$,
$G^S_{e,k}(x)$
for
$j\in \bbJ , k\in \bbK, x\in \R_+$
are the excess life distributions of the patience times.

\begin{proposition}
\label{prop:Chofinvariant}
Let $\bld{\lambda}^D \in (0,\infty)^{J}$ and $\bld{\lambda}^S \in (0,\infty)^{K}$.
For  $\bld{m}\in \mathbb{M}$, define for each $x\in\mathbb{R}_+$,
$j \in \bbJ$, and $k \in \bbK$:
\begin{equation}\label{eq:InvDe}
\eta^{D,\star}_j(dx) := \lambda_j^D (1-G^D_j(x)) dx,
\end{equation}
\begin{equation}\label{eq:InvSu}
\eta^{S,\star}_k(dx) := \lambda_k^S (1-G^S_k(x)) dx,
\end{equation}
\begin{equation}\label{eq:InvDeQ}
q^{D,\star}_j (\bld m) :=
\begin{cases}
	\frac{\lambda_j^D}{\theta^D_j}, & \mbox{ if }
\sum_{k\in \suppset_j} m_{jk}=0, \\
\frac{\lambda_j^D}{\theta^D_j}	G^D_{e,j} \Big(
	(G^D_j)^{-1}
	\Big( 1-\frac{ \sum_{k\in \suppset_j} m_{jk}}{\lambda_j^D}\Big) \Big),
	& \mbox{ if }  \sum_{k\in \suppset_j} m_{jk} \in (0,\lambda_j^D],
\end{cases}
\end{equation}
\begin{equation}\label{eq:InvSuI}
q^{S,\star}_k(\bld m) :=
\begin{cases}
	\frac{\lambda_k^S}{\theta^S_k}, & \mbox{ if }
 \sum_{j\in \demset_k} m_{jk}=0, \\
\frac{\lambda_k^S}{\theta^S_k}	G^S_{e,k} \left(
	(G^S_k)^{-1}
	\Big(1-\frac{ \sum_{j\in \demset_k} m_{jk}}{\lambda_k^S}\right)\Big),
	& \mbox{ if }  \sum_{j\in \demset_k} m_{jk} \in (0,\lambda_k^S].
\end{cases}
\end{equation}
Then, $(\bld{q}^{D,\star}(\bld m), \bld{q}^{S,\star} (\bld m), \bld{\eta}^{D,\star}, \bld{\eta}^{S,\star}) \in\overline{\mathbb{Y}}$
is in the invariant manifold
$\mathcal{I}_{\bld{\lambda}}$ for any $\bld m \in \mathbb{M}$.  Conversely, given $\bld m \in \mathbb{M}$, if $(\bld{q}^{D}(\bld m), \bld{q}^{S}(\bld m) , \bld{\eta}^{D}, \bld{\eta}^{S}) \in \mathcal{I}_{\bld{\lambda}}$, then
\[
(\bld{q}^{D}(\bld m), \bld{q}^{S}(\bld m) , \bld{\eta}^{D}, \bld{\eta}^{S}) = (\bld{q}^{D,\star}(\bld m), \bld{q}^{S,\star}(\bld m) , \bld{\eta}^{D,\star}, \bld{\eta}^{S,\star}).
\]
\end{proposition}

Proposition~\ref{prop:Chofinvariant}  shows that there is a nonlinear relationship between the matching rate matrix
$\bld{m} \in \mathbb{M}$ and the invariant fluid queue-lengths that depends on the patience time distributions.  This is not surprising given the nonlinear relationship between the fluid queue-lengths and the patience time distribution in many server queues with abandonment (see Theorem 3.1 in \cite{whitt2006}, Theorem 3.2 in \cite{zhang2013}, Theorem 5.5 in \cite{kang2012asymptotic}, and Theorem 1 in \cite{PW:2019}).

The long time behavior of the fluid model is characterized by the invariant manifold identified in Proposition~\ref{prop:Chofinvariant}.
\begin{theorem}\label{Th:steadySMB}
For each $j \in \bbJ$ and $k \in \bbK$, assume $h_j^D$ and $h_k^S$ are bounded functions.  Suppose that $\bld{\lambda}^D \in (0,\infty)^{J}$ and $\bld{\lambda}^S \in (0,\infty)^{K}$, and
$(\ovl{\bld{Q}}^D,\ovl{\bld{Q}}^S,
\ovl{\bld{\eta}}^D, \ovl{\bld{\eta}}^S )
\in \mc \amy{\mathcal{C}}(\overline{\mathbb{Y}})$ is a fluid model solution for $\left( \bld{\lambda}^D, \bld{\lambda}^S \right)$.  Then\footnote{\amy{Note that  $\lim_{t \rightarrow \infty} \left( \ovl{\bld{\eta}}^{D}(t), \ovl{\bld{\eta}}^{S}(t)\right) = \left(\bld{\eta}^{D,\star}, \bld{\eta}^{S,\star}\right)$  denotes the weak convergence of the measures; that is,
for each $j \in \bbJ$ and all bounded $f_j \in \mathcal{C}([0,H_j^D)),
\lim_{t \rightarrow \infty} \lran {f_j,\overline{\eta}^D_j(t)} = \lran{f_j, \eta^{D,\star}_j}$ and
for each $k \in \bbK$ and all bounded $f_k \in \mathcal{C}([0,H_k^S)),
\lim_{t \rightarrow \infty} \lran {f_k,\overline{\eta}^S_k(t)} = \lran{f_k, \eta^{S,\star}_k}$
.}},
\begin{equation} \label{thrm:fluid-long-time-behavior}
\lim_{t\to \infty} \left(\ovl{\bld{Q}}^D(t), \ovl{\bld{Q}}^S(t),
\ovl{\bld{\eta}}^{D}(t), \ovl{\bld{\eta}}^{S}(t)\right)
=
\left(\bld{q}^{D,\star}(\bld m), \bld{q}^{S,\star}(\bld m), \bld{\eta}^{D,\star}, \bld{\eta}^{S,\star}\right).
\end{equation}
\end{theorem}
}

\section{Performance analysis}\label{sec:Performance Analysis}
In this section, we study the performance of the proposed matching policies in a high-volume setting when the goal is to optimize \eqref{eq:Profit}. We first derive an upper bound on the objective function in a high-volume setting in Section~\ref{sec:HVS}, then show that the matching-rate-based  and priority-ordering policies achieve this upper bound asymptotically in Section~\ref{sec:High volume}.

\amy{Throughout this section, we assume the system is initially empty.  To accomodate more general initial conditions, we would replicate  Assumption 3 in~\cite{aveklouris2021MM}, and then observe that the fluid model in Section 4 accomodates all initial conditions allowed in that assumption.  We chose not to present Assumption 3 in~\cite{aveklouris2021MM} here, and instead assume an empty initial system, to eliminate the notational overhead.
 }

\subsection{High-volume setting}\label{sec:HVS}
We consider a high volume setting in which the arrival rates of demand and supply grow large. \amy{Specifically,
 we consider a family of systems indexed by $n\in\mathbb{N},$ where $n$ tends to infinity, and let the arrival rates in the $n$th system be}
\begin{equation} \label{eq:n-system-arrival-rates}
\lambda_j^{D,n} = n \lambda_j \mbox{ for all } j \in \bbJ \mbox{ and }
 \lambda_k^{S,n} = n \lambda_k^S \mbox{ for all } k \in \bbK.
\end{equation}
Otherwise, each system $n$ has
 the same basic structure as that of the system described in
Section~\ref{sec:Model description}.  Note that the patience time distributions do not scale with $n$,
 which is consistent with the existing literature; see, for example, \cite{kang2010fluid}.

\amy{We append a superscript $n$ to the system processes associated with the system having arrival rates as given in (\ref{eq:n-system-arrival-rates}).  Then, $M^n$ is the matching process associated with the system having arrival rates in (\ref{eq:n-system-arrival-rates}) and $V_{\bld{M}^n}$ is the objective value in \eqref{eq:Profit} under $M^n$.  Also, $\left( \bld{\alpha}^{D,n}, \bld{\alpha}^{S,n}, \bld{Q}^{D,n},\bld{Q}^{S,n}, \bld{\eta}^{D,n},\bld{\eta}^{S,n} \right)$ is the associated state process.}

\amy{In the high-volume setting, the fluid model given in Section~\ref{sec:fluid model} serves as an approximation to the discrete-event stochastic model given in Section~\ref{sec:Model description}, and the large time behavior of the fluid model established in Theorem~\ref{Th:steadySMB} motivates the construction of the MP in \eqref{eq:max_A_w_v3G} in Section~\ref{sec:DRMUP}.  Hence, in high-volume, we expect that the optimal objective value of the MP provides an upper bound on the value of any admissible policy given in  \eqref{eq:Profit}, and the next result validates this intuition.
 Note that the objective function value should be of order $n$ because arrival rates are of order $n$, and so the objective function in the below is multiplied by $1/n$.}
\begin{theorem}\label{prop:upper}
\amy{For each $j \in \bbJ$ and $k \in \bbK$, assume $h_j^D$ and $h_k^S$ are bounded functions.}
Let $\bld \mstar$ be an optimal solution to \eqref{eq:max_A_w_v3G}.
\amy{If $\{ \bld{M}^n \}_{n \in \mathbb{N} }$ is a sequence of matching processes such that as $n \rightarrow \infty$
\begin{equation} \label{eq:MassumpNoHolding}
    \sup_{0 \leq t \leq T } \left| \frac{M_{jk}^{n}(t)}{n} -
m_{jk} t \right| \rightarrow 0,
\end{equation}
in probability for some $\bld{m} \in \mathbb{M}$ and for all $T\ge 0$, }, then
\begin{equation} \label{eq:max-achievable-profit-fluid}
\limsup_{t \rightarrow \infty}
\limsup_{n \rightarrow \infty} \frac{1}{nt}
 V_{\bld{M}^n}(t)
  \leq
 \sum_{(j,k) \in \arcset}
 \p_{jk} \mstar_{jk}-
 \sum_{j \in \bbJ} c_{j}^D q_j^{D,\star}( \bld \mstar)-
\sum_{k \in \bbK} c_{k}^S q_k^{S,\star} (\bld \mstar),
\end{equation}
almost surely.
\end{theorem}
\noindent \amy{The condition in Theorem~\ref{prop:upper} ensures that the scaled cumulative matches made builds up linearly in time, which is an assumption in the fluid model equations \eqref{eq:QFD} and \eqref{eq:QFS}.  }

Having established an upper bound, we investigate the existence of a policy that achieves this upper bound in the high-volume setting.
\begin{definition}
Let $\bld \mstar$ be an optimal solution to \eqref{eq:max_A_w_v3G}.
\amy{When $h_j^D$ and $h_k^S$ are bounded functions for  each $j \in \bbJ$ and $k \in \bbK$,}
an admissible policy $\bld{M}^n(\cdot)$ is \textbf{asymptotically optimal}
if
\begin{equation*}
\lim_{t\rightarrow \infty} \lim_{n\rightarrow \infty}
\frac{1}{nt}\frac{V_{\bld{M}^{n}}(t)}
 {\sum_{(j,k) \in \arcset}
 \p_{jk} \mstar_{jk}-
 \sum_{j \in \bbJ} c_{j}^D q_j^{D,\star}( \bld \mstar)-
\sum_{k \in \bbK} c_{k}^S q_k^{S,\star} (\bld \mstar)}
 =1,
\end{equation*}
in probability.
\end{definition}
In the remainder of Section~\ref{sec:Performance Analysis}, we study the asymptotic behavior of our proposed discrete review matching policies. To do so, we also need to  scale the review period length as follows
\[
    l^n =
     \frac{1}{n^{2/3}} l.
\]
The intuition behind the scaling of the length of the review period is that it should be long enough to allow us to make the most valuable matches as dictated by an MP optimal solution, but also short enough that the reneging does not hurt. \amy{The scaling of the review period length is not unique, and is connected to the number of moments assumed on the inter-arrival times.  The assumption that the inter-arrival time random variables have finite 5th moment is used to bound the number of arrivals that occur during each review period\footnote{The minimal assumption necessary should be finite $(2+\epsilon)$ moments, and then the review period length should be defined in terms of $\epsilon$.  We chose $\epsilon=3$ and avoided having $\epsilon$ appear in the definition of $l^n$. }.      }


\subsection{Asymptotic optimality of the matching policies} \label{sec:High volume}
The goal of this section is to show that the proposed matching policies defined in Section~\ref{sec:DRMUP} are asymptotically optimal. We start the analysis by examining the matching-rate-based policy which was introduced in Section~\ref{sec:matching rate}.

Given a feasible point $\bld{m}$ to \eqref{eq:max_A_w_v3G} and $t\ge 0$, the number of matches for the matching-rate-based policy in the $n$th system at discrete review period
$i \in \{1,\ldots,
\lfloor t/l^n \rfloor \}$
is given by
\begin{equation}\label{eq:defmat}
\begin{split}
\mathcal{M}_{ijk}^{r,n} := \left \lfloor
nm_{jk}
\min\left(l^n, \frac{\qdn_j(il^n-)}{\lambda_j^{D,n}},
\frac{Q_k^{S,n}(il^n-)}{\lambda_k^{S,n}} \right)
\right \rfloor
=
\left \lfloor
m_{jk}
\min\left( n^{1/3}l, \frac{\qdn_j(il^n-)}{\lambda_j^D},
\frac{Q_k^{S, n}(il^n-)}{\lambda^S_{k}} \right)
\right \rfloor.
\end{split}
\end{equation}
Note that the instantaneous matching rates $m_{jk}$ are scaled by $n$ in \eqref{eq:defmat} because they remain feasible to \eqref{eq:max_A_w_v3G} when the arrival rates are scaled by a factor of $n$. The cumulative number of matches until time $t$ in the $n$th system is given by
\begin{equation}\label{eq:defmat2}
M_{jk}^{r,n}(t) :=
 \sum_{i=1}
 ^{ \lfloor t/l^n \rfloor }
  \mathcal{M}_{ijk}^{r,n}.
\end{equation}

Our goal is to show that the matching-rate-based policy is asymptotically optimal.
The following result shows that the scaled number of matches approaches asymptotically any feasible point of MP, and it is one of the keys to prove the asymptotic optimality.
\begin{theorem} \label{proposition:DRconvergence2}
Let $\bld{m}$ be a feasible point to \eqref{eq:max_A_w_v3G}. Under the matching-rate-based policy \eqref{eq:defmat} and \eqref{eq:defmat2}, as $n \rightarrow \infty$,
\[
    \sup_{0 \leq t \leq T } \left| \frac{M_{jk}^{r,n}(t)}{n} -
m_{jk} t \right| \rightarrow 0,
\]
in probability for all $T\ge 0$, $(j,k) \in \arcset$.
\end{theorem}
\noindent Although the detailed proof of Theorem~\ref{proposition:DRconvergence2}  is given in the appendix, we provide the reader with a brief outline here to establish the basic template used in the proofs of this section. The proof proceeds in three basic steps: i) establish an upper bound on the amount of reneging during a review period, ii) use the reneging upper bound to derive a lower bound on the number of matches made during a review period, and iii) show that this lower bound on matches made approaches the desired feasible point of the MP.

\amy{Theorem~\ref{proposition:DRconvergence2} holds for an optimal solution $\bld{m}^\star$ to \eqref{eq:max_A_w_v3G}, and guarantees that the condition required in Theorem~\ref{prop:upper} is satisfied.
Moreover, Theorem 3 in~\cite{aveklouris2021MM} establishes that a fluid model solution arises as a limit point of the fluid scaled state descriptors\footnote{\amy{Since arrival rates are of order $n$, the state descriptor should be divided by $n$}.} $\left( \bld{Q}^{D,n}/n,\bld{Q}^{S,n}/n, \bld{\eta}^{D,n}/n,\bld{\eta}^{S,n}/n \right)$, and Theorem~\ref{Th:steadySMB} ensures that that fluid model solution approaches the invariant point defined in Proposition~\ref{prop:Chofinvariant} for $\bld{m}^\star$.  We conclude that the matching-rate based policy is asymptotically optimal.}

\begin{corollary}\label{col:AOMR}
\amy{For each $j \in \bbJ$ and $k \in \bbK$, assume $h_j^D$ and $h_k^S$ are bounded functions.}
The matching-rate-based policy is asymptotically optimal for an optimal solution, $\bld{m}^\star$, to \eqref{eq:max_A_w_v3G}.
\end{corollary}

Having studied the asymptotic performance of the matching-rate-based policy, we now move to the asymptotic optimality of the priority-ordering policy introduced in Section~\ref{sec:priotiry}. In this case, the number of matches made in
the $n$th system between type $j\in \bbJ$ demand and type $k\in \bbK$ supply for $(j,k)\in\mathcal{P}_h$, $h=0,1,\ldots,H+1$ at the end of review period  $i \in \{1,\ldots,\lfloor t/l^n \rfloor \}$ is given by
\begin{equation}\label{eq:priorityhn}
\begin{split}
\mathcal{M}_{ijk}^{p, n} :=
\min\left( \qdn_j(il^n-)-
\sum_{k':(j,k')\in \mathcal{Q}_{h-1}(\bld{m}^\star)}
\mathcal{M}_{ijk'}^{p,n},Q_k^{S,n}(il^n-)
-\sum_{j':(j',k)\in \mathcal{Q}_{h-1}(\bld{m}^\star)}
\mathcal{M}_{ij'k}^{p,n}
 \right),
\end{split}
\end{equation}
where an empty sum is defined to be zero. Hence, the aggregate number of matches at time $t\geq 0$ is given by
\begin{equation}\label{eq:PR-matching-cumulativen}
M_{jk}^{p,n}(t) :=
\sum_{i=1}^{ \lfloor t/l^^n \rfloor } \mathcal{M}_{ijk}^{p,n}.
\end{equation}

\amy{The following result states that the priority-ordering policy asymptotically achieves the optimal matching rates, and its corollary, which is analogous to Corollary~\ref{col:AOMR} to Theorem~\ref{proposition:DRconvergence2}, shows that the priority-ordering policy is asymptotically optimal.}
\begin{theorem}\label{Thm:aopriority}
Suppose that there exists an optimal extreme point solution to \eqref{eq:max_A_w_v3G} denoted by $\bld m^{\star}$  and let the priority sets be given by Algorithm~\ref{alg:priority_sets}.
Under the priority-ordering matching policy
\eqref{eq:priorityhn} and \eqref{eq:PR-matching-cumulativen},
as $n \rightarrow \infty$,
\[
    \sup_{0 \leq t \leq T } \left| \frac{M_{jk}^{p,n}(t)}{n} -
m_{jk}^{\star} t \right| \rightarrow 0,
\]
in probability for all $T\ge 0$, $(j,k) \in \arcset$.
\end{theorem}
\noindent The proof of the last theorem in based on a two level induction. We first proceed in an induction on the priority sets, and then for any set we proceed in a second induction on the review periods inside this set. For each review period, we show that the number of matches made is close to the optimal matching rates given by the MP. Now, the asymptotic optimality of the priority-ordering policy follows in the same way as in the matching-rate-based policy.
\begin{corollary}\label{col:AOPOP}
\amy{For each $j \in \bbJ$ and $k \in \bbK$, assume $h_j^D$ and $h_k^S$ are bounded functions.}
The priority-ordering policy is asymptotically optimal for an optimal extreme point solution, $\bld{m}^\star$, to \eqref{eq:max_A_w_v3G}.
\end{corollary}

Theorems~\ref{proposition:DRconvergence2} and \ref{Thm:aopriority} suggest that the quantity $m_{jk}^{\star} t$ can be seen as an approximation of the mean number of matches at time $t$ for $j\in \bbJ$ and $k \in \bbK$. 

\begin{remark} \label{rem:bounded-hazard}
\amy{The matching-rate-based policy and priority-ordering policy can be defined for any patience time distributions that have a density, a finite mean, and a strictly increasing cdf, because these are the requirements needed to define the fluid model in Definition~\ref{def:fluid-model} and its invariant states in Proposition~\ref{prop:Chofinvariant}.
Even though Corollary~\ref{col:AOMR} to Theorem~\ref{proposition:DRconvergence2}  and Corollary~\ref{col:AOPOP} to  Theorem~\ref{Thm:aopriority} require bounded hazard functions, we can still simulate the performance of the policies when the hazard functions are unbounded (as, for example, is true for the uniform distribution), and we do this in Section~\ref{sec:SimJK} below.}
\end{remark}

\section{Simulation Experiments} \label{sec:SimJK}

In this section, we present simulation results on the behavior of the aforementioned matching policies \amy{in matching models with small numbers of nodes that have non-trivial behavior}. We consider the objective \eqref{eq:Profit} of the matching policies, and illustrate its behavior in several different parameter regimes.  
Specifically, we consider two broad sets of simulations; the first assesses policy performance in a network as the scaling parameter $n$ grows, the second investigates how to practically set the review period length in a system for a fixed $n$.

\subsection{The effect of increasing arrival rates\\}

\begin{figure}[H]
\centering
\begin{subfigure}[t]{0.30\textwidth}
  \centering
  \includegraphics[width=\linewidth]{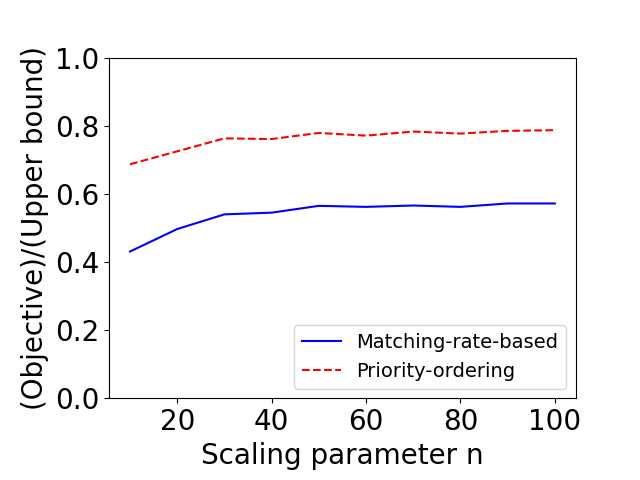}
  \caption{$l^n=0.3$ for all $n$.}
  \label{fig:GammaSiml01}
\end{subfigure}
\begin{subfigure}[t]{0.30\textwidth}
  \centering
  \includegraphics[width=\linewidth]{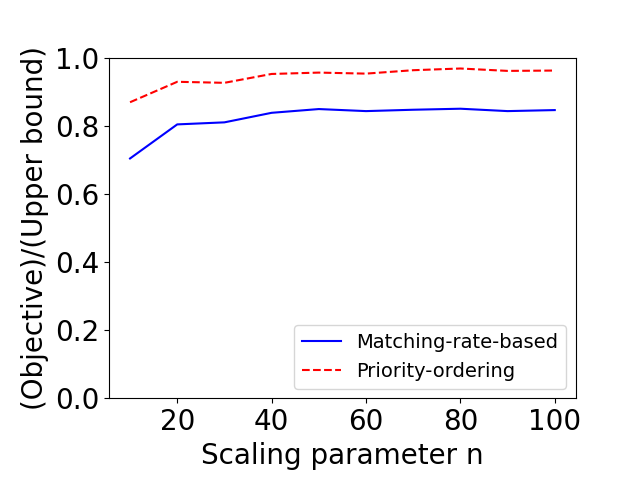}
  \caption{$l^n=0.1$ for all $n$.}
  \label{fig:GammaSiml01}
\end{subfigure}
\begin{subfigure}[t]{.30\textwidth}
  \centering
  \includegraphics[width=\linewidth]{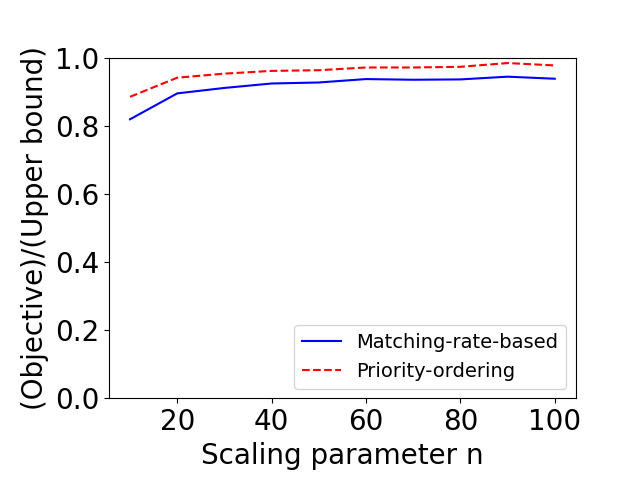}
  \caption{$l^n=0.01$ for all $n$.}
  \label{fig:GammaSiml001}
\end{subfigure}
\caption{Gamma distributed patience times with mean 1/3 and variance 1/27.}
\label{fig:GammaSim}
\end{figure}

\begin{figure}[H]
\centering
\begin{subfigure}[t]{0.30\textwidth}
  \centering
  \includegraphics[width=\linewidth]{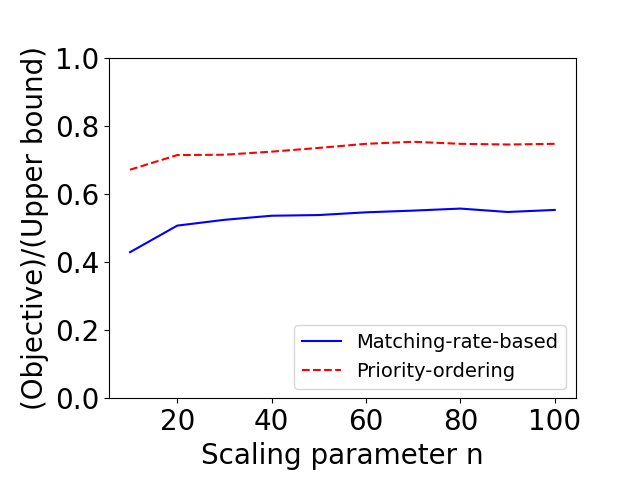}
  \caption{$l^n=0.3$  for all $n$.}
  \label{fig:UniformSiml01}
\end{subfigure}
\begin{subfigure}[t]{0.30\textwidth}
  \centering
  \includegraphics[width=\linewidth]{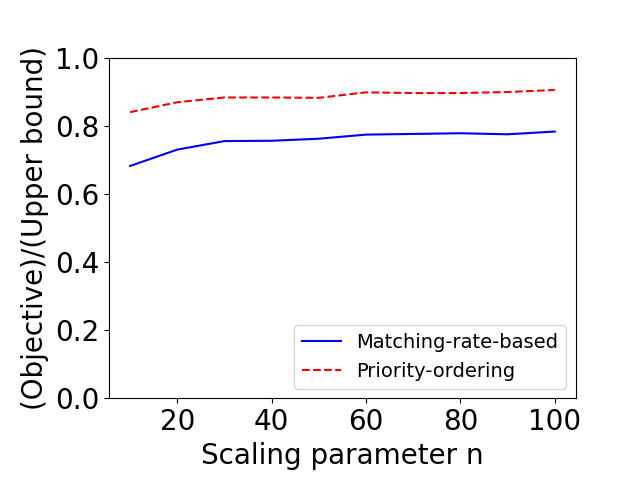}
  \caption{$l^n=0.1$ for all $n$.}
  \label{fig:UniformSiml01}
\end{subfigure}
\begin{subfigure}[t]{.30\textwidth}
  \centering
  \includegraphics[width=\linewidth]{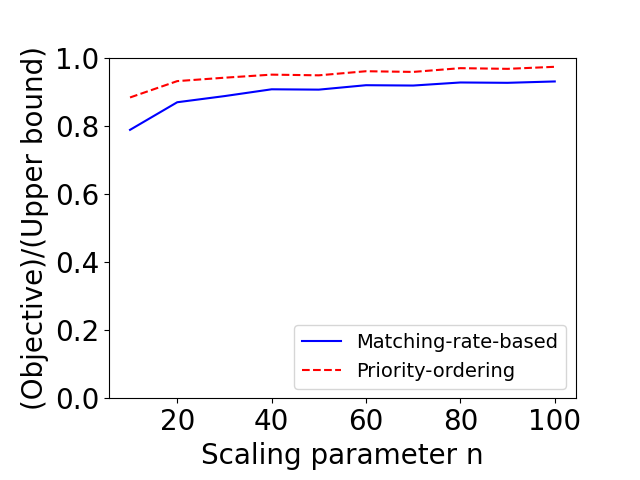}
  \caption{$l^n=0.01$ for all $n$.}
  \label{fig:UniformSiml001}
\end{subfigure}
\caption{Uniform distributed patience times with mean 1/3 and variance 1/27.}
\label{fig:UniformSim}
\end{figure}

For our first set of simulations, we consider a model with four demand and four supply nodes and Poisson arrivals with rates $\bld{\lambda}^D=(3, 2, 1, 3)$ and $\bld{\lambda}^S=(2, 2, 2, 2)$. Further, we assume that the system is initially empty and we denote the value vector for each demand node by $\bld v_{j}= (v_{j1}, \ldots, v_{jK})$ for each $j\in \bbJ$. The values at each edge are given by the vectors: $\bld v_{1} = (1, 2, 3, 1)$, $\bld v_{2}= (1, 1, 1, 1)$, $ \bld v_{3}=(2, 1, 1, 2)$, and
$\bld v_{4}= (3, 3, 2, 1)$. The holding costs for demand and supply nodes are
$\bld c^D=(1, 2, 1, 2)$ and $\bld c^S=(2, 1, 2, 1)$, respectively.
All the nodes have the same patience time distribution with the same mean and variance, i.e., $\theta^D_j = \theta^S_k = \theta$. We consider  two different distributions of the patience times both having an increasing hazard rate function; uniform in $[0,2/\theta]$ and gamma with shape parameter equal to 3 and scale parameter equal to $1/3\theta$, i.e., the mean patience time is $1/\theta$.
We fix the time horizon $t=100$ in \eqref{eq:Profit}.
Figures~\ref{fig:GammaSim} and \ref{fig:UniformSim} show the ratio of average objective of the matching policy to the optimal objective of \eqref{eq:max_A_w_v3G}
for the two policies as a function of the scaling parameter  $n$  for various values of the discrete review period $l^n$.
To illustrate the impact of changing both the review period length, $l^n$, and the scaling parameter, $n$, in each chart of this section we hold the review period length $l^n$ constant while letting the scaling parameter $n$ grow (i.e., we do not let the review period length grow with the arrival rate as it did in Section \ref{sec:HVS}). In other words, in each chart we hold the discrete review length constant while letting the arrival rates $\boldsymbol{\lambda}^{D,n}$ and $\boldsymbol{\lambda}^{S,n}$ grow large. \amy{We do this because our chosen scaling for $l^n$ is not unique, and we wanted to illustrate more broadly the review period length impact.}
Below, we summarize the main observations from the simulations.


\textbf{The matching policies have better performance for small review period length.}
\amy{Figures~\ref{fig:GammaSim} and \ref{fig:UniformSim} show that the policies perform better (are much closer to the upper bound) when the review period length is small, a phenomenon we  explore further in the simulations of the next subsection.  Moreover, Figures~\ref{fig:GammaSiml001} and~\ref{fig:UniformSiml001} show that the policy performance can be good for small arrival rates,}
 even though the asymptotic optimality results (Theorem~\ref{proposition:DRconvergence2} and \ref{Thm:aopriority}) hold for large arrival rates.

\textbf{The priority-ordering policy seems to perform better than the matching-rate-based policy.}
In all figures, the priority-ordering policy behaves better than the matching-rate-based policy. This can be explained intuitively as follows: First, the priority-ordering policy makes matches at each review period according to the priority of the edges that identifies the most valuable edges taking into account the holding costs as well. Second, the priority-ordering policy exhausts the demand or supply at each review period in contrast to the matching-rate-based policy. Based on this observation, we focus on the priority-ordering policy in our next simulation.


\subsection{The effect of decreasing review period length\\}
Next, we investigate the impact of the review period on the priority-ordering policy performance for a fixed $n$ (i.e., a non-asymptotic/high-volume setting).
While we characterized an asymptotically optimal rate for the review period as a function of $n$ in Section 5.1, we show here that there is an optimal review period length in a finite system that balances the trade-off between matching quickly to avoid reneging, and waiting to gather enough agents to make better matches. We numerically investigate this, and demonstrate that both the network structure and \amy{patience time} distributions play a critical role in the optimal review period length.
%
%
We analyze two models with different network structures. To focus on the review period, we minimize the number of parameters in the systems we consider.



\begin{figure}[ht]
	\centering
	\begin{subfigure}{.49\textwidth}
		\centering
		\begin{tikzpicture}[font=\sf,scale=0.55]
		\tikzstyle{every node}=[font=\scriptsize]
		
		\node[draw,circle,label=below:{$k=1$}] (n0) at ({2*sqrt(3)},2*1){};
		\node[draw,circle,label=below:{$j=1$}] (n1) at ({-2*sqrt(3)},2*1.5){};
		\node[draw,circle,label=below:{$j=2$}] (n2) at ({-2*sqrt(3)},2*0.5){};
		\node[draw,circle,label=below:{$k=2$}] (n3) at ({2*sqrt(3)},2*0){};
		
		\draw[thick,-latex] (n0) -- (n1);
		\draw[thick,-latex] (n2) -- (n0);
		\draw[thick,-latex] (n0) -- (n2);
		\draw[thick,-latex] (n2) -- (n3);
		
		\draw[red,dashed,thick,-latex] (-6,2*1.5) .. controls ($(n1)+(-0.2cm,0cm)$) and ($(n1)-(1cm,0cm)$) .. (-4,2*1.5);
		
		\draw[green,dashed,thick,-latex] (-6,2*0.5) .. controls ($(n2)+(-0.2cm,0cm)$) and ($(n2)-(1cm,0cm)$) .. (-4,2*0.5);
		
		\draw[green,dashed,thick,-latex] (6,2*1) .. controls ($(n0)+(0.2cm,0cm)$) and ($(n0)-(-1cm,0cm)$) .. (4,2*1);
		
		\draw[red,dashed,thick,-latex] (6,2*0) .. controls ($(n3)+(0.2cm,0cm)$) and ($(n3)-(-1cm,0cm)$) .. (4,2*0);
		
		\node at (0,3) {$v_{11} = 0.1$};
		\node at (-0.3,1.8) {$v_{21} = 1$};
		\node at (0,0) {$v_{22} = 0.1$};
		
		\node at ($(n0) + (1.5cm, 0.3cm)$) {$\lambda^S_{1}$};
		\node at ($(n1) + (-1.5cm, 0.3cm)$) {$\lambda^D_{1}$};
		\node at ($(n2) + (-1.5cm, 0.3cm)$) {$\lambda^D_{2}$};
		\node at ($(n3) + (1.5cm, 0.3cm)$) {$\lambda^S_{2}$};
		
		\node at ($(n0) + (2cm, -0.3cm)$) {$U(0, 2/\theta_1^S)$};
		\node at ($(n1) + (-2cm, -0.3cm)$) {$U(0, 2/\theta_1^D)$};
		\node at ($(n2) + (-2cm, -0.3cm)$) {$ U(0, 2/\theta_2^D)$};
		\node at ($(n3) + (2cm, -0.3cm)$) {$U(0, 2/\theta_2^S)$};
		
		
		\end{tikzpicture}
		\caption{Demand with Supply Alternatives}
		\label{Model 1 - draw}
	\end{subfigure}%
	\begin{subfigure}{.49\textwidth}
		\vspace{.8\baselineskip}
		\begin{tikzpicture}[font=\sf,scale=0.55]
		\tikzstyle{every node}=[font=\scriptsize]
		
		\node[draw,circle,label=below:{$k=1$}] (n0) at ({2*sqrt(3)},2*1){};
		\node[draw,circle,label=below:{$j=1$}] (n1) at ({-2*sqrt(3)},2*1.5){};
		\node[draw,circle,label=below:{$j=2$}] (n2) at ({-2*sqrt(3)},2*0.5){};
		
		\draw[thick,-latex] (n0) -- (n1);
		\draw[thick,-latex] (n2) -- (n0);
		\draw[thick,-latex] (n0) -- (n2);
		
		\draw[red,dashed,thick,-latex] (-6,2*1.5) .. controls ($(n1)+(-0.2cm,0cm)$) and ($(n1)-(1cm,0cm)$) .. (-4,2*1.5);
		
		\draw[green,dashed,thick,-latex] (-6,2*0.5) .. controls ($(n2)+(-0.2cm,0cm)$) and ($(n2)-(1cm,0cm)$) .. (-4,2*0.5);
		
		\draw[green,dashed,thick,-latex] (6,2*1) .. controls ($(n0)+(0.2cm,0cm)$) and ($(n0)-(-1cm,0cm)$) .. (4,2*1);
		
		\node at (0,3) {$v_{11} = 0.1$};
		\node at (-0.3,1.8) {$v_{21} = 1$};
		
		\node at ($(n0) + (1.5cm, 0.3cm)$) {$\lambda^S_{1}$};
		\node at ($(n1) + (-1.5cm, 0.3cm)$) {$\lambda^D_{1}$};
		\node at ($(n2) + (-1.5cm, 0.3cm)$) {$\lambda^D_{2}$};
		
		\node at ($(n0) + (2cm, -0.3cm)$) {$ U(0, 2/\theta_1^S)$};
		\node at ($(n1) + (-2cm, -0.3cm)$) {$ U(0, 2/\theta_1^D)$};
		\node at ($(n2) + (-2cm, -0.3cm)$) {$ U(0, 2/\theta_2^D)$};
		
		
		\end{tikzpicture}
		\caption{Demand without Supply Alternatives}
		\label{Model 2 - draw}
	\end{subfigure}
\caption{Two Network Models}
\label{fig:Models_1_and_2}
\end{figure}

The two networks are depicted in Figures \ref{Model 1 - draw} and \ref{Model 2 - draw} and are quite similar; Figure \ref{Model 2 - draw} is simply Figure \ref{Model 1 - draw} with supply node 2 removed. Thus, we describe the inputs for Figure \ref{Model 1 - draw} in detail, and note that Figure \ref{Model 2 - draw} has all the same inputs excluding those for supply node 2. We compare these two models to illustrate how the network structure impacts the optimal review period length.

We now describe Figure \ref{Model 1 - draw}, which consists of two demand and two supply nodes. We simulate three different cases with increasing patience times for each model, which we observe
impacts the importance of the review period length.
Each node has patience times that are uniformly distributed in $[0, 2/\theta_j^D]$ and $[0, 2/\theta_k^S]$. Defining reneging rate vectors $\boldsymbol{\theta^D} = (\theta_1^D, \cdots, \theta_J^D)$ and $\boldsymbol{\theta^S} = (\theta_1^S, \cdots, \theta_K^S)$, we consider the following three cases: i) $\boldsymbol{\theta^D} = (0.2, 2)$ and $\boldsymbol{\theta^S} = (2, 0.2)$, ii) $\boldsymbol{\theta^D} = (0.1, 0.2)$ and $\boldsymbol{\theta^S} = (0.2, 0.1)$, and iii) $\boldsymbol{\theta^D} = (0.05, 0.1)$ and $\boldsymbol{\theta^S} = (0.1, 0.05)$.
For each of the three cases, we consider the same Poisson arrival rates $\boldsymbol{\lambda}^D = (1, 0.1)$ and $\boldsymbol{\lambda}^S = (0.1, 1)$, for the demand and supply nodes, respectively. We assume that the system is initially empty, consider $\boldsymbol{v}_1 = (0.1)$ and $\boldsymbol{v}_2 = (1, 0.1)$ as the value vectors for each demand node, and assume no holding costs.  
We fix the time horizon $t = 1000$ in \eqref{eq:Profit} and run 120 sample paths for each case.

\subsubsection{Demand without supply alternatives}
We begin by simulating the priority ordering policy for the network in Figure \ref{Model 2 - draw}, which has all the same inputs as Figure \ref{Model 1 - draw} with supply node 2 removed. For each case in this network, the MP solution simply makes all matches between demand node 2 and supply node 1, which we call \textit{high-value matches}. 
Thus, the priority sets are $\mathcal{P}_0 = \{(2,1)\}$ and $\mathcal{P}_1 = \{(1,1)\}$, and we call the second set \textit{low-value matches}.
Note that for this network structure, if a worker at supply node $1$ makes a low-value match with demand node $1$, then we lose the ability to make a high-value match between that worker and a customer at demand node $2$ in the future. 
Therefore, for this model, we refer to workers at the supply node $1$ and customers at the demand node $2$ as \textit{high-value agents}, and customers at the demand node $1$ as \textit{low-value agents}.

\begin{figure}[!htb]
	\centering
	\begin{subfigure}[t]{0.32\textwidth}
		\centering
		\includegraphics[width=\linewidth]{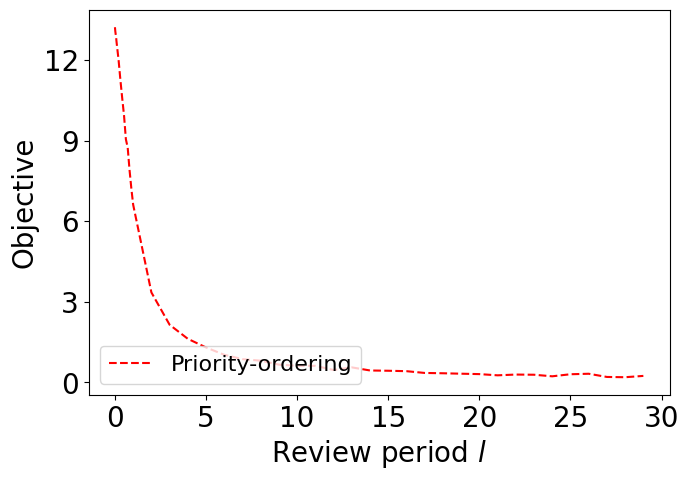}
		\caption{$\boldsymbol{\theta^D} = (0.2, 2)$, $\boldsymbol{\theta^S} = 2$.}
		\label{fig:Model2_Results_a}
	\end{subfigure}%
\begin{subfigure}[t]{0.32\textwidth}
	\centering
	\includegraphics[width=\linewidth]{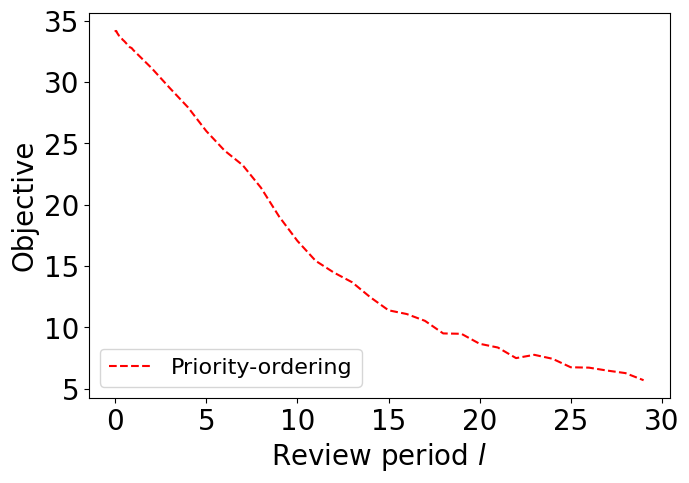}
	\caption{$\boldsymbol{\theta^D} = (0.1, 0.2)$, $\boldsymbol{\theta^S} = 0.2$.}
	\label{fig:Model2_Results_b}
\end{subfigure}
\begin{subfigure}[t]{0.32\textwidth}
	\centering
	\includegraphics[width=\linewidth]{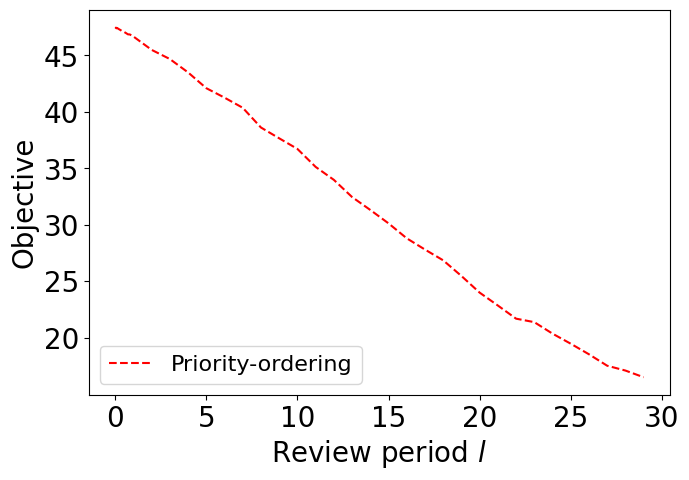}
	\caption{$\boldsymbol{\theta^D} = (0.05, 0.1)$, $\boldsymbol{\theta^S} = 0.1$.}
	\label{fig:Model2_Results_c}
\end{subfigure}
	\caption{Performance of the priority-ordering policy for Model 2.}
	\label{fig:Model2_Results}
\end{figure}

\begin{figure}[!htb]
	\centering
	\begin{subfigure}[t]{0.32\textwidth}
		\centering
		\includegraphics[width=\linewidth]{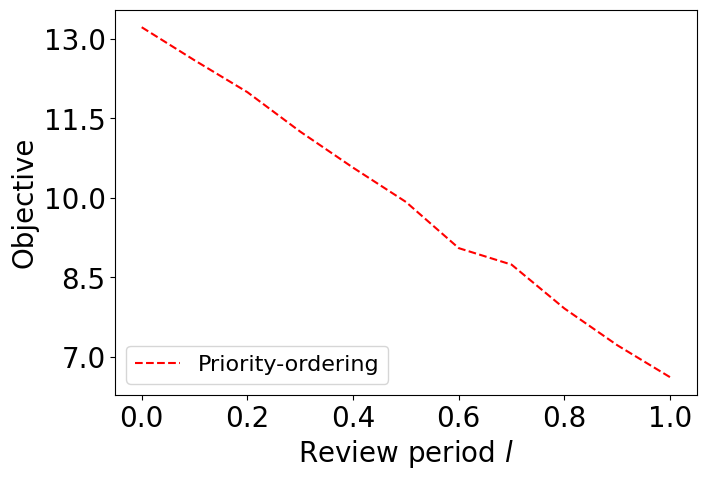}
		\caption{$\boldsymbol{\theta^D} = (0.2, 2)$, $\boldsymbol{\theta^S} = 2$.}
		\label{fig:Model2_Results_zoom_a}
	\end{subfigure}%
	\begin{subfigure}[t]{0.32\textwidth}
		\centering
		\includegraphics[width=\linewidth]{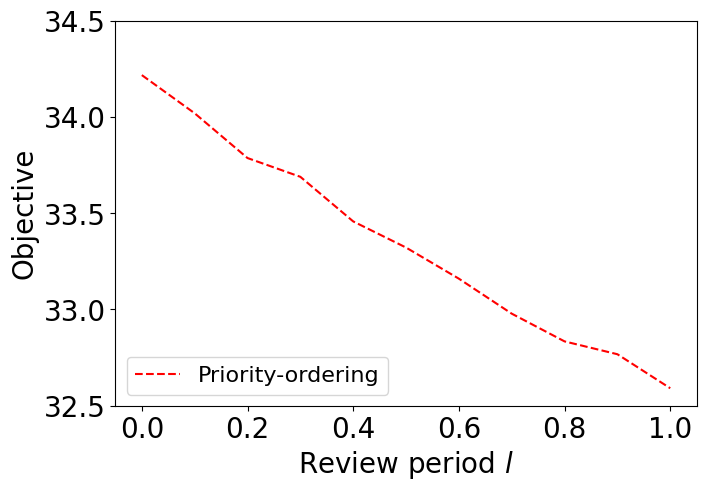}
		\caption{$\boldsymbol{\theta^D} = (0.1, 0.2)$, $\boldsymbol{\theta^S} = 0.2$.}
		\label{fig:Model2_Results_zoom_b}
	\end{subfigure}
	\begin{subfigure}[t]{0.32\textwidth}
		\centering
		\includegraphics[width=\linewidth]{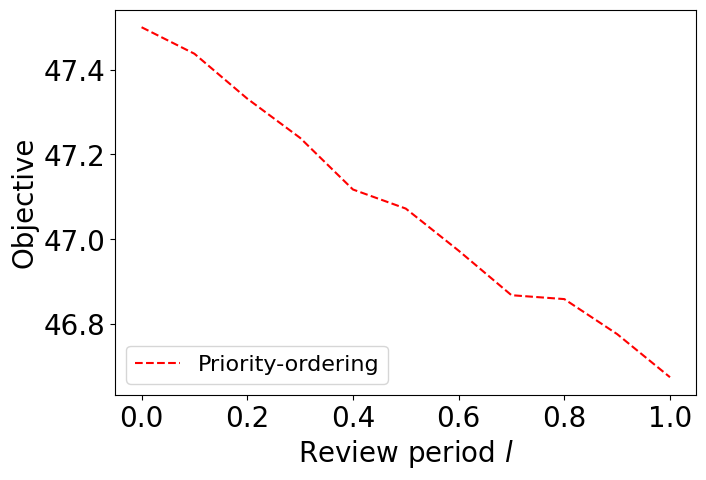}
		\caption{$\boldsymbol{\theta^D} = (0.05, 0.1)$, $\boldsymbol{\theta^S} = 0.1$.}
		\label{fig:Model2_Results_zoom_c}
	\end{subfigure}
	\caption{Performance of the priority-ordering policy for Model 2 (zoom in the x-axis between the values 0 - 1).}
	\label{fig:Model2_Results_zoom}
\end{figure}

%
We observe in this network that shorter review periods are better, regardless of the patience time distributions, as illustrated in Figures \ref{fig:Model2_Results} and \ref{fig:Model2_Results_zoom} (the figures of this subsection show the objective value over the time horizon on the y-axis; we do not use the ratio compared to the upper bound because we are not considering the asymptotic regime).
This is because customers at demand node 2 have no alternative match, so these agents will stay in the system until they renege or make a high value match with supply node 1, regardless of the review period length. Thus, with this network structure, the main focus of the policy should be on making matches quickly.


Following this line of thought, we may conjecture that for Figure \ref{Model 2 - draw}, the review period length should be zero, meaning we make matches immediately when agents arrive, if possible. We test this conjecture with an additional simulation for Figure \ref{Model 2 - draw} that focuses on review period lengths on a smaller scale,  $l \in [0,1]$. In particular, we simulate the priority policy for $l \in \{0.1, 0.2, \dots, 1\}$, and also simulate a priority policy with $l=0$, i.e., a match is made (if possible) following the priority ordering each time an agent arrives to the system. These simulation results for Figure \ref{Model 2 - draw} are depicted in Figure \ref{fig:Model2_Results_zoom}, which demonstrate that the same intuition holds at this smaller scale as the review period converges to zero: Shorter review periods are better because there is no benefit to waiting longer to gather more valuable matches in this system. This observation begs the question: Are shorter review periods always better? We explore this question in the next set of simulations on Figure \ref{Model 1 - draw}.

\subsubsection{Demand with supply alternatives}

We now simulate the priority ordering policy for Figure \ref{Model 1 - draw}, which has an additional supply node compared to Figure \ref{Model 2 - draw}, that we shall see makes the length of the review period more important. For each case of this network, the MP solution simply makes all matches between demand node 2 and supply node 1, which we call \textit{high-value matches}. Thus, the priority sets are $\mathcal{P}_0 = \{(2,1)\}$ and $\mathcal{P}_1 = \{(1,1), (2,2)\}$, and we call the second set \textit{low-value matches}.
%
%
With a similar reasoning to the one for Figure \ref{Model 2 - draw}, for this model we refer to customers at the demand node $2$ and workers at the supply node $1$ as \textit{high-value agents}, and customers at the demand node $1$ and workers at the supply node $2$ as \textit{low-value agents}.

%
\begin{figure}[!htb]
	\centering
	\begin{subfigure}[t]{0.32\textwidth}
		\centering
		\includegraphics[width=\linewidth]{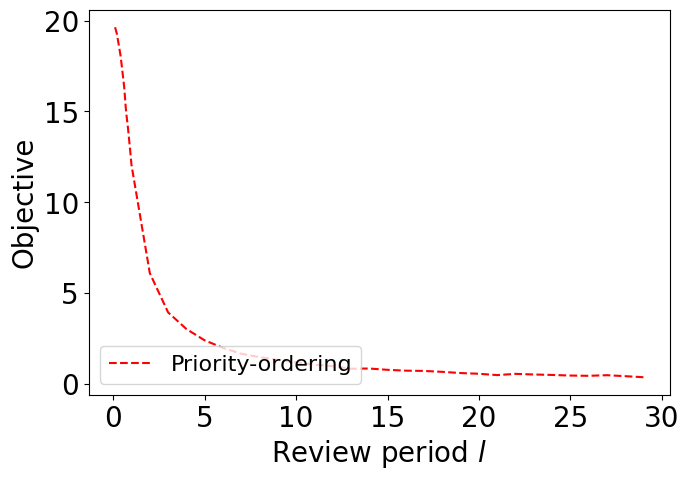}
		\caption{$\boldsymbol{\theta^D} = (0.2, 2)$, $\boldsymbol{\theta^S} = (2,0.2)$.}
		\label{fig:Model1_Results_a}
	\end{subfigure}%
	\begin{subfigure}[t]{0.32\textwidth}
		\centering
		\includegraphics[width=\linewidth]{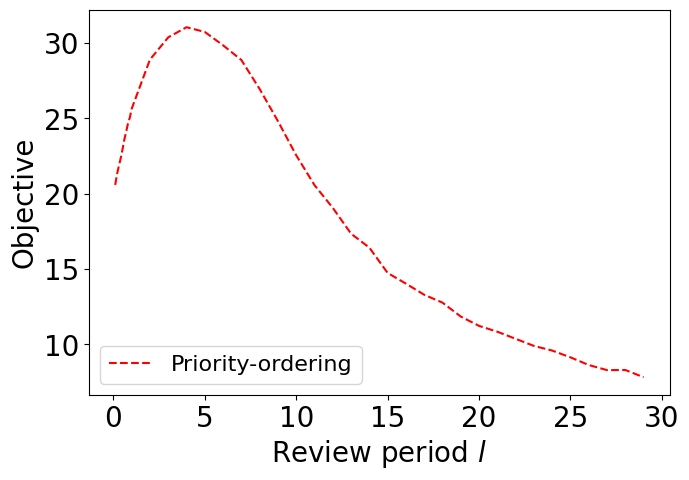}
		\caption{$\boldsymbol{\theta^D} = (0.1, 0.2)$, $\boldsymbol{\theta^S} = (0.2,0.1)$.}
		\label{fig:Model1_Results_b}
	\end{subfigure}
	\begin{subfigure}[t]{0.32\textwidth}
		\centering
		\includegraphics[width=\linewidth]{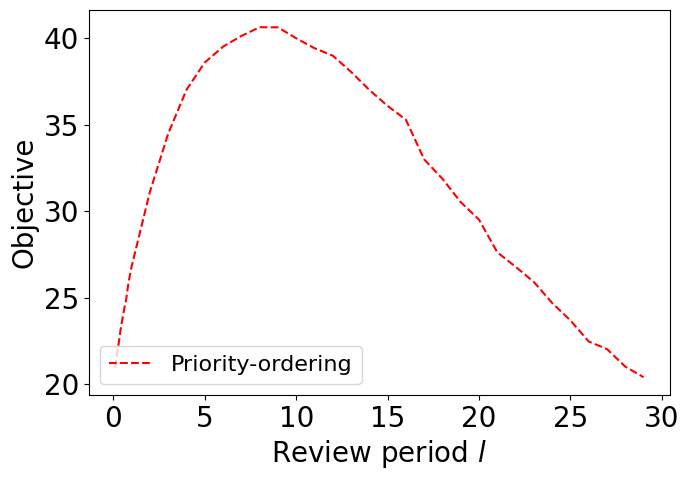}
		\caption{$\boldsymbol{\theta^D} = (0.05, 0.1)$, $\boldsymbol{\theta^S} = (0.1,0.05)$.}
		\label{fig:Model1_Results_c}
	\end{subfigure}
	\caption{Performance of the priority-ordering policy for Model 1.}
	\label{fig:Model1_Results}
\end{figure}
%
We analyze the impact of the length of the review period ($l$) on the performance of the priority-ordering policy for each of the three cases, and present the results in Figure \ref{fig:Model1_Results}.
Below, we summarize the main insights from the simulations.

\textbf{If the high-value \levi{agent}s' patience is too low, we do not allow valuable matches to arise, in which case, the smaller the review period, the better.}
In the first case, the short patience times (average of half a unit of time) and low arrival rates (average of one arrival every ten time units) of the high-value agents make having both high-value demand and supply agents at the end of a review period unlikely.
Therefore, in this case, waiting to gather valuable matches is unlikely to be fruitful, and the policy should instead focus on at least making low-value matches quickly before these agents renege. This is achieved with shorter review periods, as low-value agents are abundant in the system due to their long patience times (average of five time units) and high arrival rates (average of one arrival every time unit). Thus we see in the left panel of Figure \ref{fig:Model1_Results} that the need to quickly match agents dominates in this case, and shorter review periods give better policy performance.\textsuperscript{1}\footnotetext[1]{Note as the review period approaches $0$, the objective corresponds to the fluid solution value when no high-value matches take place and the time horizon is $t = 1000$.}

\textbf{If the high-value agents' patience is long enough to allow high-value matches to arise, then there is an optimal positive value for the review period length.}
In cases two and three, the high-value agents are patient enough to be present at the same time with higher likelihood.
This is seen in the middle and right panels of Figure \ref{fig:Model1_Results}, where, initially, increasing the review period allows gathering enough high-value agents to make more high value matches, thus improving policy performance.
This is different from Figure \ref{Model 2 - draw} because of the presence of supply node 2, which is able to ``steal" high value customers away from demand node 2 and degrade performance if review periods are too short. Thus, in Figure \ref{Model 1 - draw} there is a true trade-off between making matches quickly and waiting to gather valuable matches, due to the presense of low value agents that can divert agents away from future high value matches.

\textbf{The optimal value for the length of the review period depends on the high-value agents' mean patience times.} Figure \ref{fig:Model1_Results} (b) and Figure \ref{fig:Model1_Results} (c) show the bigger the high-value agents' mean patience times, the bigger the length of the optimal review period is. Intuitively, when \levi{agents} have small patience times, we need to decide matches faster, not to lose the valuable matches. On the other hand, when high-value agents are willing to wait in the system, there is no rush to make matches and lose the opportunity to make high-value matches. 
The discrete review length should strike a balance between the arrival and reneging rates.

In summary, both the network and the patience \levi{time} distribution critically impact the review period length. 
\amy{In general,} the patience times need to be long enough relative to the arrival rates and review period length to accumulate enough valuable matches in a review period; \amy{however,}  in network structures that do not allow ``stealing" high-value agents along low-value arcs, shorter review periods \amy{may {\em always}} be better. \amy{An interesting question for future research (that is beyond the scope of this paper) is to characterize the network structures for which an asymptotically optimal policy can make matches immediately when agents arrive, as is the case numerically in Figure~\ref{fig:Model2_Results_zoom}, and looks like the case in
  Figure~\ref{fig:Model1_Results_a} for appropriate conditions on the patience time distributions.  A policy that makes matches immediately when agents arrive is a fully greedy policy. }

\section{A matching policy for a model without holding costs}
\label{sec:PSB}

The matching policies studied until now use the knowledge of the arrival rate parameters and the patience time distributions. This is because both policies use the information provided by the MP \eqref{eq:max_A_w_v3G} which naturally depends on the invariant queue-lengths. The latter uses the information of the arrival rate parameters and the patience time distributions.

However, there is an insensitivity property on the fraction of reneging \levi{agents}.  Table~\ref{table:fractionR1} shows the fraction of reneging \levi{agents}  obtained by the simulations and the corresponding values calculated by the fluid model equations (denoted by $r_S^\star \amy{:= \lambda^D - \min(\lambda^D,\lambda^S)}$ and $r_D^\star \amy{:= \lambda^S - \min(\lambda^D,\lambda^S})$) in a network with
one demand and one supply node that makes every possible match. We observe that for exponential and gamma patience times the fraction of reneging \levi{agents}  remains almost the same, which suggests an insensitivity property, similar to ones that appear in \cite{atar2019large} and \cite{puha2019scheduling}.
\begin{table}[!h]
\begin{center}
\caption{The fraction of reneging population for exponential (right) and gamma (left, rate parameter=0.7) distributions with unit means. The parameters are $T=100$, $\lambda^D=1$, and $n=100$.}
\label{table:fractionR1}
\begin{tabular}{c| c c c c|| c c c c }
& $R^D(T)/A^D(T)$& $r_D^\star$ & $R^S(T)/A^S(T)$  & $r_S^\star$ &  $R^D(T)/A^D(T)$& $r_D^\star$& $R^S(T)/A^S(T)$  & $r_S^\star$   \\ \hline
  $\lambda^S=0.5$  & 0.4942 & 0.5 &1.4e-5       & 0 & 0.4847&0.5 &0&0
  \\
  $\lambda^S=0.9$  & 0.1193 & 0.1 &0.0022  & 0 & 0.0959&0.1 & 0.0221& 0 \\
  $\lambda^S=1.0$  & 0.0586 & 0   &0.0585  & 0 & 0.0130&0   & 0.0581& 0  \\
  $\lambda^S=1.2$  & 0.0066 & 0   &0.173  & 0.2 & 0   &0   & 0.1698& 0.2\\
  $\lambda^S=1.5$  & 4.3e-5      & 0   &0.3314  & 0.5 & 0   &0   & 0.3328& 0.5 \\
\end{tabular}\end{center}
\end{table}

In this section, we sacrifice the generality of the objective function  \eqref{eq:Profit} to take advantage of the insensitivity property observed in Table~\ref{table:fractionR1}.
If we set $c_j^D= c_k^S=0$ for each $j\in \bbJ$ and $k\in \bbK$, then \eqref{eq:Profit} becomes the cumulative matching value of the platform, and hence it is not affected by the number of customers/workers that renege. In other words, it is expected that the cumulative matching value is not affected by the patience time distribution due to the insensitivity property. This will allow us to propose a discrete review matching policy that does not use the arrival rate and patience time distributions information.

In the case without holding costs, the MP takes the following form:
\begin{align}\label{eq:max_A_w_v3}
\begin{aligned}
\max \ &
\sum_{(j,k) \in \arcset}
 \p_{jk} m_{jk} \\
\text{s.t.} \ & \sum_{j\in \demset_k} m_{jk} \le \lambda^S_{k}, \ k \in \bbK,\\
&\sum_{k \in \suppset_j}  m_{jk} \le \lambda_j^D, \ j \in \bbJ, \\
& m_{jk} \ge 0,  \ (j,k) \in \arcset.
\end{aligned}
\end{align}

The number of matches made at each discrete review time point can be decided by solving an optimization problem with an objective to maximize the matching value and constraints that respect the amount of demand and supply available. For any $i\ge 1$, let
$(\mathcal{M}_{ijk}^{\levi{LP},\star}: (j,k) \in \arcset)$ be given by an optimal solution to the following optimization problem:
\begin{align}\label{LP:MatchingPo}
\begin{split}
\max \ & \sum_{(j,k) \in \arcset} \p_{jk} \mathcal{M}_{ijk}^{\levi{LP}} \\
\text{s.t.} \ &  \sum_{j\in \demset_k}   \mathcal{M}_{ijk}^{\levi{LP}} \le Q_k^S(il-), k \in \bbK,\\
&\sum_{k \in \suppset_j}   \mathcal{M}_{ijk}^{\levi{LP}} \le \qd_j(il-), j \in \bbJ,  \\
& \mathcal{M}_{ijk}^{\levi{LP}} \in \mathbb{Z}_+, \mbox{ for all }(j,k) \in \arcset.
\end{split}
\end{align}
The quantity $\mathcal{M}_{ijk}^{\levi{LP},\star}$ is the number of matches between type $j$ demand and type $k$ supply, and the cumulative number of matches for $(j,k) \in \arcset$, and $t\ge 0$, is given by
\begin{equation}\label{eq:cumulativeG1}
M_{jk}^{\levi{LP}}(t) =
\sum_{i=1}^{ \lfloor t/l \rfloor } \mathcal{M}_{ijk}^{\levi{LP},\star}
\end{equation}
where we refer to the above policy as the \emph{LP-based matching policy}.
\levi{It is straightforward to see that the LP-based matching policy satisfies the admissibility conditions given at the end of Section~\ref{sec:Model description}.}

The main difference between \eqref{LP:MatchingPo} and \eqref{eq:max_A_w_v3} is that the right-hand side of the constraints in \eqref{eq:max_A_w_v3} are replaced by the queue-lengths $\qd_j(il-)\in \mathbb{N}$ and $Q_k^S(il-)\in \mathbb{N}$. The proposed policy is  myopic in the sense that the matching decisions made at each review time point are optimal given the available demand and supply but disregard the impact of future arrivals.  The hope is that when discrete review points are well-placed, the aforementioned myopicity will not have too much negative impact, and the resulting total value of matches made can be close to the  optimal value to \eqref{eq:max_A_w_v3}.

The implementation  of the LP-based matching policy does not require the knowledge of the arrival rates $\bld{\lambda}^D$, $\bld{\lambda}^S$ and the patience time distributions but it requires an optimal solution to \eqref{LP:MatchingPo} at each review point.
\amy{From Lemma~\ref{lem:unimodular} below,}
 \eqref{LP:MatchingPo} can effectively be solved as an LP on any sample path, since the queue-lengths are always integer valued.


\begin{lemma}\label{lem:unimodular}
 \amy{The constraint matrix of \eqref{LP:MatchingPo} is totally unimodular, which implies that an optimal extreme point solution is integer valued if the right hand side constraints of \eqref{LP:MatchingPo} are integer valued.}
\end{lemma}

Next we consider the asymptotic behavior of the LP-based policy introduced above in the same asymptotic regime as in Section~\ref{sec:Performance Analysis}. The number of matches in the $n$th system at a discrete review period is given by an optimal solution to
\eqref{LP:MatchingPo} replacing the right-hand sides of the constraints by $\qdn_j(il^n-)$ and $Q_k^{S,n}(il^n-)$, respectively. The next theorem states that the matching policy is asymptotically optimal.   \amy{Note that in contrast to Theorems~\ref{prop:upper}-\ref{Thm:aopriority}, the requirement that hazard rates be bounded is not needed.  This is because when holding costs are zero, there is not a need to analyze the queue-length behavior.}
\begin{theorem} \label{theorem:blindOptimal}
\amy{Assume $c_j^D = c_k^S=0$ for $j\in \bbJ$ and $k \in \bbK$, and let $\bld \mstar$ be an optimal solution to \eqref{eq:max_A_w_v3G}.
\begin{enumerate}
\item[(i)] If $\{ \bld{M}^n \}_{n \in \mathbb{N} }$ is a sequence of matching processes such that as $n \rightarrow \infty$
\begin{equation} \label{eq:ass-M-zero-HC}
    \sup_{0 \leq t \leq T } \left| \frac{M_{jk}^{n}(t)}{n} -
m_{jk} t \right| \rightarrow 0,
\end{equation}
in probability for some $\bld{m} \in \mathbb{M}$ and for all $T\ge 0$, then
\begin{equation} \label{eq:max-achievable-profit-fluid}
\limsup_{t \rightarrow \infty}
\limsup_{n \rightarrow \infty} \frac{1}{nt}
 V_{\bld{M}^n}(t)
  \leq
 \sum_{(j,k) \in \arcset}
 \p_{jk} \mstar_{jk},
\end{equation}
almost surely.
\item[(ii)] The LP-based matching policy given by \eqref{LP:MatchingPo} and  \eqref{eq:cumulativeG1} in the $n$th system is asymptotically optimal in the sense that it satisfies
\[
\lim_{t\rightarrow \infty} \lim_{n\rightarrow \infty}
\frac{1}{nt}\frac{V_{\bld{M}^{\levi{LP},n}}(t)}
 {\sum_{(j,k) \in \arcset}
 \p_{jk} \mstar_{jk}}
 =1,
\]
in probability.
\end{enumerate}}
\end{theorem}

We prove Theorem~\ref{theorem:blindOptimal} using the same three basic steps as Theorem~\ref{proposition:DRconvergence2}. However establishing a lower bound on the matches made during a review period (step ii) requires more effort since we cannot compare directly to the MP optimal solution. We overcome this challenge by leveraging a monotonicity property of the MP \eqref{eq:max_A_w_v3}.

Having shown that the LP-based policy is asymptotically optimal, a question that arises is whether the matches made under this policy also achieve the optimal matching rates. The answer is nuanced in the case when there are multiple optimal solutions to the MP, since the limit of  $\frac{\bld{M}^{\levi{LP},n}(t)}{n}$ may oscillate between them. However, we are able to show that the matching rates must approach the set of optimal MP solutions asymptotically. To this end, for a real vector $\bld{x}$ and a set $A$ in Euclidean space, denote the distance between them by $d\left(\bld{x}, A\right)$, e.g., one could consider
$d\left(\bld{x}, A\right):=\inf_{\bld{z}\in A}\|\bld{x}-\bld{z}\|_2$.
\begin{proposition}\label{theorem:convergenceLPset}
Let $S$ be the set of all optimal solutions of \eqref{eq:max_A_w_v3}. In other words, we have that
$S:=\{\bld{m}^{\star}: \bld{m}^{\star} \mbox{ is an optimal solution to } \eqref{eq:max_A_w_v3} \}$. Then, for each $t>0$ as $n \rightarrow \infty$,
\[
 d\left(\frac{\bld{M}^{\levi{LP},n}(t)}{nt}, S\right) \rightarrow 0,
\]
in probability.
\end{proposition}

A consequence of the last theorem is that an analogous result to Theorem~\ref{proposition:DRconvergence2} holds for the LP-based policy in the special case that \eqref{eq:max_A_w_v3} has a unique optimal solution.
\begin{corollary} \label{corollary:convergenceLP}
Assume \eqref{eq:max_A_w_v3} has a unique optimal solution $\bld{m}^{\star}$. We have that for each $t\ge0$ as $n \rightarrow \infty$,
\[
 \frac{M_{jk}^{\levi{LP},n}(t)}{n} \rightarrow
m_{jk}^{\star} t,
\]
in probability, for all  $(j,k) \in \arcset$.
\end{corollary}

\section{Conclusion}\label{sec:conclusion}
In this paper, we proposed and analyzed a model that takes into account three main features of
service platforms: (i) demand and supply heterogeneity, (ii) the random unfolding of arrivals over time, (iii) the non-Markovian impatience of demand and supply. These features result in a trade-off between making a less good match quickly and waiting for a better match.  The model is too complicated to solve for an optimal matching policy, and so we developed an approximating fluid model, that is accurate in high volume (when demand and supply arrival rates are large).  We used the invariant states of the fluid model to define a matching optimization problem, whose solution gave asymptotically optimal matching rates, that depend on \levi{agent} patience time distributions.  We proposed a discrete review policy to track those asymptotically optimal matching rates, and further established conditions under which a static priority ordering policy also resulted in asymptotically optimal matching rates.  Finally, we observed that when holding costs are zero, there is an insensitivity property that allows us to propose an LP-based matching policy, that also achieves asymptotically optimal matching rates, but does not depend on demand and supply patience time distributions.

\section*{Acknowledgements}
We thank Yuan Zhong for helpful discussion related to Lyapunov functions.  Financial support from the University of Chicago Booth School of Business is gratefully acknowledged.

\bibliographystyle{informs2014}
\bibliography{matching_LDsimAW}

\newpage
\appendix
\counterwithin{lemma}{section}
\counterwithin{proposition}{section}
\noindent{\Large \bf \amy{Matching Impatient and Heterogeneous Demand and Supply:  Technical Appendix}} \\
\noindent{\amy{Angelos\ Aveklouris, Levi\ DeValve, Maximiliano Stock, Amy\ R. Ward}} \\

\amy{In this technical appendix, in Section~\ref{ap:proofs} below, we provide proofs for the results stated in the main body of the manuscript titled: ``Matching Impatient and Heterogeneous Demand and Supply''.  The proofs of these results are in the order in which they appear in the main body, and are separated by the Section in which they appear in the main body.  We also provide, in Section~\ref{sec:markovChain} below, the calculations used to produce the solid lines in Figure~\ref{figure:ExpQueues}.}

\section{Proofs}
\label{ap:proofs}

\subsection{Proofs for Section~\ref{sec:DRMUP}}\label{Proofs priority}
\pfa{Proof of Theorem~\ref{lem:concave}}
We show the result when the hazard rate functions are  increasing. The case of decreasing hazard rate functions shares the same machinery.

Assume that the hazard rate functions are increasing. By \cite[Lemma~1]{puha2019scheduling}, we know that $q_j^{D,\star}(\bld{m})$ and
$q_k^{S,\star}(\bld{m})$ are concave functions of $m_j^D:=\sum_{k\in \suppset_j} m_{jk}\leq \lambda_j^D$ and
$m_k^S:=\sum_{j\in \demset_k} m_{jk}\leq \lambda_k^S$, respectively. Since $q_j^{D,\star}(\bld{m})$ is concave, we have that for any $x_j=\sum_{k \in \suppset_j} x_{jk}\leq \lambda_j^D$,
\begin{equation}\label{eq:conc}
q_j^{D,\star}(m_j^D)-q_j^{D,\star}(x_j) \leq
\frac{dq_j^{D,\star}(x_j) }{dx_j}(m_j^D-x_j).
\end{equation}
Furthermore, from (\ref{eq:InvDeQ}), we have that
\begin{equation*}
\frac{dq_j^{D,\star}(x_j) }{dx_j}=
-\frac{1}{ h^D_j\Big(
(G^D_j)^{-1}\Big( 1-\frac{ x_{j}}{\lambda_j^D}\Big)\Big)}=
\frac{\partial q_j^{D,\star}(\bld x) }{\partial x_{jk}}.
\end{equation*}
Now, by \eqref{eq:conc}, we obtain
\begin{equation*}
q_j^{D,\star}(\bld{m})-q_j^{D,\star}(\bld{x})=
q_j^{D,\star}(m_j^D)-q_j^{D,\star}(x_j) \leq
\frac{dq_j^{D,\star}(x_j) }{dx_j}(m_j^D-x_j) =
\sum_{k \in \bbK}
\frac{\partial q_j^{D,\star}(\bld x) }{\partial x_{jk}}
(m_{jk}^D-x_{jk}),
\end{equation*}
and hence $q_j^{D,\star}(\bld{m})$ is a concave function.  The proof for
 $q_k^{S,\star}(\bld{m})$ follows the same machinery.
\eof \\

\pfa{Proof of Proposition~\ref{proposition:Madmissible2}}
We prove that the queue-lengths are nonnegative. We only show it for the process $\qd_j(\cdot)$; the proof for
$Q_k^S(\cdot)$ follows in the same way. We have that  for $j \in \bbJ$ and
$m \in \{1,\ldots, \lfloor t/l \rfloor \}$, from (\ref{eq:queueDRdemand}),
\begin{equation*}
\begin{split}
\qd_j(ml)
         &= \qd_j((m-1)l)+ \ad_j(ml)-\ad_j((m-1)l)-R^D_j(ml)+R^D_j((m-1)l) - \sum_{k \in \suppset_j} \mathcal{M}_{mjk}^r\\
        &= Q_{j}^D(ml-) - \sum_{k \in \suppset_j}  \mathcal{M}_{mjk}^r .
\end{split}
\end{equation*}
Now, using the inequalities
$\mathcal{M}_{mjk}^r\leq \left \lfloor m_{jk}\frac{Q_{j}^D(ml-)}{\lambda_j^D}\right\rfloor
\le m_{jk}\frac{Q_{j}^D(ml-)}{\lambda_j^D}$ and
$\sum_{k \in \suppset_j} \frac{m_{jk}}{\lambda_j^D}\leq 1$, we have that
\begin{equation*}
\begin{split}
\qd_j(ml) &\geq  Q_{j}^D(ml-) - \sum_{k \in \suppset_j}  \mathcal{M}_{mjk}^r
\ge Q_{j}^D(ml-)-Q_{j}^D(ml-)= 0.
\end{split}
\end{equation*}

\eof

\pfa{Proof of Lemma~\ref{lem:properties}}
Let $m'_{jk} , (j,k) \in \arcset$ denote an optimal extreme point solution. By definition, an optimal extreme point solution cannot be written as a convex combination of any two distinct feasible solutions.
	
	{\bf Property 1.}  Assume the induced graph has a cycle, and renumber the nodes involved in the cycle so that it can be represented by the set of edges $C = \{(1,1),(1,2),\dots,(i,i),(i,i+1),\dots,(r,r),(r,1)\}$, where the length of the cycle is assumed to be $2r$ for some $r$ without loss of generality because the underlying graph is bipartite. Then, we will show that solution $m'_{jk} , (j,k) \in \arcset$ can be written as a convex combination of two distinct feasible solutions, contradicting the fact that it is an extreme point.
	
	First, for edges in the cycle, let $m^1_{ii} = m'_{ii} + \epsilon, i = 1,\dots, r$ and $m^1_{ii+1} = m'_{ii+1} - \epsilon, i = 1,\dots,r-1$ and $m^1_{r1} = m'_{r1} - \epsilon$ for $\epsilon>0$, while for edges not in the cycle, i.e., $(j,k) \notin C$, let $m^1_{jk} = m'_{jk}$. For $\epsilon \le \min_{(j,k) \in C}\{m'_{jk}\}$, all variables remain nonnegative, and we next argue that this solution remains feasible for the other constraints as well. Consider a demand node $i=1,\dots,r-1$ involved in the cycle $C$. We have $m^1_{ii} + m^1_{ii+1} = m'_{ii} +\epsilon + m^1_{ii+1}-\epsilon = m'_{ii} + m^1_{ii+1}$, so that the total contribution of edges $(i,i)$ and $(i,i+1)$ to the flow entering node $i$ is unchanged, hence the associated demand constraint for $i$ is satisfied (since $\bld{m}'$ was assumed feasible). Similarly, for demand node $r$ in the cycle $C$, we have $m^1_{rr} + m^1_{r1} = m'_{rr} +\epsilon + m^1_{r1}-\epsilon = m'_{rr} + m^1_{r1}$ so that the associated demand constraint for $r$ is satisfied. An identical argument demonstrates that the supply constraints are also satisfied.
	
	Next, using the same $\epsilon$, for edges in the cycle let $m^2_{ii} = m'_{ii} - \epsilon, i = 1,\dots, r$ and $m^2_{ii+1} = m'_{ii+1} + \epsilon, i = 1,\dots,r-1$ and $m^2_{r1} = m'_{r1} + \epsilon$ for $\epsilon>0$, while for edges not in the cycle, i.e., $(j,k) \notin C$, let $m^2_{jk} = m'_{jk}$. Again, this solution is feasible, and we have $\bld{m}' = 1/2\bld{m}^1 + 1/2\bld{m}^2$, hence $\bld{m}'$ is not an extreme point of \eqref{eq:max_A_w_v3G}.
	
	{\bf Property 2.} Assume that in the induced graph, a given tree has more than one node with a slack constraint. Then, pick any two of these slack nodes, consider the path between them, and renumber the nodes involved in the path so that it can be represented by the set of edges $P = \{(1,1),(1,2),\dots,(i,i),(i,i+1),\dots,(r,r)\}$. (This representation implicitly assumes that the path begins on a demand node and ends on a supply node, but the following argument applies for the analogous path representation for any combination of begining/ending node classification). By the assumption that both endpoints of the path have slack constraints, we have $\lambda^D_i > \sum_k m'_{ik}$ and $\lambda^S_r > \sum_j m'_{jr}$. Then, we will show that solution $m'_{jk} , (j,k) \in \arcset$ can be written as a convex combination of two distinct feasible solutions, contradicting the fact that it is an extreme point.
	
	First, for edges in the path, let $m^1_{ii} = m'_{ii} + \epsilon, i = 1,\dots, r$ and $m^1_{ii+1} = m'_{ii+1} - \epsilon, i = 1,\dots,r-1$ for $\epsilon>0$, while for edges not in the path, i.e., $(j,k) \notin P$, let $m^1_{jk} = m'_{jk}$. For $$\epsilon \le \min\left(\min_{(j,k) \in P}\{m'_{jk}\},\lambda^D_i - \sum_k m'_{ik}, \lambda^S_r - \sum_j m'_{jr}\right),$$ all variables remain nonnegative, we do not violate the constraints on the endpoint nodes of the paths, and feasibility of this solution for the remaining constraints along the path follows from an identical argument to that posed in the proof of Property 1.
	
	Next, using the same $\epsilon$, for edges in the path let $m^2_{ii} = m'_{ii} - \epsilon, i = 1,\dots, r$ and $m^2_{ii+1} = m'_{ii+1} + \epsilon, i = 1,\dots,r-1$, while for edges not in the path, i.e., $(j,k) \notin P$, let $m^2_{jk} = m'_{jk}$. Again, this solution is feasible, and we have $\bld{m}' = 1/2\bld{m}^1 + 1/2\bld{m}^2$, hence $\bld{m}'$ is not an extreme point of \eqref{eq:max_A_w_v3G}.
	
\eof \\

\levi{In order to prove Lemma 2 from Section \ref{sec:priotiry}, we first establish the following Lemma.}

\begin{lemma}\label{lem:priority_set_cap}
	For Algorithm~\ref{alg:priority_sets}, throughout iteration $h$ of the outer while loop beginning in Step \ref{alg:priority_sets_out_whl}, we maintain the following for each edge $(j,k) \in E$: either $d_j = \lambda_j^D - \sum_{k:(j,k) \in \calQ_{h-1}(\bld{m}^\star)}\mstar_{jk}$ and $s_k = \lambda_k^S - \sum_{j:(j,k) \in \calQ_{h-1}(\bld{m}^\star)}\mstar_{jk}$, or $(j,k) \notin C$.
\end{lemma}
\pfa{Proof}
In the proof we drop the dependence on $\bld{m}^\star$ for the set definitions since this dependence is clear from the context.
First we claim that at the beginning of iteration $h$, for all $(j,k) \in E$ we have $d_j = \lambda_j^D - \sum_{k:(j,k) \in \calQ_{h-1}}\mstar_{jk}$ and $s_k = \lambda_k^S - \sum_{j:(j,k) \in \calQ_{h-1}}\mstar_{jk}$. This holds because the only time we change $d_j$ or $s_k$ is in Step \ref{alg:priority_sets_updt_cap}, and each change to $d_j$ (resp. $s_k$) corresponds to a decrease of $\mstar_{jk'}$ (resp. $\mstar_{j'k}$) for some $(j,k') \in \calQ_{h-1}$ (resp.  $(j',k) \in \calQ_{h-1}$). Then, during iteration $h$, as long as $(j,k) \in C$, these values of $d_j$ and $s_k$ don't change, since Step \ref{alg:priority_sets_remove_nbrs} implies that if $(j,k) \in C$, then no edge $(j',k')$ with either $j'=j$ or $k'=k$ has satisfied the condition in Step \ref{alg:priority_sets_if_tight} yet, so neither $d_j$ nor $s_k$ have been updated in Step \ref{alg:priority_sets_updt_cap}.
\eof \\

\pfa{Proof of Lemma~\ref{prop:Alg1}}
In the proof we drop the dependence on $\bld{m}^\star$ for the set definitions since this dependence is clear from the context.
	The inner while loop beginning in Step \ref{alg:priority_sets_in_whl} always removes at least the element $(j,k)$ chosen in Step \ref{alg:priority_sets_choose_jk}, so this while loop terminates. The outer while loop beginning in Step \ref{alg:priority_sets_out_whl} removes the elements of $\calP_h \subseteq E$ in Step \ref{alg:priority_sets_remove_Ph}, so as long as $|\calP_h|>0$, this while loop terminates also. Step \ref{alg:priority_sets_add_Ph} adds the element $(j,k) \in C$ to $\calP_j$ as long as the condition in Step \ref{alg:priority_sets_if_tight} is met. Thus, we must show that at least one element in $C$ satisfies the condition in Step \ref{alg:priority_sets_if_tight} at the start of each iteration of the outer while loop (which are indexed by the counter $h$).
	
	To this end, consider the situation at the start of iteration $h$ of the outer while loop where $|E|>0$.  Then, we have $E = \{(j,k) \in \arcset:\mstar_{jk}>0\} \setminus \calQ_{h-1}$, and at the start of iteration $h$ we have $C=E$ and for all $(j,k) \in E$ we have $d_j = \lambda_j^D - \sum_{k:(j,k) \in \calQ_{h-1}}\mstar_{jk}$, and $s_k = \lambda_k^S - \sum_{j:(j,k) \in \calQ_{h-1}}\mstar_{jk}$ by Lemma \ref{lem:priority_set_cap}. Next, let
$E_0 = \{(j,k) \in \arcset: \mstar_{jk}>0\}$ denote the induced graph of $\bld{\mstar}$, which is an extreme point of \eqref{eq:max_A_w_v3G}, and thus Lemma \ref{lem:properties} guarantees that $G$ has the following two properties: i) it is a collection of trees (equivalently contains no cycles), and ii) each distinct tree in $G$ has at most one node with a slack constraint. Then, observe that removing the edges $\calQ_{h-1}$ from $E_0 $ to create the new graph
$F = \{(j,k) \in \arcset: \mstar_{jk}>0\}\setminus \calQ_{h-1}$ will preserve both these properties: for i) $F$ must still contain no cycles since we only removed edges (and hence is still a collection of trees), and for ii) if any tree in $F$ has more than one slack node, then these two nodes were also in the same tree in $E_0 $, contradicting Lemma \ref{lem:properties}. Then, note that $F$ must contain a tree with at least two nodes, since $E$ is its edge set and $|E|>0$, denote the tree by $T$. Since $T$ is a tree with at least two nodes, it contains at least two nodes with degree one. By property ii) of $F$, at most one of these degree one nodes has a slack constraint, so therefore at least one of them must have a tight constraint. Without loss of generality, assume this is node  $j \in \bbJ$ (if it is in $\bbK$ an identical argument applies) and let $(j,k)$ denote the edge connected to node $j$ in graph $F$. Since $j$ is tight we have $\sum_{k'} \mstar_{jk'} = \lambda_j^D$, and since $j$ is degree one in $F$ we have $\sum_{k':k' \ne k, (j,k') \notin \calQ_{h-1}} \mstar_{jk} = 0$, which together imply that $\mstar_{jk} = \lambda_j^D - \sum_{k':(j,k') \in \calQ_{h-1}}\mstar_{jk'} = d_j$, and thus $(j,k)$ meets the condition in Step \ref{alg:priority_sets_if_tight}.
	
	Finally, for the running time, each iteration of the outer while loop performs a number of operations linear in $|C|\le |E_0 |$ (since it has to consider each element of $C$, either to check the condition in Step \ref{alg:priority_sets_if_tight}, or to remove the element from $C$ in step \ref{alg:priority_sets_remove_nbrs}). The outer while loop removes at least one element from $E$ in each iteration, so the number of iterations is less than the initial size of $E$, which is $|E_0 |$.
\eof

\pfa{Proof of Proposition~\ref{prop:property extreme}}
	We proceed by induction on $h$ up to set $H$, then argue for $\calP_{H+1}$ separately. For $h=0$, consider an arc $(j,k) \in \mathcal{P}_0( \bld{m}^\star )$, i.e., in iteration 0 of the outer while loop, $(j,k)$ was chosen in Step \ref{alg:priority_sets_choose_jk} and satisfied the condition in Step \ref{alg:priority_sets_if_tight}. Without loss of generality, assume the condition satisfied in Step \ref{alg:priority_sets_if_tight} was $\mstar_{jk}= d_j$ (an identical argument holds if $\mstar_{jk}= \lambda_k^S$). By Lemma \ref{lem:priority_set_cap}, since $(j,k) \in C$ when it was chosen, we have $\mstar_{jk}= d_j = \lambda_j^D - \sum_{k:(j,k) \in \calQ_{-1}}\mstar_{jk} = \lambda_j^D$. Further, by feasibility of $\bld{m}^\star$ for \eqref{eq:max_A_w_v3G}, we have $m^\star_{jk} \le \lambda_k^S$, which implies that $\lambda_j^D \le \lambda_k^S$. Thus, by the definition of $y^p_{jk}(\bld{m}^\star)$ we have
	\begin{align*}
	y^p_{jk}(\bld{m}^\star) =
	\min\left( \lambda_j^D, \lambda_k^S \right) = \lambda_j^D = m^\star_{jk}.
	\end{align*}

Now assume that $y^p_{jk}(\bld{m}^\star) = m^\star_{jk}$ for all $(j,k) \in \mathcal{Q}_{h-1}(\bld{m}^\star)$ and consider $h$. For $(j,k) \in \mathcal{P}_h(\bld{m}^\star)$ we know that in iteration $h$ of the outer while loop, $(j,k)$ was chosen in Step \ref{alg:priority_sets_choose_jk} and satisfied the condition in Step \ref{alg:priority_sets_if_tight}. Without loss of generality, assume the condition satisfied in Step \ref{alg:priority_sets_if_tight} was $\mstar_{jk}= d_j$ (an identical argument holds if $\mstar_{jk}= \lambda_k^S$). By Lemma \ref{lem:priority_set_cap}, since $(j,k) \in C$ when it was chosen, we have $\mstar_{jk}= d_j = \lambda_j^D - \sum_{k:(j,k) \in \calQ_{h-1}}\mstar_{jk} $. Further, by feasibility of $\bld{m}^\star$ to \eqref{eq:max_A_w_v3G}, we have
$\sum_{j'\in \bbJ} m^\star_{j'k} \le \lambda_k^S$, and since $m^\star_{j'k} \ge 0$ for all $j'$, we also have $m^\star_{jk} + \sum_{\{j':(j',k)\in \mathcal{Q}_{h-1}(\bld{m}^\star)\}}
m^\star_{j'k} \le \lambda_k^S$, or equivalently $m^\star_{jk} \le \lambda_k^S - \sum_{\{j':(j',k)\in \mathcal{Q}_{h-1}(\bld{m}^\star)\}}m^\star_{j'k}$. From this it follows that $\lambda_j^D - \sum_{\{k':(j,k') \in \mathcal{Q}_{h-1}(\bld{m}^\star)\}}m'_{jk'} \le  \lambda_k^S - \sum_{\{j':(j',k)\in \mathcal{Q}_{h-1}(\bld{m}^\star)\}}m^\star_{j'k}$ and by the definition of $y^p_{jk}(\bld{m}^\star)$ and the induction hypothesis, we have
\begin{align*}
y^p_{jk}(\bld{m}^\star) &=
\min\left( \lambda_j^D - \sum_{\{k':(j,k') \in \mathcal{Q}_{h-1}(\bld{m}^\star)\}}m^\star_{jk'}, \lambda_k^S - \sum_{\{j':(j',k)\in \mathcal{Q}_{h-1}(\bld{m}')\}}m^\star_{j'k} \right) \\
&= \lambda_j^D - \sum_{\{k':(j,k') \in \mathcal{Q}_{h-1}(\bld{m}^\star)\}}m^\star_{jk'} = m^\star_{jk}.
\end{align*}

Finally, for $h = H+1$, for $(j,k) \in \mathcal{P}_{H+1}(\bld{m}^\star)$, we must have either $\sum_{\{k':(j,k')\in \mathcal{Q}_{h-1}(\bld{m}^\star)\}} m^\star_{jk'} = \lambda_j^D$ or $\sum_{\{j':(j',k)\in \mathcal{Q}_{h-1}(\bld{m}^\star)\}}  m^\star_{jk} = \lambda_k^S$, since otherwise $\bld{m}^\star$ was not optimal (since the objective is strictly increasing in $\bld{m}^\star$). Thus, by definition we have $y^p_{jk}(\bld{m}^\star) = 0 = m^\star_{jk}$.
\eof \\

\subsection{Proofs for Section~\ref{sec:fluid model}}
\label{proofs fluid model}
\pfa{Proof of Proposition~\ref{prop:Chofinvariant}}
\amy{
We first show that a tuple satisfying (\ref{eq:InvDe})-(\ref{eq:InvSuI}) is an invariant state, as given in Definition~\ref{Definition:IM} (the forward direction).  We second show that an invariant state must satisfy (\ref{eq:InvDe})-(\ref{eq:InvSuI}) (the converse direction).

{\bf Proof of Forward Direction.} Let  $(\bld{q}^{D,\star}, \bld{q}^{S,\star}, \bld{\eta}^{D,\star}, \bld{\eta}^{S,\star})$ be defined by (\ref{eq:InvDe})-(\ref{eq:InvSuI})
for a given matching rate matrix $\bld m \in \mathbb{M}$.
 Let
$(\ovl{\bld{Q}}^D, \ovl{\bld{Q}}^S,
\ovl{\bld{\eta}}^{D}, \ovl{\bld{\eta}}^{S})$
be the constant function given by
\begin{equation}\label{eq:constant}
(\ovl{\bld{Q}}^D(t), \ovl{\bld{Q}}^S(t),
\ovl{\bld{\eta}}^{D}(t), \ovl{\bld{\eta}}^{S}(t))=
(\bld{q}^{D,\star}, \bld{q}^{S,\star}, \bld{\eta}^{D,\star}, \bld{\eta}^{S,\star}), \mbox{ for all } t \geq 0.
\end{equation}
Then,
$(\ovl{\bld{Q}}^D,\ovl{\bld{Q}}^S,
\ovl{\bld{\eta}}^D, \ovl{\bld{\eta}}^S )
\in  \mathcal{C}(\overline{\mathbb{Y}})$, also noting that for  any $\bld m \in \mathbb{M}$,
\[
0 \leq q_j^{D,\star}(\bld m) \leq \lran{ 1, \eta_j^{D,\star} }, j \in \bbJ, \mbox{ and } 0 \leq q_k^{S,\star}(\bld m) \leq \lran{ 1, \eta_k^{S,\star} }, k \in \bbK.
\]

We must show that $(\ovl{\bld{Q}}^D, \ovl{\bld{Q}}^S,
\ovl{\bld{\eta}}^{D}, \ovl{\bld{\eta}}^{S})$ is a fluid model solution (see Definition~\ref{def:fluid-model}) for  $\left( \bld{\lambda}^D, \bld{\lambda}^S \right)$.
From (\ref{eq:InvDe}) and (\ref{eq:InvSu}),
\[
    \lran{ h^D_j, \eta^{D,\star}_j } = \lambda_j^D < \infty, \mbox{ for } j \in \bbJ, \mbox{ and }
    \lran{ h^S_k, \eta^{S,\star}_k } = \lambda_k^S < \infty, \mbox{ for } k \in \bbK,
\]
and so \eqref{in:abrates} holds.  Furthermore, (\ref{eq:InvDe}) and (\ref{eq:InvSu}) also imply (\ref{eq:fluid-no-atom-initial-condition}) holds. Next, for any bounded function $f\in \mathcal{C}(\R_+)$, again from (\ref{eq:InvDe}) and (\ref{eq:InvSu}),

\begin{eqnarray*}
 \int_0^{H^D_j} f(x+t)
\frac{1-G^D_j(x+t)}{1-G^D_j(x)} \eta^{D,\star}_j(dx) + \lambda_j^D \int_0^t f(t-u)(1-G^D_j(t-u))du
&  = &   \lran{ f,  \eta_j^{D,\star} }, \mbox{ for } j \in \bbJ,
\end{eqnarray*}
and
\begin{eqnarray*}
\int_0^{H^S_k} f(x+t)
\frac{1-G^S_k(x+t)}{1-G^S_k(x)} \eta^{S,\star}_k(dx) + \lambda_k^S \int_0^t f(t-u)(1-G^S_k(t-u))du
& = & \lran{ f,  \eta_k^{S,\star} }, \mbox{ for } k \in \bbK.
\end{eqnarray*}
Hence \eqref{eq:FetaD} and \eqref{eq:FetaS} hold.

We next show that  \eqref{eq:QFD} and \eqref{eq:QFS} hold,
 with  $\ovl{\bld{R}}^D$ and $\ovl{\bld{R}}^S$  given by  \eqref{eq:AlRD} and \eqref{eq:AlRS}.
 Let
 \[
\chi_j^D  =  \left( G_j^D \right)^{-1} \left( 1- \frac{\sum_{k \in \suppset_j} m_{jk}}{\lambda_j^D} \right)
 \]
 if $\sum_{k \in \suppset_j} m_{jk}>0$ and $H_j^D$ otherwise, for each $j \in \bbJ$.
 Similarly, let
 \[
\chi_k^S =  \left( G_k^S \right)^{-1} \left( 1- \frac{\sum_{j\in \demset_k}  m_{jk}}{\lambda_k^S} \right)
 \]
 if $\sum_{j\in \demset_k}  m_{jk}>0$ and $H_k^S$ otherwise, for each $k \in \bbK$.
 Then, from  \eqref{eq:AlRD} - \eqref{eq:AlRS} and (\ref{eq:InvDe})-(\ref{eq:InvSuI}),
 \begin{eqnarray*}
 \overline{R}_j^D(t) & = &
 t \int_0^{H_j^D} 1\left\{ G_{e,j}^D(x) < G_{e,j}^D \left( \left( G_j^D \right)^{-1} \left( 1- \frac{\sum_{k \in \suppset_j}m_{jk}}{\lambda_j^D} \right) \right) \right\} g_j^D(x) dx \\
 & = &  t \int_0^{\chi_j^D} \lambda_j^D g_j^D(x) dx \\
 & = & t \lambda_j^D G_j^D \left( \chi_j^D \right),
 \end{eqnarray*}
 for $j \in \bbJ$, and, similarly,
\begin{eqnarray*}
 \overline{R}_k^S(t)
 & = &  t  \lambda_k^S G_k^S \left( \chi_k^S \right),
 \end{eqnarray*}
 for $k \in \bbK$.
 Substituting for $\chi_j^D$ and $\chi_k^S$ in the above yields
 \[
 \overline{R}_j^D(t) = t \lambda_j^D \left( 1- \frac{\sum_{k \in \suppset_j} m_{jk}}{\lambda_j^D} \right), j \in \bbJ, \mbox{ and }
 \overline{R}_k^S(t) = t \lambda_k^S \left( 1- \frac{\sum_{j\in \demset_k}  m_{jk}}{\lambda_k^S} \right), k \in \bbK,
 \]
 and so \eqref{eq:QFD} and \eqref{eq:QFS} hold.

{\bf Proof of Converse Direction.}
Suppose that $(\bld{q}^{D}, \bld{q}^{S}, \overline{\bld{\eta}}^{D}, \overline{\bld{\eta}}^{S}) \in \mathcal{I}_{\bld{\lambda}}$.
 Let
$(\ovl{\bld{Q}}^D, \ovl{\bld{Q}}^S,
\ovl{\bld{\eta}}^{D}, \ovl{\bld{\eta}}^{S})$
be the constant function such that
\begin{equation}\label{eq:constant}
(\ovl{\bld{Q}}^D(t), \ovl{\bld{Q}}^S(t),
\ovl{\bld{\eta}}^{D}(t), \ovl{\bld{\eta}}^{S}(t))=
(\bld{q}^{D}, \bld{q}^{S}, \overline{\bld{\eta}}^{D}, \overline{\bld{\eta}}^{S}), \mbox{ for all } t \geq 0.
\end{equation}
By assumption, $(\ovl{\bld{Q}}^D, \ovl{\bld{Q}}^S,
\ovl{\bld{\eta}}^{D}, \ovl{\bld{\eta}}^{S})$ is a fluid model solution.  We must show that (\ref{eq:InvDe})-(\ref{eq:InvSuI}) hold.

From \eqref{eq:FetaD} and \eqref{eq:FetaS}, for any bounded function $f\in \mathcal{C}(\R_+)$
\begin{equation}\label{eq:FetaDinvariant}
\lran{ f,  \ovl{\eta}_j^{D} } =
\int_0^{H^D_j} f(x+t)
\frac{1-G^D_j(x+t)}{1-G^D_j(x)} \ovl{\eta}^{D}_j(dx) +
 \lambda_j^D \int_0^t f(t-u)(1-G^D_j(t-u))
du,
\end{equation}
for $j \in \bbJ$, and
\begin{equation}\label{eq:FetaSinvariant}
\lran{ f,  \ovl{\eta}_k^{S} } =
\int_0^{H^S_k} f(x+t)
\frac{1-G^S_k(x+t)}{1-G^S_k(x)} \ovl{\eta}^{S}_k(dx)
 +
\lambda_k^S \int_0^t f(t-u)(1-G^S_k(t-u))
du,
\end{equation}
for $k \in \bbK$.  For each $j \in \bbJ$ and $k \in \bbK$, as $t \rightarrow \infty$, the first integrals in (\ref{eq:FetaDinvariant}) and (\ref{eq:FetaSinvariant}) converge to zero by dominated convergence and the second to
$\lambda_j^D \int_0^\infty \left( 1-G_j^{D}(u) \right) f(u) du$ and $\lambda_k^S \int_0^\infty \left( 1-G_k^{S}(u) \right) f(u) du$.  Because $f$ was arbitrary, $\overline{\bld \eta}^D = \bld \eta^{D,\star}$ and $\overline{\bld \eta}^S = \bld \eta^{S,\star}$, and so (\ref{eq:InvDe}) and (\ref{eq:InvSu}) hold.

We next show that (\ref{eq:InvDeQ}) - (\ref{eq:InvSuI}) hold.
From \eqref{eq:QFD} - \eqref{eq:QFS},
\begin{equation} \label{eq:tmp1}
    \left( \lambda_j^D - \sum_{k \in \suppset_j} m_{jk} \right) t = \overline{R}_j^D(t), \mbox{ for } j \in \bbJ,
\end{equation}
and
\[
    \left( \lambda_k^S - \sum_{j\in \demset_k}  m_{jk} \right) t = \overline{R}_k^S(t), \mbox{ for } k \in \bbK.
\]
Define $\chi_j^D \in [0,H_j^D]$ and $\chi_k^S \in [0,H_k^S]$ to be the unique solutions to
\begin{equation} \label{eq:tmp3}
q_j^{D} = \int_0^{\chi_j^D} \eta_j^{D,\star} (u) du \mbox{ and } q_k^{S} = \int_0^{\chi_k^S} \eta_k^{S,\star} (u) du,
\end{equation}
noting the slight abuse of notation as $\chi_j^D$ and $\chi_k^S$ were defined differently in this proof, when proving the forward direction.
From \eqref{eq:AlRD} - \eqref{eq:AlRS} and (\ref{eq:InvDe}) - (\ref{eq:InvSu}),
\begin{equation} \label{eq:tmp2}
\overline{R}_j^D(t) = \int_0^t \lambda_j^D \int_0^{\chi_j^D} g_j^D(x) dx du \mbox{ and } \overline{R}_k^S(t) = \int_0^t \lambda_k^S \int_0^{\chi_k^S} g_k^S(x) dx du.
\end{equation}

Combining (\ref{eq:tmp1}) and (\ref{eq:tmp2})  shows that for each $j \in \bbJ$ and $k \in \bbK$
\[
\left( \lambda_j^D - \sum_{k \in \suppset_j} m_{jk} \right) = \lambda_j^D G_j^D \left( \chi_j^D \right),  \mbox{ and } \left( \lambda_k^S - \sum_{j\in \demset_k}  m_{jk} \right) =  \lambda_k^S G_k^S \left( \chi_k^S \right),
\]
or, equivalently,
\[
\chi_j^D =    \left( G_j^D \right)^{-1} \left( 1 - \frac{\sum_{k \in \suppset_j} m_{jk}}{\lambda_j^D} \right) \mbox{ if } \sum_{k \in \suppset_j} m_{jk} >0,
\]
and $\chi_j^D = H_j^D$ otherwise, and
\[
\chi_k^S = \left( G_k^S \right)^{-1} \left( 1- \frac{\sum_{j\in \demset_k}  m_{jk}}{\lambda_k^S } \right) \mbox{ if } \sum_{j\in \demset_k}  m_{jk}>0,
\]
and $\chi_k^S = H_k^S$ otherwise.
Finally, substituting for $\eta_j^{D,\star}$ and $\eta_k^{S,\star}$ in (\ref{eq:tmp3}), and also using the above, yields
\[
    q_j^{D} = \int_0^{\chi_j^D} \lambda_j^D \left( 1-G_j^D(x) \right) dx = \frac{\lambda_j^D}{\theta_j} G_{e,j}^D\left( \chi_j^D \right) = q_j^{D,\star}(\bld m), \mbox{ for } j \in \bbJ,
\]
and
\[
    q_k^{S} = \int_0^{\chi_k^S} \lambda_k \left( 1-G_k^S(x) \right) dx = \frac{\lambda_k}{\theta_k} G_{e,k}^D\left( \chi_k^S \right) = q_k^{S,\star}(\bld m), \mbox{ for } k \in \bbK,
\]
which  establishes that (\ref{eq:InvDeQ}) and (\ref{eq:InvSuI}) hold, and completes the proof.
}
\eof \\

\pfa{Proof of Theorem~\ref{Th:steadySMB}}
\amy{Recall equations \eqref{eq:FetaD} and \eqref{eq:FetaS}, repeated here for convenience.
For any bounded function $f\in \mathcal{C}(\R_+)$,
\[
\lran{ f,  \ovl{\eta}_j^{D} } =
\int_0^{H^D_j} f(x+t)
\frac{1-G^D_j(x+t)}{1-G^D_j(x)} \ovl{\eta}^{D}_j(dx) +
 \lambda_j^D \int_0^t f(t-u)(1-G^D_j(t-u))
du,
\]
for $j \in \bbJ$, and
\[
\lran{ f,  \ovl{\eta}_k^{S} } =
\int_0^{H^S_k} f(x+t)
\frac{1-G^S_k(x+t)}{1-G^S_k(x)} \ovl{\eta}^{S}_k(dx)
 +
\lambda_k^S \int_0^t f(t-u)(1-G^S_k(t-u))
du,
\]
for $k \in \bbK$.  For each $j \in \bbJ$ and $k \in \bbK$, as $t \rightarrow \infty$, the first integrals in the above converge to zero by dominated convergence and the second to
$\lambda_j^D \int_0^\infty \left( 1-G_j^{D}(u) \right) f(u) du$ and $\lambda_k^S \int_0^\infty \left( 1-G_k^{S}(u) \right) f(u) du$.  Because $f$ was arbitrary, we conclude that the measures  $\left( \ovl{\bld{\eta}}^{D}(t), \ovl{\bld{\eta}}^{S}(t)\right)$ weakly converge to $\left(\bld{\eta}^{D,\star}, \bld{\eta}^{S,\star}\right)$; i.e., that
\begin{equation} \label{eq:eta-long-term-behavior}
\lim_{t \rightarrow \infty} \left( \ovl{\bld{\eta}}^{D}(t), \ovl{\bld{\eta}}^{S}(t)\right) = \left(\bld{\eta}^{D,\star}, \bld{\eta}^{S,\star}\right)
\end{equation}

We next show that
\begin{equation} \label{eq:queue-long-term-behavior}
\lim_{t\to \infty} \left(\ovl{\bld{Q}}^D(t), \ovl{\bld{Q}}^S(t)
\right)
=
\left(\bld{q}^{D,\star}(\bld m), \bld{q}^{S,\star}(\bld m) \right).
\end{equation}
To do this, we fix $j \in \bbJ$ and show
\begin{equation} \label{eq:fixed-j}
\lim_{t \to\infty} Q_j^D(t) = q^{D,\star}_j(\bld m).
\end{equation}
The same arguments can then be used to show
\[
\lim_{t \to\infty} Q_j^D(t) = q^{D,\star}_j(\bld m) \mbox{ for each } j \in \bbJ, \mbox{ and } \lim_{t \to\infty} Q_k^S(t) = q^{S,\star}_k(\bld m) \mbox{ for each } k \in \bbK,
\]
from which we can conclude that (\ref{eq:queue-long-term-behavior}) holds, and complete the proof (because together
\eqref{eq:eta-long-term-behavior} and \eqref{eq:queue-long-term-behavior} establish \eqref{thrm:fluid-long-time-behavior}).  Below we separate the arguments to show (\ref{eq:fixed-j}) for the three cases: $\sum_{k \in \suppset_j} m_{jk} = \lambda_j^D$ (match everything), $\sum_{k \in \suppset_j} m_{jk} \in \left(0,\lambda_j^D \right)$ (partially match), and $\sum_{k \in \suppset_j} m_{jk} =0$ (match nothing). \\

\noindent \textit{ Case 1 (match everything, $\sum_{k \in \suppset_j} m_{jk} = \lambda_j^D$):} If $\lambda_j^D - \sum_{k \in \suppset_j} m_{jk} = 0$, then $q_j^{D,\star}(\bld m) = 0$  from \eqref{eq:InvDeQ}.
Next note that from \eqref{eq:AlRD} and \eqref{eq:QFD}
\begin{eqnarray}
\frac{d}{dt} \overline{Q}_j^D(t) & = & \left( \lambda_j^D - \sum_{k \in \suppset_j} m_{jk} \right) - \int_0^{H^D_j} h^D_j(x)
\ind{\ovl{\eta}_j^D(t)[0,x]< \ovl{Q}^D_j(t)}
\ovl{\eta}_j^D(t)(dx) \nonumber \\
& = & - \int_0^{H^D_j} h^D_j(x)
\ind{\ovl{\eta}_j^D(t)[0,x]< \ovl{Q}^D_j(t)}
\ovl{\eta}_j^D(t)(dx), \label{eq:Q-deriv}
\end{eqnarray}
for $t > 0$.
 Consider the Lyapunov function $V(q) = q$.
 To establish \eqref{eq:fixed-j},  this is sufficient to show
\begin{equation} \label{eq:case1}
 \frac{d}{dq}V(q) \left|_{q = \ovl{Q}_j^D(t)} \right. \times \frac{d}{dt} \overline{Q}_j^D(t) =  \frac{d}{dt} \overline{Q}_j^D(t)< 0,
\end{equation}
for any $t >0$ such that $\overline{Q}_j^D(t)  >0$.  Note that for any $t >0$,  \eqref{eq:fluid-no-atom-initial-condition}, (\ref{eq:FetaD}), and the assumption that $G_j^D$ is the distribution function of an absolutely continuous random variable, together imply that
$\ovl{\eta}_j^D(t)$ has no atoms and that $\ovl{\eta}_j^D(t)(dx) >0$ for $x \in (0,t)$.  Then, whenever $\overline{Q}_j^D(t) >0$,
\[
\int_0^{H^D_j} h^D_j(x)
\ind{\ovl{\eta}_j^D(t)[0,x]< \ovl{Q}^D_j(t)}
\ovl{\eta}_j^D(t)(dx) >0,
\]  and so \eqref{eq:Q-deriv} implies \eqref{eq:case1} holds
for any $t >0$ such that $\overline{Q}_j^D(t)  >0$.
\\

\noindent \textit{Case 2 (partially match, $\sum_{k \in \suppset_j} m_{jk} \in \left(0,\lambda_j^D \right)$): }
The proof of Case 2 requires the following claims, established at the end of this proof.
\begin{claim} \label{claim:strictIncreasing}
For any $t \in \Rp$, $\lran{\ind{[0,\chi)},\ovl{\eta}_j^D(t)}$ is strictly increasing in $\chi$ for $\chi \in [0,H_j^D)$.
\end{claim}
\noindent Claim~\ref{claim:strictIncreasing} ensures that we can define $\chi_j^D(t)$ as the unique solution to
\begin{equation} \label{eq:QQ}
 \lran{\ind{[0,\chi)},\ovl{\eta}_j^D(t)} = \ovl{Q}_j^D(t),
\end{equation}
which satisfies $\chi_j^D(t) \leq H_j^D$ from \eqref{in:QBs}.
\begin{claim} \label{claim:Rexpression}
Define for $t \in \Rp$,
\[
\epsilon_1(t) := \int_0^{H_j^D \wedge \left[ \chi_j^D(t)-t \right]^+} h_j^D(x+t) \frac{1-G_j^D(x+t)}{1-G_j^D(x)} \ovl{\eta}_j^D(0) (dx) \geq 0,
\]
and
\[
\epsilon_2(t) := \int_0^{H_j^D \wedge \left[ \chi_j^D(t)-t \right]^+} \frac{1-G_j^D(x+t)}{1-G_j^D(x)} \ovl{\eta}_j^D(0) (dx) \geq 0,
\]
and define for $q \in \Rp$
\[
\mathfrak{m}(q) := \lambda_j^D \left( 1-G_j^D \left( \left( G_{e,j} \right)^{-1} \left( \frac{\theta_j^D q}{\lambda_j^D} \right) \right) \right) \geq 0,
\]
where $\left( G_{e,j} \right)^{-1}$ is the inverse function of $G_{e,j}$ (that is well-defined because $G_{e,j}$ is strictly increasing).  Then $\mathfrak{m}(q)$
 is decreasing in $q$, satisfies
\[
\mathfrak{m}(q_j^{D,\star}(\bld{m})) = \sum_{k \in \suppset_j} m_{jk}, \mbox{ for } \bld{m} \in \mathbb{M},
\]
and, when $h_j^D$ is bounded,
\begin{equation} \label{eq:RR}
\ovl{R}_j^D(t) = \int_0^t \epsilon_1(s) ds + \lambda_j^D t - \int_0^t \mathfrak{m} \left( \ovl{Q}_j^D(s) - \epsilon_2(s) \right) ds.
\end{equation}
\end{claim}
\noindent  We use (\ref{eq:RR}) in Claim~\ref{claim:Rexpression} to substitute for $\ovl{R}_j^D(t)$ in \eqref{eq:QFD} to find
\[
\ovl{Q}_j^D(t) = \ovl{Q}_j^D(0) - \int_0^t \epsilon_1(s) ds + \int_0^t \mathfrak{m} \left( \ovl{Q}_j^D(s) - \epsilon_2(s) \right) ds - \sum_{k \in \suppset_j} m_{jk} t.
\]
Then, we take derivatives to find
\begin{equation} \label{eq:Qderiv2}
\frac{d}{dt} \ovl{Q}_j^D(t) = -\epsilon_1(t) + \mathfrak{m} \left( \ovl{Q}_j^D(t) - \epsilon_2(t) \right) - \sum_{k \in \suppset_j} m_{jk}.
\end{equation}
\begin{claim} \label{claim:epsilonToZero}
As $t\rightarrow \infty$, $\epsilon_1(t) \rightarrow 0$ and $\epsilon_2(t) \rightarrow 0$.
\end{claim}

Let $\gamma >0$. Set
\[
\delta_1(\gamma) := \mathfrak{m} \left( q^{D,\star}_j (\bld m) - \frac{\gamma}{2} \right) - \mathfrak{m} \left(  q^{D,\star}_j (\bld m)  \right) >0,
\]
and
\[
\delta_2(\gamma) := \mathfrak{m} \left( q^{D,\star}_j (\bld m)  \right) - \mathfrak{m} \left(  q^{D,\star}_j (\bld m) + \frac{\gamma}{2} \right) >0,
\]
which are both positive because $\mathfrak{m}(q)$ is decreasing in $q$ by Claim~\ref{claim:Rexpression} and
$q_j^{D,\star}(\bld{m}) < \lambda_j^D / \theta_j^D$ from \eqref{eq:InvDeQ}, since $\sum_{k \in \suppset_j} m_{jk} >0$.  From Claim~\ref{claim:epsilonToZero}, there exists $T_0 >0$ such that
\[
\epsilon_1(t) < \frac{1}{2} \delta_1(\gamma) \mbox{ and } \epsilon_2(t) < \frac{\gamma}{2} \mbox{ for all } t \geq T_0.
\]

Consider the Lyapunov function
\[
V(q) = \left| q -  q^{D,\star}_j (\bld m) \right|.
\]
To complete the proof, this is sufficient to show
\begin{equation} \label{eq:LyapunovCondition}
 \frac{d}{dq} V(q) \left|_{q = \ovl{Q}_j^D(t)} \right. \times \frac{d}{dt} \ovl{Q}_j^D(t) < 0
\end{equation}
for any $t \in \Rp$ such that $\left|  \ovl{Q}_j^D(t) -   q^{D,\star}_j (\bld m) \right| >0$.
We separate the argument to show \eqref{eq:LyapunovCondition} into the two cases, that $\ovl{Q}_j^D(t) - q_j^{D,\star} (\bld m) \leq -\gamma$ and that $\ovl{Q}_j^D(t) - q_j^{D,\star} (\bld m) \geq \gamma$, and then we take $\gamma \rightarrow 0$, to show (\ref{eq:LyapunovCondition}). \\

\noindent \textit{Case 2a: } $\ovl{Q}_j^D(t) - q_j^{D,\star} (\bld m) \leq -\gamma$, so that
\[
 \frac{d}{dq} V(q) \left|_{q = \ovl{Q}_j^D(t)} \right. \times \frac{d}{dt} \ovl{Q}_j^D(t) = -\frac{d}{dt} \ovl{Q}_j^D(t),
\]
and also
\begin{equation} \label{LyaCase2a}
\ovl{Q}_j^D(t) - \epsilon_2(t) \leq  q_j^{D,\star} (\bld m) - \frac{\gamma}{2}.
\end{equation}
From (\ref{eq:Qderiv2})
\[
-\frac{d}{dt} \ovl{Q}_j^D(t) = \epsilon_1(t) - \mathfrak{m}\left( \ovl{Q}_j^D(t) - \epsilon_2(t) \right) + \sum_{k \in \suppset_j} m_{jk}.
\]
From (\ref{LyaCase2a}) and the fact that $\mathfrak{m}(q)$ is decreasing in $q$,
\[
\mathfrak{m} \left( \ovl{Q}_j^D(t) - \epsilon_2(t) \right) \geq \mathfrak{m} \left( q_j^{D,\star} (\bld m) - \frac{\gamma}{2} \right),
\]
and from the definition of $\delta_1(\gamma)$,
\[
    \mathfrak{m} \left( q_j^{D,\star} (\bld m) - \frac{\gamma}{2} \right) = \mathfrak{m} \left( q_j^{D,\star} (\bld m) \right) + \delta_1(\gamma).
\]
Hence, recalling that $\mathfrak{m}\left(  q_j^{D,\star} (\bld m) \right) = \sum_{k \in \suppset_j} m_{jk}$ from Claim~\ref{claim:Rexpression},
\begin{equation} \label{eq:QderivativeBound1}
-\frac{d}{dt} \ovl{Q}_j^D(t) \leq \epsilon_1(t) - \delta_1(\gamma) < -\frac{1}{2} \delta_1(\gamma),
\end{equation}
for all $t \geq T_0$.

\noindent \textit{Case 2b: } $\ovl{Q}_j^D(t) - q_j^{D,\star} (\bld m) \geq \gamma$, so that
\[
 \frac{d}{dq} V(q) \left|_{q = \ovl{Q}_j^D(t)} \right. \times \frac{d}{dt} \ovl{Q}_j^D(t) = \frac{d}{dt} \ovl{Q}_j^D(t),
\]
and also
\begin{equation} \label{LyaCase2b}
\ovl{Q}_j^D(t) - \epsilon_2(t) \geq  q_j^{D,\star} (\bld m) + \frac{\gamma}{2}.
\end{equation}
From (\ref{eq:Qderiv2}) and the fact that $\epsilon_1(t) \geq 0$ for all $t \in \Rp$,
\[
\frac{d}{dt} \ovl{Q}_j^D(t) \leq \mathfrak{m}\left( \ovl{Q}_j^D(t) - \epsilon_2(t) \right) - \sum_{k \in \suppset_j} m_{jk}.
\]
From (\ref{LyaCase2b}) and the fact that $\mathfrak{m}(q)$ is decreasing in $q$,
\[
\mathfrak{m}\left( \ovl{Q}_j^D(t) - \epsilon_2(t) \right) \leq \mathfrak{m} \left( q_j^{D,\star} (\bld m) + \frac{\gamma}{2} \right),
\]
and from the definition of $\delta_2(\gamma)$,
\[
    \mathfrak{m} \left( q_j^{D,\star} (\bld m) + \frac{\gamma}{2} \right) = \mathfrak{m} \left( q_j^{D,\star} (\bld m) \right) - \delta_2(\gamma).
\]
Hence, recalling as above that $\mathfrak{m} \left( q_j^{D,\star} (\bld m) \right)  = \sum_{k \in \suppset_j} m_{jk}$ from Claim~\ref{claim:Rexpression},
\begin{equation}  \label{eq:QderivativeBound2}
\frac{d}{dt} \ovl{Q}_j^D(t) \leq \mathfrak{m} \left( q_j^{D,\star} (\bld m) \right) - \delta_2(\gamma) - \sum_{k \in \suppset_j} m_{jk} = -\delta_2(\gamma).
\end{equation}

Finally, to see \eqref{eq:LyapunovCondition} holds and complete the proof of Case 2, let $\gamma \rightarrow 0$ in (\ref{eq:QderivativeBound1}) and (\ref{eq:QderivativeBound2}), and observe that $\delta_1(\gamma) \rightarrow 0$ and $\delta_2(\gamma) \rightarrow 0$ as $\gamma \rightarrow 0$ from the definitions of $\delta_1(\gamma)$ and $\delta_2(\gamma)$. \\

\noindent \textit{Case 3 (match nothing, $\sum_{k \in \suppset_j} m_{jk} =0$): }   Since $\ovl{Q}_j^D(t) \in \Rp$ for all $t \geq 0$ and \eqref{in:QBs} holds
\[
\ovl{Q}_j^D(t) \in \left[ 0,  \lran{1,\ovl{\eta}_j^D(t)} \right].
\]
From \eqref{eq:eta-long-term-behavior} and the definition of $\eta_j^{D,\star}$ in \eqref{eq:InvDe},
\[
 \lran{1,\ovl{\eta}_j^D(t)} \rightarrow \lran{1,\eta_j^{D,\star}} = \frac{\lambda_j^D}{\theta_j^D}, \mbox{ as } t \rightarrow \infty.
\]
Then,
\[
\lim_{t \rightarrow \infty} \frac{\ovl{Q}_j^D(t)}{t} = 0,
\]
and so, from \eqref{eq:QFD}
\begin{equation} \label{eq:averageRlimit}
\lim_{t \rightarrow \infty} \frac{\ovl{R}_j^D(t)}{t} = \lambda_j^D.
\end{equation}
Also, from Claim~\ref{claim:Rexpression},
\begin{equation} \label{eq:averageR}
\frac{1}{t}\ovl{R}_j^D(t) = \frac{1}{t}\int_0^t \epsilon_1(s) ds + \lambda_j^D  - \frac{1}{t} \int_0^t \mathfrak{m} \left( \ovl{Q}_j^D(s) - \epsilon_2(s) \right) ds.
\end{equation}
Since $\epsilon_1(s) \rightarrow 0$ as $s \rightarrow \infty$ by Claim~\ref{claim:epsilonToZero}, together (\ref{eq:averageRlimit}) and (\ref{eq:averageR}) imply
\[
\frac{1}{t} \int_0^t \mathfrak{m} \left( \ovl{Q}_j^D(s) - \epsilon_2(s) \right) ds \rightarrow 0, \mbox{ as } t \rightarrow \infty.
\]
Since $\mathfrak{m}$ is decreasing by Claim~\ref{claim:Rexpression},
\[
\frac{1}{t} \int_0^t \mathfrak{m} \left( \ovl{Q}_j^D(s) \right) ds \rightarrow 0, \mbox{ as } t \rightarrow \infty.
\]
Since $\mathfrak{m}$ is also non-negative and continuous from its definition,
\[
\mathfrak{m} \left( \ovl{Q}_j^D(t) \right) \rightarrow 0, \mbox{ as } t \rightarrow \infty,
\]
which implies
\[
 \ovl{Q}_j^D(t) \rightarrow \frac{\lambda_j^D}{\theta_j^D} = q^{D,\star}_j (\bld m), \mbox{ as } t \rightarrow \infty.
\]

\pfa{Proof of Claim~\ref{claim:strictIncreasing}}
First note that \eqref{eq:FetaD} holds for every bounded, Borel measurable function $f$ and $t\geq 0$ by Corollary 4.2 in \cite{kang2010fluid}. Next, apply \eqref{eq:FetaD} to $f(x) = \ind{[0,\chi)}(x)$ to find
\[
\lran{\ind{[0,\chi)},\ovl{\eta}_j^D(t)} =  \int_0^{H_j^D \wedge \left[ \chi-t \right]^+} \frac{1-G_j^D(x+t)}{1-G_j^D(x)} \ovl{\eta}_j^D(0) (dx) + \lambda_j^D \int_0^{t\wedge \chi} \left( 1-G_j^D(w) \right) dw.
\]
The first term on the right-hand side of the above equation is strictly increasing in $\chi$ for $\chi \in [t,H_j^D)$, and the second term is strictly increasing in $\chi$ for $\chi \in [0,t]$. \\

\pfa{Proof of Claim~\ref{claim:Rexpression}}
First note that $\mathfrak{m}(q)$ is decreasing in $q$ because both $G_j^D$ and $G_{e,j}^D$ are increasing functions, and that $\mathfrak{m} \left( q_j^{D,\star} (\bld m) \right) = \sum_{k \in \suppset_j} m_{jk}$ from \eqref{eq:InvDeQ} and the definition of $\mathfrak{m}$.  Hence showing (\ref{eq:RR}) is sufficient to prove the claim.

For fixed $t \in \Rp$, let
\[
f(x) = h_j^D(x) \ind{[0,\chi_j^D(t))}(x),
\]
which is a bounded function because $h_j^D$ is bounded, and is also Borel measurable.  By Corollary 4.2 in~\cite{kang2010fluid} (with their $\overline{E}(t) = \lambda_j^D t, t \in \Rp$), \eqref{eq:FetaD} holds for every bounded Borel measurable function $f$ and $t \in \Rp$.  Substituting for $f$ in \eqref{eq:FetaD} shows
\[
\lran{f, \ovl{\eta}_j^D(t)} = \epsilon_1(t) + \lambda_j^D \int_0^{\chi_j^D(t) \wedge t} h_j^D(\zeta) \left( 1-G_j^D(\zeta) \right) d\zeta.
\]
Next, from \eqref{eq:AlRD}, the definition of $\chi_j^D(t)$, the above equality, and the relationship between $G_j^D$ and $h_j^D$,
\begin{eqnarray} \label{eq:RexpressOther}
\ovl{R}_j^D(t) & = & \int_0^t \int_0^{\chi_j^D(t)} h_j^D(x) \ovl{\eta}_j^D(u)(dx) du \\
& = & \int_0^t \lran{f, \ovl{\eta}_j^D(u)} du  \nonumber \\
& = & \int_0^t \epsilon_1(s) ds + \lambda_j^D \int_0^t \int_0^{\chi_j^D(s) \wedge s} h_j^D(\zeta) \left( 1-G_j^D(\zeta) \right) d\zeta ds \nonumber \\
& = & \int_0^t \epsilon_1(s) ds + \lambda_j^D \int_0^t G_j^D \left( \chi_j^D(s) \wedge s \right) ds. \nonumber
\end{eqnarray}

Finally, we show the last term in (\ref{eq:RexpressOther}) equals the last two terms in (\ref{eq:RR}), which is enough to complete the proof.  From \eqref{eq:QQ}, for any $t \in \Rp$,
\[
\ovl{Q}_j^D(t) = \lran{\ind{[0,\chi_j^D(t))},\ovl{\eta}_j^D(t)},
\]
and so, since $\ind{[0,\chi_j^D(t))}$ is bounded and Borel measurable, substituting for $f$ in (\ref{eq:FetaD}) shows
\begin{equation} \label{eq:queueXX}
\ovl{Q}_j^D(t) = \epsilon_2(t) + \lambda_j^D \int_0^{\chi_j^D(t) \wedge t} \left( 1-G_j^D(\zeta) \right) d\zeta.
\end{equation}
Recalling the definition of the excess life distribution $G_{e,j}^D$ near the beginning of Section~\ref{sec:Model description}, we have that
\[
\lambda_j^D \int_0^{\chi_j^D(t) \wedge t} \left( 1-G_j^D(\zeta) \right) d\zeta = \frac{\lambda_j^D}{\theta_j^D} G_{e,j}^D \left( \chi_j^D(t) \wedge t \right),
\]
and so, solving for $\chi_j^D \wedge t$ in (\ref{eq:queueXX}) show that
\[
\chi_j^D(t) \wedge t = \left( G_{e,j}^D \right)^{-1} \left( \frac{\theta_j^D}{\lambda_j^D} \left( \ovl{Q}_j^D(t) - \epsilon_2(t) \right) \right).
\]
Then, from algebra, substitution, and the definition of $\mathfrak{m}$,
\begin{eqnarray*}
\lambda_j^D \int_0^t G_j^D \left( \chi_j^D(s) \wedge s \right) ds & = &
    \lambda_j^D t - \lambda_j^D \int_0^t \left( 1-G_j^D \left( \chi_j^D(s) \wedge s \right) \right) ds \\
    & = & \lambda_j^D t - \int_0^t \mathfrak{m} \left( \ovl{Q}_j^D(s) - \epsilon_2(s) \right) ds.
\end{eqnarray*}
Substituting for the above equality in (\ref{eq:RexpressOther}) shows that
\[
\ovl{R}_j^D(t) = \int_0^t \epsilon_1(s) ds + \lambda_j^D t - \int_0^t \mathfrak{m}\left( \ovl{Q}_j^D(s) - \epsilon_2(s) \right) ds
\]

\pfa{Proof of Claim~\ref{claim:epsilonToZero}}
Since
\[
0 \leq \epsilon_2(t) \leq \int_0^{H_j^D} \frac{1-G_j^D(x+t)}{1-G_j^D(x)} \ovl{\eta}_j^D(0) (dx),
\]
and, by dominated convergence,
\[
\int_0^{H_j^D} \frac{1-G_j^D(x+t)}{1-G_j^D(x)} \ovl{\eta}_j^D(0) (dx) \rightarrow 0,
\]
as $t \rightarrow \infty$, we conclude that $\epsilon_2(t) \rightarrow 0$ as $t \rightarrow \infty$.  Next, since $h_j^D(x) \left( 1-G_j^D(x) \right) = g_j^D(x)$ for any $x \in [0,H_j^D)$, we have that
\[
\epsilon_1(t) = \int_0^{H_j^D \wedge \left[ \chi_j^D(t) -t \right]^+} \frac{g_j^D(x+t)}{1-G_j^D(x)} \ovl{\eta}_j^D(0) (dx).
\]
Since
\[
0 \leq \epsilon_1(t) \leq \int_0^{H_j^D} \frac{g_j^D(x+t)}{1-G_j^D(x)} \ovl{\eta}_j^D(0) (dx),
\]
and, by dominated convergence,
\[
\int_0^{H_j^D} \frac{g_j^D(x+t)}{1-G_j^D(x)} \ovl{\eta}_j^D(0) (dx) \rightarrow 0,
\]
as $t \rightarrow \infty$, we conclude that $\epsilon_1(t) \rightarrow 0$ as $t \rightarrow \infty$.
}
\eof

\subsection{Proofs for Section~\ref{sec:Performance Analysis}}
\label{Proofs for Section sec:High volume}
\pfa{Proof of Theorem~\ref{prop:upper}}
\amy{
The sequence $\left\{ \left( \frac{\bld{Q}^{D,n}}{n}, \frac{\bld{Q}^{S,n}}{n}, \frac{\bld{\eta}^{D,n}}{n}, \frac{\bld{\eta}^{S,n}}{n} \right) \right\}_{n \in \mathbb{N}}$ is tight by Theorem 2 in~\cite{aveklouris2021MM}.  Furthermore, any distributional limit point $(\overline{\bld{Q}}^D, \overline{\bld{Q}}^S, \overline{\bld{\eta}}^D, \overline{\bld{\eta}}^S)$ of $\left\{ \left( \frac{\bld{Q}^{D,n}}{n}, \frac{\bld{Q}^{S,n}}{n}, \frac{\bld{\eta}^{D,n}}{n}, \frac{\bld{\eta}^{S,n}}{n} \right) \right\}_{n \in \mathbb{N}}$ is almost surely a fluid model solution for $\left( \bld{\lambda}^D, \bld{\lambda}^S \right)$ with initial state zero.  Then, on any convergent subsequence $\{ n_l\}$, by the continuous mapping theorem, for any $t \in \Rp$,
\begin{equation}\label{in:Up1}
\frac{V_{\bld{M}^{n_l}}(t)}{n_l} \Rightarrow
 \sum_{(j,k) \in \arcset}
 v_{jk} m_{jk} t -
 \sum_{j\in \bbJ} \int_{0}^{t}c_j^D \ovl Q_j^{D}(s)ds -
 \sum_{k\in \bbK} \int_{0}^{t} c_k^S \ovl Q_k^{S}(s)ds,
\end{equation}
as $n_l \rightarrow \infty$.
Next, since by Theorem~\ref{Th:steadySMB},
\[
\lim_{t \rightarrow \infty} \left( \ovl{\bld{Q}}^D(t), \ovl{\bld{Q}}^S(t) \right) = \left(\bld{q}^{D,\star}(\bld m), \bld{q}^{S,\star}(\bld m) \right),
\]
we have that
\begin{eqnarray} \label{eq:limit-t-average}
\lefteqn{lim_{t \rightarrow \infty} \frac{1}{t} \left( \sum_{(j,k) \in \arcset}
 v_{jk} m_{jk} t -
 \sum_{j\in \bbJ} \int_{0}^{t}c_j^D \ovl Q_j^{D}(s)ds -
 \sum_{k\in \bbK} \int_{0}^{t} c_k^S \ovl Q_k^{S}(s)ds \right)} \\
 & & =  \sum_{(j,k) \in \arcset}
 v_{jk} m_{jk} - \sum_{j\in \bbJ} c_j^D q_j^{D,\star}(\bld m) - \sum_{k\in \bbK}  c_k^S q_k^{S,\star}(\bld m). \nonumber
\end{eqnarray}
From the definition of $\bld{m}^\star$,
\begin{equation} \label{eq:m-star-inequality}
\sum_{(j,k) \in \arcset}
 v_{jk} m_{jk} - \sum_{j\in \bbJ} c_j^D q_j^{D,\star}(\bld m) - \sum_{k\in \bbK}  c_k^S q_k^{S,\star}(\bld m)
\leq
\sum_{(j,k) \in \arcset}
 v_{jk} m^\star_{jk} - \sum_{j\in \bbJ} c_j^D q_j^{D,\star}(\bld m^\star) - \sum_{k\in \bbK}  c_k^S q_k^{S,\star}(\bld m^\star).
\end{equation}
We conclude from (\ref{in:Up1}), (\ref{eq:limit-t-average}), and (\ref{eq:m-star-inequality}) that
\[
\lim_{t \rightarrow \infty} \lim_{n\to \infty}
\frac{1}{t} \frac{V_{\bld{M}^{n_l}}(t)}{n_l} \leq \sum_{(j,k) \in \arcset}
 v_{jk} m^\star_{jk} - \sum_{j\in \bbJ} c_j^D q_j^{D,\star}(\bld m^\star) - \sum_{k\in \bbK}  c_k^S q_k^{S,\star}(\bld m^\star).
\]
Since the above arguments hold on any convergent subsequence $\{n_l\}$, we have that \eqref{eq:max-achievable-profit-fluid} holds, which completes the proof.}
\eof

\subsubsection{Proof of Theorems~\ref{proposition:DRconvergence2} and \ref{Thm:aopriority}}

Before we present the proof of Theorems~\ref{proposition:DRconvergence2} and \ref{Thm:aopriority}, we show some preliminary results
that are also used in the proofs in Section~\ref{sec:PSB}.  In the following, $r_{jm}^D$ denotes the patience time of the $m$th customer arriving at note $j \in \bbJ$ and $r_{km}^S$ denotes the patience time of the $m$th worker arriving at node $k \in \bbK$, for $m \in \mathbb{N}$.

\begin{proposition}\label{prop:LBR}
For any $j\in \bbJ$, $k\in \bbK$, and $i \ge 1$, the following inequalities hold almost surely
\begin{equation*}
R_{j}^{D}(il)-R_{j}^{D}((i-1)l)
\le Q_{j}^D((i-1)l)+
 \sum_{m=\ad_j((i-1)l)+1}^{\ad_j(il)} \ind{r_{jm}^D\le l},
\end{equation*}
\begin{equation*}
R_{k}^{S}(il)-R_{k}^{S}((i-1)l)
\le Q_{k}^S((i-1)l)+ \sum_{m=\as_k((i-1)l)+1}^{\as_k(il)}
\ind{r_{km}^S\le l}.
\end{equation*}
\end{proposition}
\pfa{Proof}
\amy{The proof of Proposition~\ref{prop:LBR} relies on the detailed evolution equations for the discrete-event stochastic system in Section 2 in \cite{aveklouris2021MM}, and the reader may want to refer to that paper when reading this proof. The following notation is used from that paper.  First, $m_{jh}^D$ denotes the matching time of the $h$th type $j$ customer if that customer is matched, and is infinity otherwise.  Also,  the potential waiting time of the $h$th type $j$ customer at time $t \geq 0$ (potential, because matching is ignored) is
\begin{equation*}
w^D_{jh}(t):=\min \left\{ [t-e_{jh}^D]^+, r_{jh}^D\right\}.
\end{equation*}
Then, the cumulative number of type $j$ customers that renege by time $t\geq 0$ when the system is initially empty is given by
\begin{equation} \label{eq:cumulative-reneging-D}
R^D_{j}(t):=\sum_{l=1}^{\ad_j(t)}\sum_{s\in [0, t]}
 \ind{s\leq m_{jl}^D,\frac{dw^D_{jl}}{dt}(s-)>0, \frac{dw^D_{jl}}{dt}(s+)=0}
\end{equation}
Analogous expressions related to the workers can be found in \cite{aveklouris2021MM}.
 }

We show only the first inequality, the second one follows in the same way.
\amy{From the expression for the cumulative number of reneging customers in (\ref{eq:cumulative-reneging-D}) above}, we have that
\begin{equation*}
\begin{split}
R_{j}^{D}(il)-R_{j}^{D}((i-1)l)
=
\sum_{h=A^D_{j}((i-1)l)+1}^{A^D_{j}(il)}
\sum_{s\in [(i-1)l, il]}
   \ind{s\leq m_{jh}^D,\frac{dw^D_{jh}}{dt}(s-)>0,
       \frac{dw^D_{jh}}{dt}(s+)=0}\\
 +
\sum_{h=1}^{A^D_{j}((i-1)l)}
\sum_{s\in [(i-1)l, il]}
   \ind{s\leq m_{jh}^D,\frac{dw^D_{jh}}{dt}(s-)>0,
       \frac{dw^D_{jh}}{dt}(s+)=0}.
 \end{split}
\end{equation*}

Now, we have that
\begin{equation*}
\begin{split}
\sum_{h=A^D_{j}((i-1)l)+1}^{A^D_{j}(il)}
\sum_{s\in [(i-1)l, il]}
   \ind{s\leq m_{jh}^D,\frac{dw^D_{jh}}{dt}(s-)>0,
       \frac{dw^D_{jh}}{dt}(s+)=0}
\le
\sum_{h=\ad_j((i-1)l)+1}^{\ad_j(il)} \ind{r_{jh}^D\le l}.
\end{split}
\end{equation*}
To see this, observe that if
$\ind{s\leq m_{jh}^D,\frac{dw^D_{jh}}{dt}(s-)>0,
       \frac{dw^D_{jh}}{dt}(s+)=0}=1$, then
$\frac{dw^D_{jh}}{dt}(s+)=0$. That is, $s-e_{jh}^D>r_{jh}^D$. This yields
$l\ge r_{jh}^D$ by the fact that $s\le il $ and $e_{jh}^D\ge  (i-1)l $. In other words, $\ind{r_{jh}^D\le l}=1$.

To finish the proof it remains to show that
\begin{equation*}
\begin{split}
\sum_{h=1}^{A^D_{j}((i-1)l)}
\sum_{s\in [(i-1)l, il]}
   \ind{s\leq m_{jh}^D,\frac{dw^D_{jh}}{dt}(s-)>0,
       \frac{dw^D_{jh}}{dt}(s+)=0}
\le
Q^D_j((i-1)l).
\end{split}
\end{equation*}
To this end, the left hand side of the last inequality equals to
\begin{equation*}
\begin{split}
&\sum_{h=1}^{A^D_{j}((i-1)l)}
\sum_{s\in [0, il]}
   \ind{s\leq m_{jh}^D,\frac{dw^D_{jh}}{dt}(s-)>0,
       \frac{dw^D_{jh}}{dt}(s+)=0}
-
R_j^D((i-1)l)\\
& =
\sum_{h=1}^{A^D_{j}((i-1)l)}
\sum_{s\in [0, il]}
   \ind{s\leq m_{jh}^D,\frac{dw^D_{jh}}{dt}(s-)>0,
       \frac{dw^D_{jh}}{dt}(s+)=0}
-\ad_j((i-1)l) \\
& \qquad +\sum_{k \in \suppset_j}M_{jk}((i-1)l)+Q_{j}^D((i-1)l).
\end{split}
\end{equation*}
Now observe that
\begin{equation*}
\begin{split}
&\sum_{h=1}^{A^D_{j}((i-1)l)}
\sum_{s\in [0, il]}
   \ind{s\leq m_{jh}^D,\frac{dw^D_{jh}}{dt}(s-)>0,
       \frac{dw^D_{jh}}{dt}(s+)=0}
+\sum_{k \in \suppset_j}M_{jk}((i-1)l)\\
&=\sum_{h=1}^{A^D_{j}((i-1)l)}
\sum_{s\in [0, (i-1)l]}
   \ind{s\leq m_{jh}^D,\frac{dw^D_{jh}}{dt}(s-)>0,
       \frac{dw^D_{jh}}{dt}(s+)=0} \\
&\qquad+\sum_{h=1}^{A^D_{j}((i-1)l)}\sum_{s\in ((i-1)l, il]}
   \ind{s\leq m_{jh}^D,\frac{dw^D_{jh}}{dt}(s-)>0,
       \frac{dw^D_{jh}}{dt}(s+)=0}  \\
&\qquad+\sum_{k \in \suppset_j}M_{jk}((i-1)l)\\
&\le \ad_j((i-1)l).
\end{split}
\end{equation*}
To see the last inequality, note that the left-hand side counts the number of customers than arrive in time interval $[0,(i-1)l]$ that depart (due to matching or reneging) until time $(i-1)l$ and  renege in time interval $((i-1)l,il]$. Hence it cannot exceed the number of customers than arrive until time $(i-1)l$ which is $\ad_j((i-1)l)$.
\eof \\

Moreover, we need the following lemma, which is taken from \cite{ata2005}, and adapted to the below form found in Lemma 4.1 in \cite{PlambeckWard2006}.
\begin{lemma}(Ata and Kumar) \label{lemma:AtaKumar}
\amy{When the inter-arrival distributions have finite 5th moment,} for any finite constant $\alpha >0$, there exists $\beta >0$ such that
\[
\Prob{
 \max_{i \in \{1,\ldots, \lfloor T/l^n \rfloor \}}
 \max_{j \in \bbJ}
  \left| A_j^{D,n}( i l^n) - A_j^{D,n}((i-1)l^n) - n \lambda_j^D l^n \right| < \alpha n^{1/3}} \geq 1- \beta n^{-1/6}
\]
and
\[
\Prob { \max_{i \in \{1,\ldots, \lfloor T/l^n \rfloor \}} \max_{k \in \bbK} \left| \asn_k( i l^n) - \asn_k((i-1)l^n) - n \lambda^S_{k} l^n \right| < \alpha n^{1/3} } \geq 1- \beta n^{-1/6}.
\]
\end{lemma}
A consequence of Lemma~\ref{lemma:AtaKumar} is that with high probability the following hold for each $j \in \bbJ$, $k \in \bbK$, and
$i \in \{1,\ldots, \lfloor T/l^n \rfloor\}$,
\begin{equation}\label{eq:boundD1}
n \lambda_j^D l^n - \alpha n^{1/3}<A_j^{D,n}( i l^n) - A_j^{D,n}((i-1)l^n) <
 n \lambda_j^D l^n + \alpha n^{1/3}
\end{equation}
and
\begin{equation}\label{eq:boundS1}
n \lambda^S_{k} l^n - \alpha n^{1/3}<A_k^{S,n}( i l^n) - A_k^{S,n}((i-1)l^n) < n
\lambda^S_{k} l^n + \alpha n^{1/3}.
\end{equation}

\pfa{Proof of Theorem~\ref{proposition:DRconvergence2}}
Let $\epsilon>0$. Fix $n$ large enough so that
$n^{2/3} > 2 \max_{j,k}(m_{jk}) l  / \epsilon$. Denote by $\Omega_1^n$ the event such that \eqref{eq:boundD1} and \eqref{eq:boundS1} hold and by $\Omega_2^n$ the event such that both inequalities in Proposition~\ref{prop:LBR} hold. Let
$\omega \in \Omega_1^n \cap \Omega_2^n$ and
$\alpha :=\frac{ \min_{(j,k) \in \arcset}\left(\lambda_j^D, \lambda^S_{k}\right) l}{6T\min_{(j,k) \in \arcset}\left(m_{jk}\right)}\epsilon$.
We first derive a lower bound for the number of reneging \levi{agents} at any discrete review period.
By Proposition~\ref{prop:LBR},
the following upper bound of the reneging customers in the discrete review period $((i-1)l^n,il^n)$ holds almost surely: for any $j\in \bbJ$,
\begin{equation}\label{eq:boundRD}
R_{j}^{D,n}(il^n)-R_{j}^{D,n}((i-1)l^n)
\le Q_{j}^{D,n}((i-1)l^n)+
 \sum_{m=A^{D,n}_j((i-1)l^n)+1}^{A^{D,n}_j(il^n)} \ind{r_{jm}^D\le l^n},
\end{equation}
recalling that $r_{jm}^D$ denotes the patience time of the $m$th customer arriving at node $j\in \bbJ$.
Further, for the reneging workers in the discrete review period $((i-1)l^n,il^n)$, we have that for any $k\in \bbK$,
\begin{equation}\label{eq:boundRS}
R_{k}^{S,n}(il^n)-R_{k}^{S,n}((i-1)l^n)
\le Q_{k}^{S,n}((i-1)l^n)+ \sum_{m=A^{S,n}_k((i-1)l^n)+1}^{A^{S,n}_k(il^n)}
\ind{r_{km}^S\le l^n},
\end{equation}
recalling that $r_{km}^S$ denotes the patience time of the $m$th worker arriving at node $k\in \bbK$. In the sequel, we show that the second term of the right-hand side of \eqref{eq:boundRD} and \eqref{eq:boundRS} converge to zero under the high-volume setting.
Denote by $\Omega_0^n$ the event such that \eqref{eq:boundD1} holds and by $(\Omega_0^n)^c$ its complement. For simplicity, let
$\Gamma^n=\frac{1}{nl^n}\sum_{m=A^{D,n}_j((i-1)l^n)+1}^{A^{D,n}_j(il^n)}
 \ind{r_{jm}^D\le l^n}$.
For any $\delta>0$, we have that
\begin{equation}\label{eq:convRene}
\begin{split}
\Prob{
\Gamma^n\ge \delta}
&\le
\Prob{ \left\{
\Gamma^n\ge \delta\right\} \cap \Omega_0^n}+
\Prob{\left\{
\Gamma^n\ge \delta \right\}\cap (\Omega_0^n)^c}\\
&\le
\Prob{
\Gamma^n\ge \delta |\Omega_0^n}
+
\Prob{(\Omega_0^n)^c}\\
&\le
\frac{\E{
\Gamma^n |\Omega_0^n}}{\delta}
+
\beta n^{-1/6}\\
&\le
\frac{\E{
\Gamma^n 1_{\Omega_0^n}}}{\delta(1-\beta n^{-1/6})}
+
\beta n^{-1/6},
\end{split}
\end{equation}
where the third inequality follows by the conditional Markov's inequality and Lemma~\ref{lemma:AtaKumar}. To complete the proof that $\Gamma^n$ converges to zero in probability, we show that
$\E{\Gamma^n 1_{\Omega_0^n}}$ goes to zero. To this end,
\begin{equation*}
\begin{split}
0 \le \E{\frac{1}{nl^n}\sum_{m=A^{D,n}_j((i-1)l^n)+1}^{A^{D,n}_j(il^n)}
 \ind{r_{jm}^D\le l^n}1_{\Omega_0^n}}
&\le \frac{1}{nl^n} \E{\left(n \lambda_j^D l^n+\alpha n^{1/3}  \right) \ind{r_{j1}^D\le l^n}}\\
&= \frac{\left(n \lambda_j^D l^n+\alpha n^{1/3}  \right)}{nl^n} \Prob{r_{j1}^D\le l^n}\\
&= \left( \lambda_j^D + \frac{\alpha}{l}\right)
 \Prob{r_{j1}^D\le l^n}.
\end{split}
\end{equation*}
It follows that
\begin{equation*}
\begin{split}
0 \le \E{ \frac{1}{nl^n }\sum_{m=A^{D,n}_j((i-1)l^n)+1}^{A^{D,n}_j(il^n)}
 \ind{r_{jm}^D\le l^n}1_{\Omega_0^n}} \rightarrow 0,
\end{split}
\end{equation*}
as $n \rightarrow \infty$ because $l^n \to 0$. Now, by \eqref{eq:convRene}, we have that
$\lim_{n\rightarrow \infty}\Prob{\Gamma^n\ge \delta}=0$,
for any $\delta>0$.
By the last convergence and for large enough $n$, we have that
 \begin{equation}\label{eq:BoundRD2}
\begin{split}
 \frac{1}{nl^n }\sum_{m=A^{D,n}_j((i-1)l^n)+1}^{A^{D,n}_j(il^n)} \ind{r_{jm}^D\le l^n}
  <
 \frac{\lambda_j^D}
 {6T\max_{j,k}\left(m_{jk}\right)} \epsilon.
\end{split}
\end{equation}
In a similar way using \eqref{eq:boundS1}, we obtain
for large enough $n$,
 \begin{equation}\label{eq:BoundRS2}
\begin{split}
 \frac{1}{nl^n }\sum_{m=A^{S,n}_k((i-1)l^n)+1}^{A^{S,n}_k(il^n)} \ind{r_{km}^S\le l^n}
  <
  \frac{ \lambda^S_{k}}
  {6T\max_{j,k}\left(m_{jk}\right)} \epsilon.
\end{split}
\end{equation}
Define $\Omega_3^n$ the events such that \eqref{eq:BoundRD2} and
\eqref{eq:BoundRS2} hold and note that
$\lim_{n\rightarrow \infty}
\Prob{\Omega_1^n \cap \Omega_2^n  \cap \Omega_3^n }=1$. In the sequel, we take $\omega \in \Omega_1^n \cap \Omega_2^n  \cap \Omega_3^n$.

Now, we move to the second step of proof of Theorem~\ref{proposition:DRconvergence2}. We  derive a desirable lower bound on the number of matches at any discrete review period. To this end, by \eqref{eq:boundD1}, \eqref{eq:boundRD}, and \eqref{eq:BoundRD2}, we have that for $i \in \{1,\ldots, \lfloor T/l^n \rfloor \}$,
\begin{equation*}
\begin{split}
Q^{D,n}_j(il^n-)&=Q^{D,n}_j((i-1)l^n)+ A_j^{D,n}(il^n) - A_j^{D,n}((i-1)l)- R^{D,n}_j(il^n) + R^{D,n}_j((i-1)l^n)\\
&\ge
A_j^{D,n}(il^n) - A_j^{D,n}((i-1)l)-
\sum_{m=A^{D,n}_j((i-1)l^n)+1}^{A^{D,n}_j(il^n)}
\ind{r_{jm}^D\le l^n}\\
&\ge
A_j^{D,n}(il^n) - A_j^{D,n}((i-1)l)
-\frac{\lambda_j^D}
{6T \max_{j,k}\left(m_{jk}\right)}  nl^n \epsilon\\
&\ge
n \lambda_j^D l^n-\alpha n^{1/3}
-\frac{\lambda_j^D}
{6T \max_{j,k}\left(m_{jk}\right)}nl^n \epsilon\\
&=
n \lambda_j^D l^n\left(1 -\frac{\alpha}{l\lambda_j^D} -
\frac{\epsilon}
{6T\max_{j,k}\left(m_{jk}\right)} \right)\\
&\ge
n \lambda_j^D l^n\left(1 -\frac{\alpha}{l\min_{j,k}
\left(\lambda_j^D, \lambda^S_{k}\right)} -
\frac{\epsilon}
{6T\max_{j,k}\left(m_{jk}\right)} \right),
\end{split}
\end{equation*}
and by \eqref{eq:boundS1}, \eqref{eq:boundRS}, and \eqref{eq:BoundRS2}, we have that
\begin{equation*}
\begin{split}
Q_k^{S,n}(il^n-)&=Q_k^{S,n}((i-1)l^n) + A_k^{S,n}( i l^n )- A_k^{S,n}((i-1)l^n)
 - R^{S,n}_k(il^n) + R^{S,n}_k((i-1)l^n)\\
 &\ge
n \lambda^S_{k} l^n\left(1 -\frac{\alpha}{l\min_{j,k}
\left(\lambda_j^D, \lambda^S_{k}\right)} -
\frac{\epsilon}
{6T\max_{j,k}\left(m_{jk}\right)} \right).
\end{split}
\end{equation*}
Now, the inequality $\lfloor x \rfloor \ge x -1$ and the bounds for the quantities $Q^{D,n}_j(il^n-)$, $Q_k^{S,n}(il^n-)$ yield
\begin{equation*}
\begin{split}
\frac{\mathcal{M}_{ijk}^{r,n}}{nl^n} &= \frac{1}{n l^n}
\left \lfloor
m_{jk}
\min\left(\amy{n^{1/3}l}, \frac{Q^{D,n}_j(il^n-)}{\lambda_j^D},
\frac{Q_k^{S,n}(il^n-)}{\lambda^S_{k}} \right)
\right\rfloor
\\
&\ge
m_{jk} \left(1 -\frac{\alpha}{l\min_{j,k}
\left(\lambda_j^D, \lambda^S_{k}\right)} -
\frac{\epsilon}
{6T\max_{j,k}\left(m_{jk}\right)} \right)
- \frac{1}{n l^n}\\
& \ge
m_{jk} -\frac{\epsilon}{6T} -
\frac{\epsilon}
{6T}
-\frac{\epsilon}
{6T}
=
m_{jk}-\frac{\epsilon}
{2T},
\end{split}
\end{equation*}
for large enough $n$.
The second inequality follows from the definition of $\alpha$ and the fact that
$\frac{1}{n l^n}=\frac{1}{n^{1/3} l}\rightarrow 0$,
$n \rightarrow \infty$.
That is, for any $t\ge 0$ and any feasible point $\bld{m}$,
\begin{equation*}\label{in:CLM}
\begin{split}
\frac{M_{jk}^{r,n}(t)}{n}=
\frac{1}{n}
 \sum_{i=1}
 ^{ \lfloor t/l^n \rfloor }
\mathcal{M}_{ijk}^{r,n}
\ge
\lfloor t/l^n \rfloor
l^n(m_{jk}  - \frac{\epsilon}{2T})
\ge
\lfloor t/l^n \rfloor l^n
m_{jk}  - \frac{\epsilon}{2}
\ge
m_{jk} t
-m_{jk} l^n  - \frac{\epsilon}{2}
\ge
m_{jk} t- \epsilon,
\end{split}
\end{equation*}
where the last inequality follows by $n^{2/3} > 2 \max_{j,k}(m_{jk}) l  / \epsilon$.

On the other hand, we have $\frac{M_{jk}^{r,n}(t)}{n}\le m_{jk} t- \epsilon$ by definition of the matching-rate-based policy for any $0 \le t\le T$ and $\epsilon$ is arbitrary small independent of time $t$. This concludes the proof.
\eof\\

\pfa{Proof of Theorem~\ref{Thm:aopriority}}
Let $\epsilon>0$. Denote by $\Omega_1^n$ the events such that \eqref{eq:boundD1} and \eqref{eq:boundS1} hold, by $\Omega_2^n$ the events such both inequalities in Proposition~\ref{prop:LBR} hold, and by $\Omega_3^n$ the events such that \eqref{eq:BoundRD2} and \eqref{eq:BoundRS2} hold. Note that
$\lim_{n\rightarrow \infty}
\Prob{\Omega_1^n \cap \Omega_2^n  \cap \Omega_3^n }=1$
and consider $\omega \in \Omega_1^n \cap \Omega_2^n  \cap \Omega_3^n$. To show the result is enough to prove that for each $(j,k) \in \arcset$, $i\ge 1$, and large enough
$n \in \mathbb{N}$,
\begin{equation}\label{in:bounds}
 n l^n (m_{jk}^\star-\frac{\epsilon}{T}) \le \calM_{ijk}^{p,n} \le
 n l^n (m_{jk}^\star+\frac{\epsilon}{T}).
\end{equation}
Then, we have that for $n \in \mathbb{N}$,
\begin{equation*}
\begin{split}
\frac{M_{jk}^{p,n}(t)}{n}=
\frac{1}{n}
 \sum_{i=1}
 ^{ \lfloor t/l^n \rfloor }
\mathcal{M}_{ijk}^{p,n}
\ge
\lfloor t/l^n \rfloor
l^n(m_{jk}^\star  - \frac{\epsilon}{T})
\ge
 (t/l^n-1) l^n
(m_{jk}^\star  - \frac{\epsilon}{T})
\ge
m_{jk}^\star t- \epsilon,
\end{split}
\end{equation*}
because $l^n=1/n^{1/3}\to 0$ as $n\to \infty$. On the other hand, we have that
\begin{equation*}
\begin{split}
\frac{M_{jk}^{p,n}(t)}{n}=
\frac{1}{n}
 \sum_{i=1}
 ^{ \lfloor t/l^n \rfloor }
\mathcal{M}_{ijk}^{p,n}
\le
\lfloor t/l^n \rfloor
l^n(m_{jk}^\star  + \frac{\epsilon}{T})
\le
 m_{jk}^\star t  +\epsilon.
\end{split}
\end{equation*}
We next focus on proving \eqref{in:bounds}. Define for any
$(j,k)\in \mathcal{P}_{h}(\bld{m}^\star)$ with $h = 0, \dots, H$,
$$\delta_{jk}^h:=\left|\lambda_j^D- \sum_{k':(j,k')\in \mathcal{Q}_{h-1}(\bld{m}^\star)} m_{jk'}^\star
-\lambda_k^S+\sum_{j':(j',k)\in \mathcal{Q}_{h-1}(\bld{m}^\star)} m_{j'k}^\star\right|$$
and
$\delta^+:= \begin{cases}
              \min\limits_{j,k,h}(\delta_{jk}^h: \delta_{jk}^h>0) , & \mbox{if } \exists\ \delta_{jk}^h>0, \\
              0, & \mbox{otherwise}.
            \end{cases}$.
Further, define
$\delta:=\begin{cases}
   \min(\frac{\delta^+}{2^{H+1}+1},\frac{\epsilon}{2^{H+1}T}), & \mbox{if } \delta^+>0, \\
   \frac{\epsilon}{2^{H+1}T}, & \mbox{otherwise}.
 \end{cases}$.
We shall show that for $(j,k)\in \mathcal{P}_{h}(\bld{m}^\star) $ and $i\ge 1$,
\begin{equation}\label{in:bounds2}
 n l^n (m_{jk}^\star-2^{h+1}\delta) \le \calM_{ijk}^{p,n} \le
 n l^n (m_{jk}^\star+2^{h+1}\delta),
\end{equation}
then \eqref{in:bounds} follows by noticing that
 $2^{h+1}\delta \le 2^{h+1}\frac{\epsilon}{2^{H+1}T} \le\frac{\epsilon}{T} $.  We proceed in a two level induction: i) on the priority sets $\mathcal{P}_{h}(\bld{m}^\star) $, $h = 0, \dots, H+1$ and ii) on the review periods $i \in \{1,\ldots, \lfloor T/l^n \rfloor \}$. In particular, we first proceed in an induction on the priority sets. Then, for any set we proceed in a second induction on the review periods inside this set.

Let $(j,k) \in \mathcal{P}_{0}(\bld{m}^\star) $ and without loss of generality assume that $\lambda_j^D\le \lambda_k^S$. That is, $m_{jk}^\star=\lambda_j^D$ and $m_{jk'}^\star =0$ for $k'\neq k$. By definition of the first priority set, we have that for any $i\ge 1$,
$\calM_{ijk}^{p,n}= \min(Q^{D,n}_j(il^n-), Q^{S,n}_k(il^n-))$. Also, notice that node $j$ can only be connected to nodes $k' \neq k $ such that $m_{jk'}=0$ because otherwise $(j,k)$ cannot be a priority arc. In other words, $(j,k')\in \mathcal{P}_{H+1}(\bld{m}^\star)$ for $k' \neq k $.
First, using the same ideas as in the proof of Theorem~\ref{proposition:DRconvergence2}, we derive for large enough $n$ and each $i\ge 1$,
\begin{align}
Q^{D,n}_j(il^n-)\ge n l^n(\lambda_j^D -\delta ),\label{in:lowerQ}\\
Q_k^{S,n}(il^n-)\ge n l^n(\lambda^S_{k} - \delta)\label{in:lowerI}.
\end{align}
The latter inequalities yield a lower bound of the number of matches for each $i
\ge i$.
\begin{equation*}\label{in:larger}
\calM_{ijk}^{p,n} \ge n l^n (\min(\lambda_j^D,\lambda_k^S)-\delta)
= n l^n (m_{jk}^\star-\delta)
\ge n l^n (m_{jk}^\star-2\delta).
\end{equation*}

Next, we shall show that the opposite inequality holds as well.
By \eqref{eq:posQ}, \eqref{eq:posI},
Lemma~\ref{lemma:AtaKumar} choosing $\alpha=\delta$,
and observing that $R^{D,n}(il^n)-R^{D,n}((i-1)l^n)\ge 0$,
we have that for large $n$,
\begin{equation}\label{in:leftoverQ}
\begin{split}
Q^{D,n}_j(il^n-)&\le Q^{D,n}_j((i-1)l^n)+ A_j^{D,n}(il^n)-A_j^{D,n}((i-1)l^n)\\
&\le Q^{D,n}_j((i-1)l^n) + n l^n(\lambda_j^D +\delta)
\end{split}
\end{equation}
and
\begin{equation}\label{in:leftoverI}
\begin{split}
Q_k^{S,n}(il^n-)&\le Q_k^{S,n}((i-1)l^n) + A_k^{S,n}(il^n)- A_k^{S,n}((i-1)l^n)\\
&\le Q_k^{S,n}((i-1)l^n) + n l^n(\lambda^S_{k} + \delta).
\end{split}
\end{equation}
Furthermore, by \eqref{eq:posQ} and \eqref{eq:posI}, we get that
\begin{equation*}
\begin{split}
Q_j^{D,n}(il^n) &= Q_j^{D,n}(0) + A_j^{D,n}(il^n)-R_j^{D,n}(il^n)-\sum_{h=1}^{i} \sum_{k'\in \suppset_j} \calM_{hjk'}^{p,n}\\
&= Q_j^{D,n}(0) + A_j^{D,n}(il^n)-A_j^{D,n}((i-1)l^n)
-R_j^{D,n}(il^n)+R_j^{D,n}((i-1)l^n)\\
& \hspace{1cm} +A_j^{D,n}((i-1)l^n)-R_j^{D,n}((i-1)l^n)
-\sum_{h=1}^{i-1} \sum_{k'\in \suppset_j} \calM_{hjk'}^{p,n}-\sum_{k'\in \suppset_j} \calM_{ijk'}^{p,n}\\
&= Q_j^{D,n}((i-1)l^n)+A_j^{D,n}(il^n)-A_j^{D,n}((i-1)l^n)
-R_j^{D,n}(il^n)+R_j^{D,n}((i-1)l^n)-\sum_{k'\in \suppset_j} \calM_{ijk'}^{p,n}\\
 &=Q_j^{D,n}(il^n-)-\sum_{k'\in \suppset_j} \calM_{ijk'}^{p,n}.
\end{split}
\end{equation*}
Hence, we have that
\begin{equation}\label{eq:leftoverQ}
Q_j^{D,n}(il^n) = Q_j^{D,n}(il^n-) - \sum_{k'\in \suppset_j \setminus\{k\}}\calM_{ijk'}^{p,n}
- \min(Q_j^{D,n}(il^n-), Q_k^{S,n}(il^n-)).
\end{equation}
In a similar way, we obtain
\begin{equation}\label{eq:leftoverI}
Q_k^{S,n}(il^n) = Q_k^{S,n}(il^n-) - \sum_{j'\in \demset_k} \calM_{ij'k}^{p,n}.
\end{equation}
By the last equations, we have that for each $i\ge 1$ at least one of the two queue-lengths is zero. To see this, note that if
$\min(Q_j^{D,n}(il^n-), Q_k^{S,n}(il^n-))= Q_j^{D,n}(il^n-)$, then $\calM_{ijk'}^{p,n}=0$ for all $k'\neq k$ and hence
$Q_j^{D,n}(il^n)=0$. If $\min(Q_j^{D,n}(il^n-), Q_k^{S,n}(il^n-))= Q_k^{S,n}(il^n-)$, then
$\calM_{ijk}^{p,n}=Q_k^{S,n}(il^n-)$ and hence $\calM_{ij'k}^{p,n}=0$ for $j'\neq j$. That is, $Q_k^{S,n}(il^n)=0$ and $\min(Q_j^{D,n}(il^n), Q_k^{S,n}(il^n))=0$ for any $i\ge 1$. Further, by \eqref{in:leftoverQ} and \eqref{in:leftoverI}, we have that for large $n$,
\begin{equation}\label{in:upperbound0}
\begin{split}
\calM_{ijk}^{p,n} &= \min(Q_j^{D,n}(il^n-), Q_k^{S,n}(il^n-)) \\
&\le
\min\left(Q_j^{D,n}((i-1)l^n)+n l^n\lambda_j^D),
Q_k^{S,n}((i-1)l^n)+n l^n \lambda_k^S \right)+ n l^n\delta.
\end{split}
\end{equation}
Further, by the assumptions for the initial conditions, we have that
$\frac{ Q^{D,n}_j(0)}{nl^n} \to 0$ and $\frac{ Q_k^{S,n}(0)}{nl^n} \to 0$
as $n\to \infty$. That is,
 $Q^{D,n}_j(0) \le nl^n \delta$ and  $Q_k^{S,n}(0) \le nl^n \delta$.
If $\lambda_j^D=\lambda_k^S$, then \eqref{in:upperbound0} yields
\begin{equation*}
\begin{split}
M_{1jk}^{p,n} &\le
\min(Q_j^{D,n}(0),Q^{S,n}_k(0)) +n l^n(\lambda_j^D+\delta)\le
n l^n(\lambda_j^D+2\delta)= n l^n(m_{jk}^\star+2\delta)
\end{split}
\end{equation*}
and for $i>1$,
\begin{equation*}
\begin{split}
M_{ijk}^{p,n} &\le
\min(Q_j^{D,n}((i-1)l^n),Q_k^{S,n}((i-1)l^n)) +n l^n(\lambda_j^D+\delta)\le
n l^n(\lambda_j^D+\delta) \le  n l^n(m_{jk}^\star+2\delta),
\end{split}
\end{equation*}
because $\min(Q_j^{D,n}((i-1)l^n),Q_k^{S,n}((i-1)l^n))=0.$
Let $\lambda_j^D<\lambda_k^S$. We shall use induction in the review period $i$ to show that $Q_j^{D,n}(il^n-)\le Q_k^{S,n}(il^n-)$ and $Q_j^{D,n}(il^n)=0$ for each $i\ge 1$ and for large enough $n$.
We start by proving this for $i=1$ (i.e., the first step of the induction). Replacing $i=1$ in \eqref{in:leftoverQ} and \eqref{in:leftoverI},
it follows that
$Q^{D,n}_j(l^n-)\le  n l^n(\lambda_j^D + 2\delta)$
and
$Q_k^{S,n}(l^n-)\le  n l^n(\lambda^S_{k} + 2\delta)$.
By the last inequalities, \eqref{in:lowerQ}, \eqref{in:lowerI}, and observing that
$ \delta \le \frac{\lambda_k^S-\lambda_j^D}{3}$, we get
\begin{equation*}
Q^{D,n}_j(l^n-)\le n l^n(\lambda_j^D+2\delta) \le n l^n(\lambda_k^S-\delta) \le Q_k^{S,n}(l^n-),
\end{equation*}
and by \eqref{eq:leftoverQ}, $Q^{D,n}_j(l^n)=0$. Replacing now $i=2$ in \eqref{in:leftoverQ}, we obtain
$Q^{D,n}_j(2l^n-)\le  n l^n(\lambda_j^D +2 \delta)$. Using again  \eqref{in:lowerI} and \eqref{eq:leftoverQ}, we have that
\begin{equation*}
Q^{D,n}_j(2l^n-)\le n l^n(\lambda_j^D+2\delta) \le n l^n(\lambda_k^S-\delta) \le Q_k^{S,n}(2l^n-)
\end{equation*}
and $Q^{D,n}_j(2l^n)=0$. Let now $Q^{D,n}_j(zl^n-)\le Q_k^{S,n}(zl^n-)$ and $Q^{D,n}_j(zl^n)=0$ for all $1<z\le i$ for some $i>1$ be the hypothesis induction (for the induction in $i$). We shall show that the same holds for the review period $i+1$ as well.
 Using exactly the same steps as above, we conclude that $Q^{D,n}_j((i+1)l^n-)\le Q_k^{S,n}((i+1)l^n-)$ and $Q^{D,n}_j((i+1)l^n)=0$. That is, by \eqref{in:upperbound0}, we have that for each $i\ge 1$, $(j,k) \in \mathcal{P}_{0}(\bld{m}^\star) $, and large enough $n$,
\begin{equation}\label{in:smaller}
\calM_{ijk}^{p,n} \le n l^n( \lambda_j^D+2\delta)\le n l^n (m_{jk}^\star+2\delta)
\end{equation}
and hence \eqref{in:bounds2} holds.

Now, we make the hypothesis induction for the induction in priority sets.  Assume that \eqref{in:bounds2} holds
for $(j,k)\in\mathcal{P}_{r}(\bld{m}^\star) $ for all $ r= 0, \dots, h-1$, with
$h< H+1$. We shall show \eqref{in:bounds2} holds for $h$ as well. For $(j,k)\in\mathcal{P}_{h}(\bld{m}^\star) $, by definition of the priority ordering, we have that
\begin{align*}
\calM_{ijk}^{p,n} =
\min\left(Q^{D,n}_j(il^n-)-\sum_{k':(j,k')\in \mathcal{Q}_{h-1}(\bld{m}^\star)}
\calM_{ijk'}^{p,n},
Q^{S,n}_k(il^n-)-\sum_{j':(j',k)\in \mathcal{Q}_{h-1}(\bld{m}^\star)}
\calM_{ij'k}^{p,n}
\right).
\end{align*}
Without loss of generality assume that $(j,k)$ makes node $j$ tight. As in the first induction step, we observe that node $j$ can only be connected through edges
$(j,k') \notin \mathcal{Q}_{h-1}(\bld{m}^\star) $ such that  $m_{jk'}=0$ because otherwise $(j,k)$ cannot be a priority arc. In other words, $(j,k')\in \mathcal{P}_{H+1}(\bld{m}^\star)$. Furthermore, the optimal solution to the MP becomes
	\begin{align*}
	m_{jk}^\star &=
	\min\left( \lambda_j^D - \sum_{k':(j,k') \in \mathcal{Q}_{h-1}(\bld{m}^\star)}m^\star_{jk'}, \lambda_k^S - \sum_{j':(j',k)\in \mathcal{Q}_{h-1}(\bld{m}')}m^\star_{j'k} \right) \\
&= \lambda_j^D - \sum_{k':(j,k') \in \mathcal{Q}_{h-1}(\bld{m}^\star)}m^\star_{jk'}.
	\end{align*}
By \eqref{in:bounds2} and \eqref{in:lowerQ}, we derive for large enough $n$ and each $i\ge 1$,
\begin{equation}\label{in:lowerQG}
\begin{split}
Q^{D,n}_j(il^n-)-\sum_{k':(j,k')\in \mathcal{Q}_{h-1}(\bld{m}^\star)}
\calM_{ijk'}^{p,n}&\ge
n l^n(\lambda_j^D-\delta-
\sum_{k':(j,k') \in \mathcal{Q}_{h-1}(\bld{m}^\star)}m^\star_{jk'}-\sum_{z=0}^{h-1}2^{z+1}\delta)\\
&\ge
n l^n(\lambda_j^D-\sum_{k':(j,k') \in \mathcal{Q}_{h-1}(\bld{m}^\star)}m^\star_{jk'}
-(2^{h+1}-1)\delta),\\
\end{split}
\end{equation}
where we use that $\sum_{z=0}^{h-1}2^{z+1}=2^{h+1}-2$.
In a similar way,
\begin{equation}\label{in:lowerIG}
\begin{split}
Q_k^{S,n}(il^n-)-\sum_{j':(j',k)\in \mathcal{Q}_{h-1}(\bld{m}^\star)}
\calM_{ij'k}^{p,n}
\ge
n l^n(\lambda^S_{k}-\sum_{j':(j',k) \in \mathcal{Q}_{h-1}(\bld{m}^\star)}m^\star_{j'k}
-(2^{h+1}-1)\delta).
\end{split}
\end{equation}
The last two inequalities yield a lower bound of the number of matches for each
$i\ge 1$,
\begin{equation*}\label{in:largerG}
\calM_{ijk}^{p,n} \ge n l^n (m_{jk}^\star-(2^{h+1}-1)\delta)
\ge n l^n (m_{jk}^\star-2^{h+1}\delta).
\end{equation*}

We now move in proving the opposite inequality.
First observe that \eqref{in:leftoverQ} and \eqref{in:leftoverI} continue to hold.
Using a similar way as in the first induction step, we obtain
\begin{equation}\label{eq:leftoverQG}
Q^{D,n}_j(il^n) = Q^{D,n}_j(il^n-)
-\sum_{k':(j,k')\in \mathcal{Q}_{h-1}(\bld{m}^\star)}
\calM_{ijk'}^{p,n}
-\calM_{ijk}^{p,n}
\end{equation}
and
 \begin{equation}\label{eq:leftoverIG}
Q^{S,n}_k(il^n) = Q^{S,n}_k(il^n-)
-\sum_{j':(j',k)\in \mathcal{Q}_{h-1}(\bld{m}^\star)}
\calM_{ij'k}^{p,n}
-\calM_{ijk}^{p,n}.
\end{equation}
Further, by \eqref{eq:leftoverQG}, \eqref{eq:leftoverIG}, and the  definition of the priority ordering, we have that
$\min(Q^{D,n}_j(il^n), Q^{S,n}_k(il^n))=0$ for any $i\ge 1$. By \eqref{in:bounds2}, \eqref{in:leftoverQ}, and  \eqref{in:leftoverI}, we get the following inequality for large enough $n$,
 \begin{equation}\label{in:upperbound0G}
\begin{split}
\calM_{ijk}^{p,n} &=
\min\Bigg(Q^{D,n}_j(il^n-)-\sum_{k':(j,k')\in \mathcal{Q}_{h-1}(\bld{m}^\star)}
\calM_{ijk'}^{p,n},
Q^{S,n}_k(il^n-)-\sum_{j':(j',k)\in \mathcal{Q}_{h-1}(\bld{m}^\star)}
\calM_{ij'k}^{p,n}
\Bigg) \\
&\le
\min\Bigg(Q^{D,n}_j((i-1)l^n)+n l^n(\lambda_j^D+\delta)
-n l^n(\sum_{k':(j,k')\in \mathcal{Q}_{h-1}(\bld{m}^\star)} m_{jk'}^\star- (2^{h+1}-2)\delta),\\
&\hspace{1.5cm}Q^{S,n}_k((i-1)l^n)+n l^n(\lambda_k^S+\delta)
-n l^n(\sum_{j':(j',k)\in \mathcal{Q}_{h-1}(\bld{m}^\star)} m_{j'k}^\star-(2^{h+1}-2)\delta)\Bigg)\\
&\le
\min\Bigg(Q^{D,n}_j((i-1)l^n)+n l^n(\lambda_j^D- \sum_{k':(j,k')\in \mathcal{Q}_{h-1}(\bld{m}^\star)} m_{jk'}^\star),\\
&\hspace{1.5cm}Q^{S,n}_k((i-1)l^n)+n l^n(\lambda_k^S-\sum_{j':(j',k)\in \mathcal{Q}_{h-1}(\bld{m}^\star)} m_{j'k}^\star)\Bigg)+ n l^n(2^{h+1}-1)\delta.
\end{split}
\end{equation}
If
$\lambda_j^D- \sum_{k':(j,k')\in \mathcal{Q}_{h-1}(\bld{m}^\star)} m_{jk'}^\star
=\lambda_k^S-\sum_{j':(j',k)\in \mathcal{Q}_{h-1}(\bld{m}^\star)} m_{j'k}^\star$,
 then \eqref{in:upperbound0G} and the assumptions for the initial conditions yield
\begin{equation*}
\begin{split}
M_{ijk}^{p,n} &\le
\min(Q^{D,n}_j((i-1)l^n),Q^{S,n}_k((i-1)l^n)) +n l^n(m_{jk}^\star+2^{h+1}\delta)
=n l^n(m_{jk}^\star+2^{h+1}\delta).
\end{split}
\end{equation*}
Let $\lambda_j^D- \sum_{k':(j,k')\in \mathcal{Q}_{h-1}(\bld{m}^\star)} m_{jk'}^\star
<\lambda_k^S-\sum_{j':(j',k)\in \mathcal{Q}_{h-1}(\bld{m}^\star)} m_{j'k}^\star$. We shall use the same ideas as in the first induction step to show the upper bound in \eqref{in:bounds2}.
Assume that
$$Q^{D,n}_j(zl^n-)-\sum_{k':(j,k')\in \mathcal{Q}_{h-1}(\bld{m}^\star)}
\calM_{zjk'}^{p,n}\le Q^{S,n}_k(zl^n-)
-\sum_{j':(j',k)\in \mathcal{Q}_{h-1}(\bld{m}^\star)}
\calM_{zj'k}^{p,n}$$
and $Q^{D,n}_j(zl^n)=0$ for all
$z\le i$ for some $i> 1$.
By \eqref{in:bounds2} for $(j,k')\in\mathcal{Q}_{h-1}(\bld{m}^\star)$, \eqref{in:leftoverQ}, and \eqref{in:lowerIG}, and observing that
$\delta\le \frac{\lambda_k^S- \sum_{j':(j',k)\in \mathcal{Q}_{h-1}(\bld{m}^\star)} m_{j'k}^\star
-\lambda_j^D+\sum_{k':(j,k')\in \mathcal{Q}_{h-1}(\bld{m}^\star)} m_{jk'}^\star}{2^{h+1}+1},$
we have that
\begin{equation*}
\begin{split}
Q^{D,n}_j((i+1)l^n-)-\sum_{k':(j,k')\in \mathcal{Q}_{h-1}(\bld{m}^\star)}
\calM_{i+1jk'}^{p,n}
&\le
n l^n(\lambda_j^D- \sum_{k':(j,k')\in \mathcal{Q}_{h-1}(\bld{m}^\star)}m_{jk'}^\star+2^{h+1}\delta)\\
&\le
 n l^n(\lambda_k^S-\sum_{j':(j',k)\in \mathcal{Q}_{h-1}(\bld{m}^\star)} m_{j'k}^\star-(2^{h+1}-1)\delta)\\
&\le Q^{S,n}_k((i+1)l^n-)
 -\sum_{j':(j',k)\in \mathcal{Q}_{h-1}(\bld{m}^\star)}
\calM_{i+1j'k}^{p,n},
\end{split}
\end{equation*}
and hence $Q^{D,n}_j((i+1)l^n)=0$ by \eqref{eq:leftoverQG}.
That is, by \eqref{in:upperbound0G}  for each $i\ge 1$, $(j,k) \in \mathcal{P}_{h}(\bld{m}^\star) $, and large enough $n$,
\begin{equation}\label{in:smaller}
\calM_{ijk}^{p,n} \le n l^n( \lambda_j^D-\sum_{k':(j,k')\in \mathcal{Q}_{h-1}(\bld{m}^\star)}m_{jk'}^\star+2^{h+1}\delta)
= n l^n (m_{jk}^\star+2^{h+1}\delta)
\end{equation}
and hence \eqref{in:bounds2} holds.

Last, it is clear from the previous steps that \eqref{in:bounds2} holds for
$(j,k) \in \mathcal{P}_{H+1}(\bld{m}^\star) $ since $\calM_{ijk}^{p,n}=m_{jk}^\star=0$ for large $n$. That is,  \eqref{in:bounds2} holds for any $(j,k) \in \mathcal{P}_{h}(\bld{m}^\star) $ with $h=0,\dots,H+1$.
\eof

\subsection{Proofs for Section~\ref{sec:PSB}}
\label{Proofs sec:blindOptimal}
Before we move to the main part of the proof of Theorem~\ref{theorem:blindOptimal} and Proposition~\ref{theorem:convergenceLPset}, we show some helpful properties of the value of matching problem \eqref{eq:max_A_w_v3}. For these proofs, let $\mathbf{e}_j$ denote the $j^{th}$ unit vector and $\mathbf{e}$ denote the vector of all ones.

Let $\bar{\p}:=\max_{(j,k)\in\arcset}(\p_{jk})$. For two vectors $\bld{x}\in \R^J$ and $\bld{b}\in \R^K$, let $m^\star_{jk}(\bld{x},\bld{b})$ be an optimal solution to \eqref{eq:max_A_w_v3} if we replace $\bld{\lambda}^D$ and  $\bld{\lambda}^S$ by $\bld{x}$ and $\bld{b}$, respectively.
Define the value of the matching problem
\begin{align*}
F(\bld{x},\bld{b}) := \sum_{(j,k)\in\arcset} \p_{jk} m^\star_{jk}(\bld{x},\bld{b}).
\end{align*}
\
\begin{lemma}\label{lemma:LC}
Let $\bld{\lambda}^D\in \R_+^J$ and $\bld{\lambda}^S\in \R_+^K$. There exists
a Lipschitz continuous mapping
$\bld{m}^\star:\R_+^J \times \R_+^K \rightarrow \R_+^{J} \times \R_+^K$ such that $\bld{m}^\star(\bld{\lambda}^D,\bld{\lambda}^S)$ is an optimal solution to  \eqref{eq:max_A_w_v3}\footnote{
Lemma~\ref{lemma:LC} implies that there exists an optimal solution to \eqref{eq:max_A_w_v3} that is Lipschitz continuous. The later holds for any norm $\|\cdot\|$ in $\R_+^J \times \R_+^K$ as all norms in this space (and in general in any vector space with finite dimension) are equivalent. In the sequel, we freely use the norm that suits our approach.}.
\end{lemma}
\pfa{Proof}
We begin by observing that the feasible region of \eqref{eq:max_A_w_v3} is
nonempty, closed, and convex. The remainder of the proof follows by using the Lipschitz selection theorem
\cite[Theorem~9.4.3]{aubin} and \cite[Theorem~10.5]{schrijver} exactly as in \cite[Proposition~2]{bassamboo2006}.
\eof

\begin{lemma}\label{lem:prof_bnd}
The value of \eqref{eq:max_A_w_v3} is nondecreasing in each of its arguments and the following inequalities hold, for any $\alpha \ge 0$,
	\begin{align*}
	& F(\bld{x}+\alpha\mathbf{e}_j,\bld{b})
- F(\bld{x},\bld{b}) \in [0,\bar{\p} JK C \alpha], \ \forall j\in\bbJ, \\
	& F(\bld{x},\bld{b}+\alpha\mathbf{e}_k) - F(\bld{\lambda}^D,\bld{\lambda}^S) \in [0,\bar{\p} JK C \alpha], \ \forall k\in\bbK,
	\end{align*}
where the constant $C$ is the Lipschitz constant in Lemma~\ref{lemma:LC}.
\end{lemma}
\pfa{Proof}
Clearly, $F(\cdot,\cdot)$ is nondecreasing in each of its arguments by the definition of the MP. We shall show that the rate of change is bounded. To this end, observe that by Lemma~\ref{lemma:LC} there exists a Lipschitz continuous optimal solution to  \eqref{eq:max_A_w_v3} such that
\begin{equation*}
\left\| \bld{m}^\star(\bld{x}+\alpha\mathbf{e}_j,\bld{b})
-\bld{m}^\star(\bld{x},\bld{b})\right\|_\infty
 \leq
 C \left\| (\bld{x}+\alpha\mathbf{e}_j,\bld{b})
 -(\bld{x},\bld{b})\right\|_\infty \leq C\alpha.
\end{equation*}
The latter leads to
\begin{equation*}
\left| m^\star_{jk}(\bld{x}+\alpha\mathbf{e}_j,\bld{b})
-m^\star_{jk}(\bld{x},\bld{b}) \right|
\le C\alpha,
\end{equation*}
for any $(j,k)\in\arcset$.
Now, using the last inequality, we have that
\begin{equation*}
  \begin{split}
0 \le F(\bld{x}+\alpha\mathbf{e}_j,\bld{b})
- F(\bld{x},\bld{b}) &=
\sum_{(j,k)\in\arcset} \p_{jk}
(m^\star_{jk}(\bld{x}+\alpha\mathbf{e}_j,\bld{b})-m^\star_{jk}(\bld{x},\bld{b}))\\
& \le \sum_{(j,k)\in\arcset} \p_{jk}
\left|m^\star_{jk}(\bld{x}+\alpha\mathbf{e}_j,\bld{b})-m^\star_{jk}(\bld{x},\bld{b})\right|\\
&\le \bar{\p} JK C \alpha.
  \end{split}
\end{equation*}
The second inequality follows analogously. Note that both inequalities do not depend on the particular optimal solution $\bld{m}^\star(\cdot,\cdot)$, that is Lipschitz continuous, as any solution achieves the same value.
\eof

\begin{lemma}\label{lemma:monotoneBl}
	For any $\epsilon > 0$, and any $\mathbf{A} \ge \bld{x}-\epsilon\mathbf{e}$ and $\mathbf{B} \ge \bld{b}-\epsilon\mathbf{e}$,
	\begin{align*}
	& F(\mathbf{A},\mathbf{B}) \ge F(\bld{x},\bld{b})
      - C_1 \epsilon,
	\end{align*}
where $C_1:=(J+K) JK C \bar{\p}$.
\end{lemma}
\pfa{Proof}
	First, by a telescoping sum (and assuming an empty sum evaluates to zero), we observe that
	\begin{align*}
	& F(\bld{x},\bld{b}) - F(\bld{x}-\epsilon\mathbf{e},\bld{b}-\epsilon\mathbf{e}) \\
	& =  \sum_{j = 1}^J \Big[F(\bld{x}-\epsilon \sum_{s=1}^{j-1}\mathbf{e}_s , \bld{b}) - F(\bld{x}-\epsilon\sum_{s=1}^j\mathbf{e}_s,\bld{b})\Big] + \sum_{k = 1}^K \Big[F(\bld{x} -\epsilon\mathbf{e}, \bld{b} -\epsilon \sum_{t=1}^{k-1}\mathbf{e}_t) - F(\bld{x} -\epsilon\mathbf{e}, \bld{b}-\epsilon\sum_{t=1}^k\mathbf{e}_t)\Big] \\
	& \le (J+K) JK C \bar{\p}\epsilon,
	\end{align*} where the inequality follows from the upper bound of Lemma \ref{lem:prof_bnd} for each incremental change in the objective value. Then, since the lower bound in Lemma \ref{lem:prof_bnd} implies $F(\cdot,\cdot)$ is monotone nondecreasing, we have
	\begin{align*}
	F(\mathbf{A},\mathbf{B}) \ge F(\bld{x}-\epsilon\mathbf{e},\bld{b}-\epsilon\mathbf{e}) \ge F(\bld{x},\bld{b}) - (J+K) JK C \bar{\p}\epsilon,
	\end{align*} completing the proof.
\eof

\pfa{Proof \amy{of Lemma~\ref{lem:unimodular}}}
Totally unimodular matrices with elements $a_{il}$ are characterized by \cite{tamir1976} as those where every subset of rows $R$ can be partitioned into two subsets, $R_1\cup R_2 = R$, such that $\sum_{i \in R_1} a_{il} - \sum_{i \in R_2} a_{il} \in \{-1,0,1\}$ for all columns $l$. The constraint matrix of \eqref{LP:MatchingPo} has a row for each $k \in \bbK$ and each $j \in \bbJ$ and a column for every variable $y_{jk}$ for $(j,k)\in\arcset$. For a given $j \in \bbJ$ representing a row and $ (j',k')\in\arcset$ representing a column, the corresponding element of the constraint matrix is given by
	\begin{align*}
	a_{(j);(j'k')} =
	\begin{cases}
	1, & j' = j \\
	0, & j' \ne j
	\end{cases},
	\end{align*} representing that only those $y_{j'k'}$ with $j' = j$ are included in the sum for the constraint corresponding to $j$. Similarly, for a given $k \in \bbK$ and $ (j',k')\in\arcset$, the corresponding constraint matrix element is
	\begin{align*}
	a_{(k);(j'k')} =
	\begin{cases}
	1, & k' = k \\
	0, & k' \ne k
	\end{cases}.
	\end{align*}
	Any subset of the rows of the constraint matrix corresponds to the union of some subsets of $\bbJ$ and $\bbK$, i.e. $J' \cup K'$ where $J' \subseteq \bbJ$ and $K' \subseteq \bbK$. Therefore, given a subset of rows $J'\cup K'$, define the required partition as $R_1 = J'$ and $R_2 = K'$. Then, for a given variable $y_{j'k'}$,
	\begin{align*}
	\sum_{j \in J'} a_{(j);(j'k')} - \sum_{k \in K'} a_{(k);(j'k')} =
	\begin{cases}
	-1, & j' \notin J' \text{ and } k' \in K', \\
	0, &  (j' \in K' \text{ and } k' \in K') \text{ or } (j' \notin J' \text{ and } k \notin K'), \\
	1, &  j' \in J' \text{ and } k' \notin K', \\
	\end{cases}
	\end{align*}
	This completes the proof that the constraint matrix of \eqref{LP:MatchingPo} is totally unimodular. Then, by \cite{tamir1976}, if the right hand side constraints of \eqref{LP:MatchingPo} are integer valued, any optimal extreme point solution to \eqref{LP:MatchingPo} is also integer valued. A similar unimodularity property is proven in \cite{DPW2019}.
\eof \\

\pfa{Proof of Theorem~\ref{theorem:blindOptimal}}
\amy{We first prove part (i), and then (ii).} \\

\amy{\noindent {\bf Proof of (i):}
The assumption \eqref{eq:ass-M-zero-HC} and the fact that $\bld{m}^\star$ is an optimal solution to
\eqref{eq:max_A_w_v3G} imply \eqref{eq:max-achievable-profit-fluid}.
} \\

\amy{\noindent {\bf Proof of (ii):}}
Let $\epsilon>0$ and $\alpha=\frac{l \epsilon}{C_1T }$, where $C_1$ is given in Lemma~\ref{lemma:monotoneBl}.
\amy{The same arguments to show \eqref{eq:BoundRD2} and \eqref{eq:BoundRS2} in Theorem~\ref{proposition:DRconvergence2} show that we can} choose large enough $n$,
such that
 \begin{equation}\label{eq:BoundRD3}
\begin{split}
 \frac{1}{nl^n }\sum_{m=A^{D,n}_j((i-1)l^n)+1}^{A^{D,n}_j(il^n)} \ind{r_{jm}^D\le l^n}
  <
  \frac{\epsilon}{C_1T}
\end{split}
\end{equation}
and
 \begin{equation}\label{eq:BoundRS4}
\begin{split}
 \frac{1}{nl^n }\sum_{m=A^{S,n}_k((i-1)l^n)+1}^{A^{S,n}_k(il^n)} \ind{r_{km}^S\le l^n}
  <
 \frac{\epsilon}{C_1T}.
\end{split}
\end{equation}
Define $\Omega_4^n$ the events such that \eqref{eq:BoundRD3} and
\eqref{eq:BoundRS4} hold. In the sequel, we take \amy{a sample path} $\omega \in \Omega_1^n \cap \Omega_2^n  \cap \Omega_4^n$, \amy{where $\Omega_1^n$ and $\Omega_2^n$ are as defined in the proof of Theorem~\ref{proposition:DRconvergence2}}.

The next step is to show  the following inequality  for large enough $n$,
\begin{equation}\label{in:ineF}
\amy{\frac{1}{n}}V_{\bld{M}^{\levi{LP},n}}(t)\ge \sum_{\amy{(j,k) \in \arcset}}  \p_{jk} m_{jk}^{\star}t-2\epsilon.
\end{equation}
Following similar steps as in  proof of Theorem~\ref{proposition:DRconvergence2}, we obtain for
$i \in \{1,\ldots, \lfloor t/l \rfloor \}$,
\begin{equation*}
\frac{Q^{D,n}_j(il^n-)}{nl^n}\ge \lambda_j^D
-\frac{\epsilon}{C_1t}
\end{equation*}
and
\begin{equation*}
\frac{Q_k^{S,n}(il^n-)}{nl^n}\ge \lambda^S_{k}
-\frac{\epsilon}{C_1t}.
\end{equation*}
Applying now Lemma~\ref{lemma:monotoneBl}, we have that for
$i \in \{1,\ldots, \lfloor t/l^n \rfloor \}$,
\begin{equation*}\label{in:boundF}
F\left(\frac{\bld{Q}^{D,n}(il^n-)}{nl^n},
\frac{\bld{Q}^{S,n}(il^n-)}{nl^n}\right)
\ge
F(\bld{\lambda}^D,\bld{\lambda}^S) - \frac{\epsilon}{t}.
\end{equation*}
By the properties of linear programming, we have that $F(\cdot,\cdot)$ is homogenous of degree one. That is,
\begin{equation*}\label{in:boundF}
F\left(\frac{\bld{Q}^{D,n}(il^n-)}{nl^n},
\frac{\bld{Q}^{S,n}(il^n-)}{nl^n}\right)
=\frac{1}{n l^n}
F\left(\bld{Q}^{D,n}(il^n-),\bld{Q}^{S,n}(il^n-)\right).
\end{equation*}
Now observe that the optimization problem  \eqref{LP:MatchingPo}, which decides the matches for the LP-based policy, and the MP \eqref{eq:max_A_w_v3} have the same objective function and the same constraint matrix, which is totally unimodular by Lemma~\ref{lem:unimodular}. This implies the well known fact that these programs have the same optimal value if they have the same right hand side constraints, and if these constraints are integer valued (i.e., the objective value of the integer program \eqref{LP:MatchingPo} is equal to the objective value of its LP relaxation \eqref{eq:max_A_w_v3}; see for example \cite[Theorem 13.2]{papadimitriou1998}. Combining all the above together, we derive
\begin{equation*}
\begin{split}
 \sum_{(j,k) \in \arcset} \p_{jk} \frac{M_{jk}^{\levi{LP},n}(t)}{n} &
=
\sum_{i=1}^{\lfloor t/l^n \rfloor} \frac{F\left(\bld{Q}^{D,n}(il^n-),\bld{Q}^{S,n}(il^n-)\right)}{n}\\
&=
\sum_{i=1}^{\lfloor t/l^n \rfloor} l^n F\left(\frac{\bld{Q}^{D,n}(il^n-)}{nl^n},\frac{\bld{Q}^{S,n}(il^n-)}{nl^n}\right)\\
&\ge
l^n \lfloor t/l^n \rfloor (F(\bld{\lambda}^D,\bld{\lambda}^S) -\frac{\epsilon}{t})\\
&\ge
 t F(\bld{\lambda}^D,\bld{\lambda}^S)-F(\bld{\lambda}^D,\bld{\lambda}^S)l^n - \epsilon .
\end{split}
\end{equation*}
By definition of $F(\cdot,\cdot)$ and for $n$ such that
$n^{2/3}>F(\bld{\lambda}^D,\bld{\lambda}^S)l/\epsilon$, we have that
\begin{equation*}
\begin{split}
\frac{V_{\bld{M}^{\levi{LP},n}}(t)}{n}=\sum_{(j,k) \in \arcset} \p_{jk} \frac{M_{jk}^{\levi{LP},n}(t)}{n}
\ge
\sum_{\amy{(j,k) \in \arcset}}  \p_{jk} m_{jk}^{\star}t - 2\epsilon,
\end{split}
\end{equation*}
and
hence \eqref{in:ineF} follows.

 \amy{Finally, since both \eqref{in:ineF} and \eqref{eq:max-achievable-profit-fluid} hold for any sample path $\omega \in \Omega_1^n \cap \Omega_2^n \cap \Omega_4^n$ and $\lim_{n\rightarrow \infty}
\Prob{\Omega_1^n \cap \Omega_2^n  \cap \Omega_4^n }=1$, letting $\epsilon \rightarrow 0$ completes the proof.}

\eof \\

\pfa{Proof of Proposition~\ref{theorem:convergenceLPset}}
Without loss of generality assume that $t=1$ otherwise fix a $t>0$ and appropriately scale  \eqref{eq:max_A_w_v3} by $t$. First, we note that the set of optimal solutions $S$ is a closed set. To see this observe that if it is not, then we can always increase one of the components of $\bld{y}^{\star}$ and achieve a better optimal value.

\levi{
 Fix a sample path such that Theorem~\ref{theorem:blindOptimal} holds. First, observe that for all  $j \in \bbJ$ and $k \in \bbK$, we have that
 $\frac{\bld{M}^{\levi{LP},n}(1)}{n} \leq \frac{\bld{A}^{D,n}(1)}{n}$ and $\frac{\bld{A}^{D,n}(1)}{n} \rightarrow \bld{\lambda}^D$ almost surely. In other words, $\frac{\bld{M}^{\levi{LP},n}(1)}{n}$ is a bounded sequence in $\R^{J\times K}$. Now, by
 Bolzano–Weierstrass theorem there exists a convergent subsequence, which we denote by $\frac{\bld{M}^{\levi{LP},n_l}(1)}{n_l} $, such that
 $\frac{\bld{M}^{\levi{LP},n_l}(1)}{n_l}  \rightarrow \bld{x}$. If for all such convergent subsequences the subsequential limits lie in the closed set $S$ (i.e., $\bld{x} \in S$), then the proof follows.

 So we proceed by contradiction: assume there exists a convergent subsequence such that $\frac{\bld{M}^{\levi{LP},n_l}(1)}{n_l}  \rightarrow \bld{x} \notin S$. }
%
Then, note that by the admissibility of the LP-based policy any limit should be feasible to \eqref{eq:max_A_w_v3}. To see this, observe that for all $j \in \bbJ$ and $k \in \bbK$ the following hold almost surely by the nonnegativity of the queue-lengths
\begin{equation*}
\sum_{k \in \suppset_j}
 \frac{M_{jk}^{\levi{LP},n_l}(1)}{n_l} \le  \frac{A_j^{D,n_l}(1)}{n_l }
\mbox{ and  }
\sum_{j\in \demset_k}
 \frac{M_{jk}^{\levi{LP},n_l}(1)}{n_l } \le  \frac{A_k^{S,n_l}(1)}{n_l },
\end{equation*}
and by law of large numbers
$\lim_{n_l\rightarrow \infty} \frac{A_j^{D,n_l}(1)}{n_l }= \lambda_j^D$,
$\lim_{n_l\rightarrow \infty} \frac{A_k^{S,n_l}(1)}{n_l }= \lambda_k^S$.
Now, we have that
\begin{equation*}
  \sum_{(j,k) \in \arcset}
 \p_{jk} \frac{M_{jk}^{\levi{LP},n_l}(1)}{n_l}
 \rightarrow
 \sum_{(j,k) \in \arcset}
 \p_{jk} x_{jk},
\end{equation*}
and by Theorem~\ref{theorem:blindOptimal}
$ \sum_{(j,k) \in \arcset}
 \p_{jk} x_{jk}=
 \sum_{(j,k) \in \arcset}
 \p_{jk} m_{jk}^{\star}  $. In other words, the vectors $\bld{x}$ and $\bld{m}^{\star} $ achieve the same optimal value of  \eqref{eq:max_A_w_v3} and so $\bld{x}$ is an optimal solution to  \eqref{eq:max_A_w_v3}. On the other hand,
 $\bld{x} \notin S$ which yields a contradiction.
\eof

\section{Exact calculations of the mean queue-lengths in Figure~\ref{figure:ExpQueues}}\label{sec:markovChain}
For convenience in the notation of this appendix, we let $\lambda = \lambda^D$ denote the demand arrival rate, and $\mu = \lambda^S$ denote the supply arrival rate. Let $X(\infty)$ be a Markov chain with transition rates are shown in Figure~\ref{figure:ExpQueues}.
\begin{figure}[H]
\centering
  \includegraphics[width=0.8\textwidth]{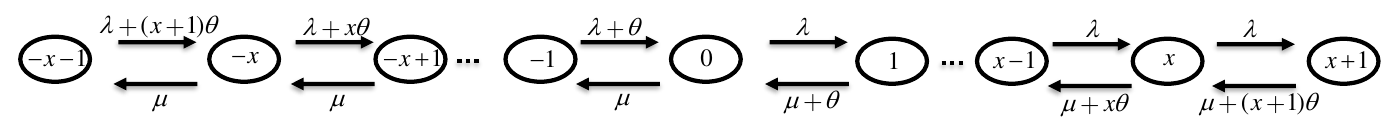}
\caption{Transition rates of the Markov chain.}
\label{figure:transition}
\end{figure}
Its probability distribution $\pi(\cdot)$ is the (unique) solution  of the balance equations that are written as follows:  for $x>0$
\begin{equation*}
(\lambda+\mu+x \theta) \pi(x) = \lambda \pi(x-1) + (\mu+(x+1) \theta) \pi(x+1),
\end{equation*}
\begin{equation*}
(\lambda+\mu+x \theta) \pi(-x) = \mu \pi(-x+1) + (\lambda+(x+1) \theta) \pi(-x-1),
\end{equation*}
and $x=0$,
\begin{equation*}
(\lambda+\mu) \pi(0) = (\lambda+\theta) \pi(-1)+ (\mu+\theta) \pi(1).
\end{equation*}
The solution to the balance equations in given by
\begin{equation}\label{eq:solutionMarkov}
 \pi(x) =
\begin{cases}
  \frac{\lambda^x}{\prod_{j=1}^{x}(\mu+j \theta)} \pi(0) & \mbox{if } x>0, \\
  \frac{\mu^{-x}}{\prod_{j=1}^{-x}(\lambda+j \theta)} \pi(0) & \mbox{if }x<0,
\end{cases}
\end{equation}
where
\begin{equation*}
\pi(0) =
\left(1+ \sum_{x=1}^{\infty}\frac{\lambda^x}{\prod_{j=1}^{x}(\mu+j \theta)}
+ \sum_{x=1}^{\infty}\frac{\mu^x}{\prod_{j=1}^{x}(\lambda+j \theta)} \right)^{-1}.
\end{equation*}
To see this, observe that for $x>0$,
\begin{equation*}
\begin{split}
\lambda \pi(x-1) + (\mu+(x+1) \theta) \pi(x+1) & =
\left( \lambda \frac{\lambda^{x-1}}{\prod_{j=1}^{x-1}(\mu+j \theta)} +
(\mu+(x+1) \theta) \frac{\lambda^{x+1}}{\prod_{j=1}^{x+1}(\mu+j \theta)}
\right)
\pi(0)\\
& =
\left((\mu+x \theta)  \frac{\lambda^{x}}{\prod_{j=1}^{x}(\mu+j \theta)} +
 \lambda \frac{\lambda^{x}}{\prod_{j=1}^{x}(\mu+j \theta)}
\right)
\pi(0)\\
& = ( \lambda + \mu+x \theta ) \pi(x).
\end{split}
\end{equation*}
In a similar way, for $x>0$ we have that
\begin{equation*}
\begin{split}
\mu \pi(-x+1) + (\lambda+(x+1) \theta) \pi(-x-1) & =
\left( \mu \frac{\mu^{x-1}}{\prod_{j=1}^{x-1}(\lambda+j \theta)} +
(\lambda+(x+1) \theta) \frac{\mu^{x+1}}{\prod_{j=1}^{x+1}(\lambda+j \theta)}
\right)
\pi(0)\\
& =
\left((\lambda+x \theta)  \frac{\mu^{x}}{\prod_{j=1}^{x}(\lambda+j \theta)} +
 \mu \frac{\mu^{x}}{\prod_{j=1}^{x}(\lambda+j \theta)}
\right)
\pi(0)\\
& = ( \lambda + \mu+x \theta ) \pi(-x).
\end{split}
\end{equation*}
Last, we have that
\begin{equation*}
\begin{split}
(\lambda+\theta) \pi(-1)+ (\mu+\theta) \pi(1) & =
\left( (\lambda+\theta) \frac{\mu}{\lambda+ \theta} +
(\lambda+\theta) \frac{\lambda}{\mu+\theta}
\right)
\pi(0)\\
& =(\lambda+\mu)
\pi(0).
\end{split}
\end{equation*}
In other words, \eqref{eq:solutionMarkov} satisfies the balance equations.
Now, adapting \cite{ward03}, we derive the following formulas
\begin{equation*}
\prod_{j=1}^{x}(\mu+j \theta)= \frac{1}{\mu}  \prod_{j=0}^{x}(\mu+j \theta)=
\frac{\theta^{x+1} \Gamma(\frac{\mu}{\theta}+x+1)}{\mu \Gamma(\frac{\mu}{\theta})},
\end{equation*}
\begin{equation*}
\prod_{j=1}^{x}(\lambda+j \theta)= \frac{1}{\lambda}  \prod_{j=0}^{x}(\lambda+j \theta)=
\frac{\theta^{x+1} \Gamma(\frac{\lambda}{\theta}+x+1)}{\lambda \Gamma(\frac{\lambda}{\theta})},
\end{equation*}
\begin{equation*}
\sum_{x=1}^{\infty}\frac{\lambda^x}{\prod_{j=1}^{x}(\mu+j \theta)}=
\gamma\left( \frac{\lambda}{\theta}, \frac{\mu}{\theta}\right)
\left(\frac{\lambda}{\theta}\right)^{-\mu/ \theta} \frac{\mu}{\theta}
e^{\lambda/\theta}
-1,
\end{equation*}
\begin{equation*}
\sum_{x=1}^{\infty}\frac{\mu^x}{\prod_{j=1}^{x}(\lambda+j \theta)}=
\gamma\left( \frac{\mu}{\theta}, \frac{\lambda}{\theta}\right)
\left(\frac{\mu}{\theta}\right)^{-\lambda/ \theta}  \frac{\lambda}{\theta}
e^{\mu/\theta}
-1,
\end{equation*}
where $\gamma(x,y):=\int_{0}^{x} t^{y-1}e^{-t}dt$ and $\Gamma(y):=\gamma(\infty,y)$. The first two equations follow directly by \cite{ward03}. We show the third equation and the forth one follows by symmetry. By the power series expansion of the incomplete gamma function we have that
\begin{equation*}
\begin{split}
\gamma\left( \frac{\lambda}{\theta}, \frac{\mu}{\theta}\right) & =
\left(\frac{\lambda}{\theta}\right)^{\mu/ \theta}
e^{-\lambda/\theta}
\sum_{x=0}^{\infty}
\frac{(\lambda/\theta)^x}{\prod_{j=0}^{x}(\mu/\theta+j)}\\
& =
\left(\frac{\lambda}{\theta}\right)^{\mu/ \theta}
e^{-\lambda/\theta}
\sum_{x=0}^{\infty}
\frac{\lambda^x \theta}{\prod_{j=0}^{x}(\mu+ j \theta)}\\
& =
\left(\frac{\lambda}{\theta}\right)^{\mu/ \theta}
e^{-\lambda/\theta} \frac{\theta}{\mu}
\sum_{x=0}^{\infty}
\frac{\lambda^x }{\prod_{j=1}^{x}(\mu+ j \theta)}\\
& =
\left(\frac{\lambda}{\theta}\right)^{\mu/ \theta}
e^{-\lambda/\theta} \frac{\theta}{\mu}
\left ( \sum_{x=1}^{\infty}
\frac{\lambda^x}{\prod_{j=1}^{x}(\mu+ j \theta)} +1 \right),
\end{split}
\end{equation*}
where we define an empty product to be one.

Let $\pi^n(\cdot)$ denote the distribution of $X^n(\infty)$ with rates $n \lambda$ and $n \mu$. The means of $\frac{X^{+,n}(\infty)}{n}$ and $\frac{X^{-,n}(\infty)}{n}$ are given by the following expressions
\begin{equation*}
\begin{split}
\frac{1}{n}\E{X^{+,n}(\infty)} &= \frac{1}{n}
\sum_{x=-\infty}^{\infty} \max(x,0) \pi^n(x)
=
\frac{1}{n}\sum_{x=1}^{\infty}
x \frac{(n\lambda)^x}{\prod_{j=1}^{x}(n\mu+j \theta)} \pi^n(0)\\
&= \frac{ \mu}{\theta}
\Gamma\left(\frac{n \mu}{\theta}\right) \pi^n(0)
\sum_{x=1}^{\infty}
x \frac{(n\lambda/\theta)^x}{\Gamma\left(\frac{n\mu}{\theta}+x+1\right)}
\end{split}
\end{equation*}
and
\begin{equation*}
\begin{split}
\frac{1}{n}\E{X^{-,n}(\infty)} =
\frac{ \lambda}{\theta}
\Gamma\left(\frac{n \lambda}{\theta}\right) \pi^n(0)
\sum_{x=1}^{\infty}
x \frac{(n\mu/\theta)^x}{\Gamma\left(\frac{n\lambda}{\theta}+x+1\right)},
\end{split}
\end{equation*}
where
\begin{equation*}
\pi^n(0) =
\left(
\gamma\left( \frac{n\lambda}{\theta}, \frac{n\mu}{\theta}\right)
\left(\frac{n\lambda}{\theta}\right)^{-n\mu/ \theta} \frac{n\mu}{\theta}
e^{n\lambda/\theta}
+
\gamma\left( \frac{n\mu}{\theta}, \frac{n\lambda}{\theta}\right)
\left(\frac{n\mu}{\theta}\right)^{-n\lambda/ \theta}  \frac{n\lambda}{\theta}
e^{n\mu/\theta}
-1
 \right)^{-1}.
\end{equation*}
When $\mu/\theta:=a\in \mathbb{N}$ and
$\lambda/\theta:=b\in \mathbb{N}$, we have that
\begin{equation*}
\sum_{x=1}^{\infty}
x\frac{(n\mu/\theta)^x}{\Gamma\left(\frac{n\lambda}{\theta}+x+1\right)}
=
\frac{b^{-an}n^{-an+1}
\left(-ae^{bn}\gamma(an+1,bn)
+b^{bn}e^{bn}\gamma(an+1,bn)+ b^{an+1} n^{an}\right) }
{(a n)!}.
\end{equation*}
In the special case where $\lambda=\mu=\theta$,  by observing that $\Gamma\left(n\right)=(n-1)!$ and $\sum_{x=1}^{\infty}
x\frac{n^x}
{\Gamma\left(n+x+1\right)}=\frac{1}{(n-1)!}$, the above expressions can be simplified further as follows
\begin{equation*}
\begin{split}
\frac{1}{n}\E{X^{+,n}(\infty)} = \frac{1}{n}\E{X^{-,n}(\infty)}
=
(n-1)! \pi^n(0)
\sum_{x=1}^{\infty}
 \frac{x n^x}{(n+x)!}
=\pi^n(0),
\end{split}
\end{equation*}
where
$\pi^n(0) =
\left(
2 \gamma\left( n, n\right)
n^{-n+1}
e^{n}
-1
 \right)^{-1}$
and
$\gamma(x,n):=\int_{0}^{x} t^{n-1}e^{-t}dt=
(n-1)!\left(1- e^{-x} \sum_{k=0}^{n-1}\frac{x^k}{k!}
 \right)$.
Moreover, the mean queue-lengths are exactly
$\frac{1}{n}\E{Q^{D,n}(\infty)}=\frac{1}{n}\E{X^{+,n}(\infty)}$ and $\frac{1}{n}\E{Q^{S,n}(\infty)}=\frac{1}{n}\E{X^{-,n}(\infty)}.$

\end{document}